\documentclass[12pt]{amsart}
\usepackage{a4wide}
\usepackage{enumerate}
\usepackage{graphicx}
\usepackage{color}
\usepackage{amsmath}
\usepackage{comment}
\usepackage{soul}
\allowdisplaybreaks

\newtheorem{theorem}{Theorem}[section]
\newtheorem{corollary}[theorem]{Corollary}

\newtheorem{lemma}[theorem]{Lemma}
\newtheorem{proposition}[theorem]{Proposition}

\theoremstyle{definition}
\newtheorem{definition}[theorem]{Definition}
\newtheorem{example}[theorem]{Example}
\newtheorem{remark}[theorem]{Remark}

\newcommand{\ep}{\varepsilon}
\newcommand{\eps}[1]{{#1}_{\varepsilon}}

\let\pa\partial  
\let\del\partial  
\let\na\nabla

\newcommand{\R}{{\mathbb R}} 
\newcommand{\T}{{\mathbb T}}
\newcommand{\F}{{\mathcal F}}
\newcommand{\FrhoE}{\mathcal F^0}
\newcommand{\velo}{u}

\newcommand{\N}[1]{\left|#1\right|}
\newcommand{\NN}[1]{\left\|#1\right\|}
\newcommand{\NNN}[1]{\left|\!\left|\!\left|#1\right|\!\right|\!\right|}

\title[Semiconductor Boltzmann-Dirac-Benney EQ.]{Semiconductor Boltzmann-Dirac-Benney equation with a BGK-type collision operator:\\ existence of solutions vs.~ill-posedness}


\email{marcel.braukhoff@asc.tuwien.ac.at}

\thanks{The  author was partially funded by the Austrian Science Fund (FWF) project F 65.} 

\begin{document}
	\maketitle
	
	% Enter the first author's name and address:
	\centerline{\scshape Marcel Braukhoff}
	\medskip
	
	{\footnotesize
		% please put the address of the first author
		\centerline{Institute for Analysis and Scientific Computing}
		\centerline{Vienna University of Technology}
		\centerline{Wiedner Hauptstrasse 8-10, 1040 Wien, Austria}
	} % Do not forget to end the {\footnotesize by the sign }
	
	%\medskip
	%
	%\centerline{\scshape First-name2 last-name2 and First-name3
	%last-name3}
	%\medskip
	%{\footnotesize
	% % please put the address of the second  and third author
	% \centerline{ First line of the address of the second author}
	%   \centerline{Other lines}
	%   \centerline{Springfield, MO 65810, USA}
	%}
	
	\bigskip
	
	% The name of the associate editor will be entered by an editorial staff
	% "Communicated by the associate editor name" is not needed for special issue.
%	\centerline{(Communicated by the associate editor name)}

	%The abstract of your paper
	\begin{abstract}
		A semiconductor Boltzmann equation with a non-linear BGK-type collision operator
		is analyzed for a cloud of ultracold atoms in an optical lattice:
		\[
		\partial_t f + \nabla_p\epsilon(p)\cdot\nabla_x f - \nabla_x n_f\cdot\nabla_p f =   n_f(1- n_f)(\mathcal{F}_f-f), \quad x\in\mathbb{R}^d, p\in\mathbb{T}^d, t>0.
		\]
		This system contains an interaction potential $n_f(x,t):=\int_{\mathbb{T}^d}f(x,p,t)dp$ being significantly more singular than the Coulomb potential, which is used in the Vlasov-Poisson system. This causes major structural difficulties in the analysis. Furthermore, we call $\epsilon(p) = -\sum_{i=1}^d$  $\cos(2\pi p_i)$ the dispersion relation and $\mathcal{F}_f$ denotes the Fermi--Dirac equilibrium distribution, which depends non-linearly on $f$ in this context. 
		
		In a dilute plasma---without collisions (r.h.s$.=0$)---this system is closely related to the Vlasov--Dirac--Benney equation. It is shown for analytic initial data that the semiconductor Boltzmann equation possesses a local, analytic solution. Here, we exploit the techniques of Mouhout and Villani by using Gevrey-type norms which vary over time. In addition, it is proved that this equation is locally ill-posed in Sobolev spaces close to some Fermi--Dirac equilibrium distribution functions.
	\end{abstract}
	
	\section{Introduction}
	In the last decades, the theory of charge transport in semiconductors has become a thriving field in applied mathematics. Due to the complexity of semiconductors consisting of some $ 10^{23}$ atoms, there are several effective equations describing different phenomenological properties of semiconductors. Recently, the description of charge transport in semiconductors was extended by an experimental model \cite{SHR12}: a cloud of ultra\-cold atoms in an optical lattice. In this model, the ultra\-cold atoms stand for the charged electrons and the optical lattice describes the periodic potential of the crystal, formed by the ions of the semiconductor. Using the interference of optical laser beams, the atoms are trapped in an optical standing wave \cite{Blo05}. In contrast to a solid lattice, the geometry of an optical lattice as well as the strength of the potential can easily be changed during the experiment. Moreover, the time scale slows down to milliseconds while working with temperatures of a few nanokelvin. Therefore, this experimental model is particularly suited to understand the physical behavior of solid materials and of great interest. In addition, it may have the potential to accomplish quantum information processors \cite{Jak04} as well as very precise atomic clocks \cite{ADH15}.
	
	The main difference between a cloud of ultra\-cold atoms and a system of electrons is the interaction potential. Assuming that the atoms are uncharged, the interaction potential is significantly more singular than the Coulomb potential of the electrons causing major structural difficulties in the analysis.

	In this paper we investigate the ill-posedness of the following Boltzmann equation for the distribution function
	$f(x,p,t)$,
	\begin{equation}\label{1.be}
	\pa_t f + u\cdot\na_x f + \na_x V_f\cdot\na_p f = Q(f),
	\end{equation}
	where $x\in\R^d$ is the spatial variable, $p$ is the crystal momentum,
	defined on the $d$-dimensional torus $\T^d$ with unit measure, 
	and $t>0$ is the time. The velocity $u$ is defined by $u(p)=\na_p\epsilon(p)$
	with the energy $\epsilon(p)$,
	$V_f(x,t)$ is the lattice potential, and $Q(f)$ is the collision operator. 
	Compared to the standard semiconductor Boltzmann equation, there are two major differences.
	
	First, we assume that the dispersion relation, i.e.~the band energy, is given by 
	\begin{equation}\label{1.eps}
	\epsilon(p) = -2\epsilon_0\sum_{i=1}^d \cos(2\pi p_i), \quad p\in\T^d,
	\end{equation}
	where $\epsilon_0$ denotes the tunneling rate of a particle from one lattice site to a neighboring one \cite{RMR10}. This dispersion relation is typically used in semiconductor physics as for an approximation of the lowest band \cite{AsMe77}. In contrast to this, a parabolic band structure is given by $\epsilon(p)=\frac12|p|^2$ \cite{Jue09}, which also occurs in kinetic gas theory as the microscopic kinetic energy of free particles.
	
	Second, the potential $V_f$ is supposed to be proportional to the particle density $n_f=\int_{\T^d}fdp$ with 
	\begin{equation}\label{1.V_f}
	V_f(x,t) = Un_f(x,t)=U\int_{\T^d}f(x,p,t)dp, \quad x,\in \R^d,p\in\T^d,t>0.
	\end{equation}
	Here, $U\neq0$ describes the strength of the on-site interaction
	between spin-up and spin-down components \cite{SHR12}. However,  in semiconductor
	physics, the interaction potential is often given by the Coulomb potential $\Phi_f$ of the electric field which fulfills $\Delta \Phi_f= n_f$ \cite{Jue09}. Due to this Poisson equation, the Coulomb potential is more regular than the particle density $n_f$ in contrast to the potential $V_f$ defined in \eqref{1.V_f}. Therefore, we expect a more ``singular behavior'' of \eqref{1.be} compared to the standard semiconductor Boltzmann equation; see the discussion below.

	Similar to \cite{SHR12}, we use the following relaxation-time approximation  
	\begin{equation}\label{eq.Qf}
	Q(f) = \gamma n_f(1-\eta n_f)(\F_f-f)
	\end{equation}
	for the collision operator,
	where $1/\gamma>0$ denotes the relaxation time
	and   \[\F_f(x,p,t) = \big(\eta + \exp(-\lambda_0(x,t)-\lambda_1(x,t)\epsilon(p))\big)^{-1},
	\quad x\in\R^d,\ p\in\T^d,\ t>0\]
	is  the generalized Fermi-Dirac distribution function depending on $f$ through the Lagrange multipliers $(\lambda_0,\lambda_1)$: We define $\lambda_0$ and $\lambda_1$ by the mass and energy constraints
	\[
	\int_{\T^d}(\F_f-f)dp=0,\quad\int_{\T^d}(\F_f-f)\epsilon(p)dp=0.
	\]
	Note that $\eta=1$ leads to the original Fermi-Dirac distribution as in \cite{SHR12} and $\eta=0$ entails that  $\F_f$ equals the Maxwell-Boltzmann distribution.
	% Throughout this article, we write
	%We may consider $\F_f$ as a function of $(\lambda_0,\lambda_1)$ and write
	%$\F(\lambda_0,\lambda_1;p)=[\eta+\exp((-\lambda_0-\lambda_1\epsilon(p))]^{-1}$.
	
	Physically, $\lambda_1$ can be interpreted as the negative inverse (absolute) temperature,
	while $\lambda_0$ is related to the so-called chemical potential \cite{Jue09}.
	Since the dispersion relation is bounded, the equilibrium $\F_f$ is well-defined and integrable for all
	$\lambda_1\in\R$, which includes negative absolute temperatures. These negative absolute temperatures can actual be realized in experiments with 
	ultracold atoms \cite{RMR10}. Negative temperatures occur in
	equilibrated (quantum) systems that are characterized by an inverted population of
	energy states. The thermodynamical implications of negative temperatures 
	are discussed in \cite{Ram56}.

	So far, there are some results for this type of equation using $\epsilon(p)=\frac12|p|^2$ and that $Q(f)$ either vanishes or is quadratic in $f$:
	
	Combining this with the Vlasov equation yields the Vlasov-Dirac-Benney equation
	\begin{equation}\label{eq: intro: Vlasov3}
	\begin{aligned}
	\pa_t f(x,\velo,t)+ \velo\cdot\na_x f(x,\velo,t)-\na\rho_f(x,t)\cdot\na_\velo f(x,\velo,t) &=0
	\end{aligned}
	\end{equation}
	for $x\in\R^d,\velo\in\R^d$ and $t >0$. In spatial dimension one, this equation can be used to describe the density of fusion plasma in a strong magnetic field in direction of the field \cite{BaNo12}.
	%   Note the similarity to the left-hand side of the semiconductor Boltzmann equation from \eqref{eq: intro: Boltzmann algemein}.
	It can be derived as a limit of a scaled non-linear Schr\"odinger equation \cite{BaBe16}. 
	Comparing the Vlasov-Poisson equation to the Equation \eqref{eq: intro: Vlasov3}, we see that the interaction potential $\Phi$ is long ranged (i.e., the support is the whole space) in contrast to the delta distribution with $\mathrm{supp}(\delta_0)=\{0\}$. Therefore, we can understand \eqref{eq: intro: Vlasov3} as a version of the classical Vlasov-Poisson system with a short-ranged Dirac potential, which motivated the ``Dirac" in the name of the Vlasov-Dirac-Benney equation. The name Benney is due to its relation to the Benney equation in dimension one (for details see \cite{BaBe13}).
	
	However, the analysis of a Vlasov-Dirac-Benney equation is more delicate as in \cite{JaNo11} only local in time solvability was shown for analytic initial data in spatial dimension one. Moreover, it is shown in \cite{BaBe13} that this system is not locally weakly $(H^m-H^1)$ well-posed in the sense of Hadamard. In \cite{HaNg15} it is shown that the Vlasov-Dirac-Benney equation is ill-posed in $d=3$, requiring that the spatial domain is restricted to the $3$-dimensional torus $\T^3$. More precisely, they show that the flow of solutions does not belong to $C^\alpha(H^{s,m}(\R^3\times\T^3),L^2(\R^3\times\T^3))$ for any $s\geq0,\alpha\in(0,1]$ and $m\in\mathbb N_0$. Here, $H^{s,m}(\R^3\times\T^3)$ denotes the weighted Sobolev space of order $s$ with weight $(x,\velo)\mapsto \langle\velo\rangle^m:=(1+\N{\velo}^2)^{m/2}$. Even more precisely, they prove that there exist a stationary solution $\mu=\mu(\velo)$ of \eqref{eq: intro: Vlasov3} and a family of solutions $(f_\varepsilon)_{\varepsilon>0}$, times $t_\varepsilon= O(\varepsilon\N{\log\varepsilon})$ and $(x_0,\velo_0)\in\T^3\times\R^3$ such that
	
	\begin{equation*}
	\lim_{\varepsilon\to0}\frac{\NN{f_\varepsilon-\mu}_{L^2([0,t_\varepsilon]\times B_\varepsilon(x_0)\times B_\epsilon(\velo_0))}}{\NN{\langle\velo\rangle^m(f_\varepsilon|_{t=0}-\mu)}^\alpha_{H^s(\T^3_x\times\R^3_\velo)}}=\infty,
	\end{equation*}
	where $B_\varepsilon(x_0)$ denotes the ball with radius $\varepsilon$ centered at $x_0$.  In addition, \cite{HaNg15} covers also equation \eqref{eq: intro: Vlasov3} with a non vanishing r.h.s.: The authors consider
	\[\pa_t f+ \velo\cdot\na_x f-\na\rho_f(x,t)\cdot\na_\velo f= Q(f,f)\]
	for a bilinear operator $Q$.
	
	Moreover, the Vlasov-Dirac-Benney equation can also be derived by a quasi-neutral limit of the Vlasov-Poisson equation \cite{HaRo16}. Han-Kwan and Rousset are also able to provide uniform estimates on the solution of the scaled Vlasov-Poisson equation. By taking the quasi-neutral limit, they prove the existence of a unique local solution $f\in C([0,T], H^{2m-1,2r}(\R^3\times\T^3))$ of the Vlasov-Dirac-Benney equation. For this, they require that the initial data $f_0\in  H^{2m,2r}(\R^3\times\T^3)$ satisfies the Penrose stability condition 
	
	\begin{equation*}%\label{Penrose}
	\inf_{x\in\mathbb T^d}\inf_{(\gamma,\tau,\eta)\in(0,\infty)\times\mathbb R\times\mathbb R^d\setminus\{0\}}\bigg|1-\int_0^\infty e^{-(\gamma+i\tau)s}\frac{i\eta}{1+|\eta|^2}\cdot\mathcal (F_v\nabla_ vf)(x,\eta s)ds\bigg|>0,
	\end{equation*}
	where $\mathcal{F}_v$ denotes the Fourier Transform in $v$.
	
	\subsection*{Focus of this article}
	
	We introduce a concrete BGK-type collision operator (see Equation \eqref{eq.Qf}) arising from semiconductor physics \cite{SHR12}, which depends nonlinearly on $f$. 
	Since a Vlasov equation with collisions is in general called a semiconductor Boltzmann equation, we may call our system a semiconductor Boltzmann-Dirac-Benney equation with a BGK-type collision operator: 
	
	Let $\gamma>0$, $U\neq0$, we consider
	\begin{equation}\label{3.be}
	\pa_t f + u(p)\cdot\na_x f - U\na_x n_f\cdot\na_p f = \gamma n_f(1-\eta n_f)(\F_f-f)
	\end{equation}
	with $f(x,p,0)=f_0(x,p)$, where $\F_f(x,p,t) = \big(\eta + \exp(-\lambda_0(x,t)-\lambda_1(x,t)\epsilon(p))\big)^{-1}$, for $x\in\R^d,\ p\in\T^d$ and $t>0$. Here, $\lambda_0,\lambda_1$ shall be chosen in such a way that
	
	\begin{equation}\label{constraints}
	n_f(x,t)=n_{\F_f}(x,t)\quad\mbox{and}\quad E_f(x,t)=E_{\F_f}(x,t),
	\end{equation}
	where $n_f(x,t):=\int_{\T^d}f(x,p,t)dp$ and $E_f(x,t):=\int_{\T^d}\epsilon(p)f(x,p,t)dp$. Moreover, we have $u(p)=\nabla_p\epsilon(p)$ with 
	
	\begin{equation*}
	\epsilon(p) = -2\epsilon_0\sum_{i=1}^d \cos(2\pi p_i), \quad p\in\T^d,
	\end{equation*}
	for some $\epsilon_0>0$.

	In the first theorem, we prove the local existence of a solution for analytic initial data. It therefore extends the existence results of \cite{JaNo11} and \cite{HaNg15} to our setting.
	\begin{theorem}\label{thm.BGK.local.solution.bitorus}
		Let $\eta>0$, $\gamma\geq0,$ $U\neq0$ and $f_0:\T^d\times\T^d\to(0,\eta^{-1})$ be analytic. Then there exists a time $T>0$ such that \eqref{3.be} admits a unique analytic solution $f:\T^d\times\T^d\times[0,T)\to\R$ with $f(x,p,0)=f_0(x,p)$.
	\end{theorem}
	%J\"ungel proves in Proposition 4.6 in \cite{Jue09} that the zero-set of this BGK-type colission operator are given by the genrealized Fermi-Dirac distributions $F_\lambda(p) = (\eta + \exp(-\lambda_0-\lambda_1\epsilon(p)))^{-1}$. 
	Physically, the BGK-collision operator shall drive the system into an equilibrium given by the generalized Fermi-Dirac distribution and one would expect some nicer results than in \cite{HaNg15}. However, the following theorem tells us that this is not always the case since some Fermi-Dirac equilibria are unstable, leading to an ill-posedness result.
	\begin{theorem}\label{thm.illposed}
		Let $k\in\mathbb N$, $\theta>0$ and $\gamma>0$, $U\neq0$. There exist $\bar{\lambda}\in\R^2$ and a time $\tau>0$ and such that there exist solutions $f_\delta:\R^d_x\times\T^d_p\times[0,\tau]\to [1,\eta^{-1}]$ of \eqref{3.be} such that
		
		\begin{equation*}
		\lim_{\delta\to0}\frac{\NN{f_\delta(\cdot,\cdot,t)-F_{\bar{\lambda}}}_{L^1(B_\delta(x,p))}}{\NN{f_\delta(\cdot,\cdot,0)-F_{\bar{\lambda}}}^\theta_{W^{k,\infty}(\R^3_x\times\T^3_p)}}=\infty\quad\mbox{for all }x\in\R^d,p\in\T^d,t\in(0,\tau),
		\end{equation*}
		where $F_{\bar{\lambda}}(p):=1/(\eta+e^{-\bar{\lambda}_0-\bar{\lambda}_1\epsilon(p)})$ is a steady-state solution of \eqref{3.be}.
	\end{theorem} 
	\begin{remark}
		The theorem can easily be extended to all $\gamma\in\R$. A sufficient condition for the critical $\bar{\lambda}$ is given by
		\begin{equation}\label{eq.Bed.for.lambda,c2}
		1<U \bar\lambda_1\int_{\T^d}F_{\bar{\lambda}}(p)(1-\eta F_{\bar \lambda}(p))dp.
		\end{equation}
		It is still an open problem, whether this condition is necessary. However, a similar condition also appears in a different context of semiconductor physics for ultra cold atoms: In \cite{BrJu17}, a formal drift-diffusion limit of \eqref{3.be} was considered. The formal analysis indicates degeneracies of the limiting diffusion equation, whenever
		
		\begin{equation*}%\label{eq.Bed.for.lambda,c2}
		1=U \bar\lambda_1\int_{\T^d}F_{\bar{\lambda}}(p)(1-\eta F_{\bar \lambda}(p))dp.
		\end{equation*}
	\end{remark}
	Now we would also like to be able to treat the full space $\R^d$ in the space variable instead of the periodic case. In a realistic physical experiment, the most part of the particle cloud is localized at the origin meaning that the density distribution tends to zero as $|x|\to\infty$. These functions have to be treated with caution since the Fermi-Dirac distributions $\F_f$ are not analytic in $f=0$ as we can see in the following remark.
	
	\begin{remark}
		According to the definition of the BGK-collision operator, $\F_f$ is uniquely determined by the constraints from \eqref{constraints} and can be rewritten as a function $\FrhoE:U\subset\R^2\times\T^d\to[0,\eta^{-1}]$ with
		\[\FrhoE(n_f(x,t),E_f(x,t);p)=\F_f(x,p,t).\]
		For this function, one can compute that
		\begin{equation*}
		\partial_E^2\FrhoE(n,0;p)=\frac{1-2\eta n}{8\epsilon_0^4d^2n(1-\eta n)}(\epsilon(p)^2-2\epsilon_0^2d)
		\end{equation*}
		see \cite{Bra17} section 5.5. Thus, we can see that the second derivative has a singularity in $n=0$ (and in $n=\eta^{-1}$). In particular, there exist a $g=g(n,p)$ with $g(0,\cdot)\neq0$ such that
		\begin{equation*}
		\partial_n^i\partial_E^2\FrhoE(n,0;p)=\frac{g(n,p)}{n^{i+1}(1-\eta n)^{i+1}}.
		\end{equation*}
		Clearly, this implies that $\FrhoE$ is not analytic in $(n,E)=0$. Fortunately, we are only interested in the composition of $\FrhoE$ with $n_f$ and $E_f$. The idea is to assume enough regularity on $f$ such that $\F_f$ is analytic.
		
		This leads to a first version of the local existence theorem for the whole space:
		
		%	hypothesis for admissible  particle densities $n$ and the energy densities $E$: 
		%	we suppose there exist $C_0>0$ and $\nu\geq0$ such that 
		%	
		%	\begin{align}\label{hyp-nE}
		%	|\del_x^a n(x)|+|\del_x^a E(x)|\leq C_0 n(x)(1-\eta n(x))a!\nu^{-a}
		%	\end{align}
		%	for all $0\neq a\in\mathbb N_0^d$ and all $x\in\R^d$.
		
		%	 Therefore, we require additional hypothesis in contrast to Theorem \ref{thm.BGK.local.solution.bitorus}.
	\end{remark}
	
	\begin{theorem}\label{thm.BGK.local.solution0}	Let $\eta>0$, $\gamma\geq0,$ $U\geq0$ and $\lambda^0=(\lambda_0^0,\lambda_1^0):\R^d\to\R^2$ be analytic such that $\lambda_1^0\in L^\infty(\R^d)$ and let
		
		\begin{equation*}
		F_{\lambda^0}(x,p):=\frac{1}{\eta+e^{-\lambda_0^0(x)-\lambda_1^0(x)\epsilon(p)}}.
		\end{equation*}
		Moreover, we suppose there exist $C_0>0$ and $\nu\geq0$ such that 
		
		\begin{equation}\label{hyp-nE}
		|\del_x^a n_{F_{\lambda^0}}(x)|+|\del_x^a E_{F_{\lambda^0}}(x)|\leq C_0 n_{F_{\lambda^0}}(x)(1-\eta n_{F_{\lambda^0}}(x))a!\nu^{-|a|}
		\end{equation}
		for all $0\neq a\in\mathbb N_0^d$ and all $x\in\R^d$. Then there exists $T>0$ such that \eqref{3.be} admits a unique analytic solution $f:\R^d\times\T^d\times[0,T)\to\R$ with $f(x,p,0)=F_{\lambda^0}(x,p)$.
	\end{theorem}
	
	\begin{example}
		In this version of the local existence result, we allow also initial data which may approach zero as $|x|\to\infty$. Let $\lambda_1^0=0$ and
		
		\begin{align*}
		\lambda_0^0(x):=-\log (1+x^2).
		\end{align*}
		Then 
		
		\[F_{\lambda^0}(x)=\frac{1}{\eta+1+x^2}\]
		and hence $E_{F_{\lambda^0}}$ vanishes and $n_{F_{\lambda^0}}=F_{\lambda^0}(x)$. We will prove in example \ref{exappendix} in the appendix that
		
		\begin{align*}
		\N{F_{\lambda^0}^{(a)}(x)}
		&\leq 	\frac{a!}{\nu^a}F_{\lambda^0}(x)\qquad\mbox{for }\nu=\frac12\min\{\sqrt{\eta},1\}
		%	\leq \frac{\eta+1}{\eta} F_{\lambda^0}(x)(1-\eta F_{\lambda^0}(x))
		\end{align*}
		Using that $n_{F_{\lambda^0}}=F_{\lambda^0}(x)\leq1/(1+\eta)$ yields
		
		\begin{align*}
		\N{n_{F_{\lambda^0}}^{(a)}(x)}
		&
		\leq		
		%F_{\lambda^0}(x)
		\frac{\eta+1}{\eta} n_{F_{\lambda^0}}(x)(1-\eta n_{F_{\lambda^0}}(x))a!\nu^{-|a|}.
		\end{align*}	
		Finally, we can conclude that 
		
		\[F_{\lambda^0}(x)=\frac{1}{\eta+1+x^2}\]
		satisfies the hypothesis of the foregoing theorem. Thus, there exists $T>0$ such that \eqref{3.be} admits a unique analytic solution $f:\R^d\times\T^d\times[0,T)\to\R$ with $f(x,p,0)=F_{\lambda^0}(x,p)$.
	\end{example}
	%If we also want admit functions which may approach zero, which was not possible for the solutions on the torus given by Theorem \ref{thm.BGK.local.solution.bitorus}.  Therefore, we require additional hypothesis in contrast to Theorem \ref{thm.BGK.local.solution.bitorus}.
	
	%Nevertheless, we also prove a well-posedness for analytic norms.
	Note that \eqref{hyp-nE} is a local conditions for the particle and energy densities. This is a consequence of the fact that the BGK-collision operator is local in space. 
	%Therefore, it is reasonable to work with analytic semi-norms, which are local in space in order to treat the condition from \eqref{hyp-nE}.
	%%%%%%%%%%%%%%%%%%%
	%%%%%%%%%%%%%%%%%%%
	%%%%%%%%%%%%%%%%%%%
	%%%%%%%%%%%%%%%%%%%
	%%%%%%%%%%%%%%%%%%%
	%%%%%%%%%%%%%%%%%%%
	%%%%%%%%%%%%%%%%%%%
	%\begin{definition}
	%	For $d\in\mathbb N$ let $f:\R^d\times\T^d\to \R$ being analytic in $(x,p)$ for all $p\in \T^d$. Given $x\in\R^d,\nu\geq0$ we define the semi-norm
	%	
	%	\begin{equation*}
	%	|f|_{\dot{C}^\nu_x}:= \sum_{\substack{a,b\in\mathbb N_0^d\\a+b>0}} \frac{\nu^{|a+b|}}{a!b!}\int_{\T^d}|\del_x^a\del_p^b f(x,p)|dp.
	%%	\sum_{\substack{a,b\in\mathbb N_0\\a+b>0}} \frac{\nu^{a+b}}{a!b!}\int_{\T^d}\bigg(\sum_{\alpha\in\mathbb N_0^d,|\alpha|=a}\sum_{\beta\in\mathbb N_0^d|\beta|=b}\frac{a!b!}{\alpha!\beta!}|\del_x^\alpha\del_p^\beta f(x,p)|^2\bigg)^{\frac12}dp\bigg|_{x=y}.
	%	\end{equation*}
	%\end{definition}
	
	%%%%%%%%%%%%%%%%%%%
	%%%%%%%%%%%%%%%%%%%
	%%%%%%%%%%%%%%%%%%%
	%%%%%%%%%%%%%%%%%%%
	%%%%%%%%%%%%%%%%%%%
	%%%%%%%%%%%%%%%%%%%
	%%%%%%%%%%%%%%%%%%%
	%In the following theorem, we also admit functions which may approach zero, which was not possible for the solutions on the torus given by Theorem \ref{thm.BGK.local.solution.bitorus}. These functions have to be treated with caution since the Fermi-Dirac distributions $\F_f$ are not analytic in $f=0$. Therefore, we require additional hypothesis in contrast to Theorem \ref{thm.BGK.local.solution.bitorus}.
	
	\begin{theorem}\label{thm.BGK.local.solution}
		Let $\eta>0$, $\gamma\geq0,$ $U\geq0$ and let $\lambda^0=(\lambda_0^0,\lambda_1^0):\R^d\to\R^2$, $C_0$, $\nu$ be as in Theorem \ref{thm.BGK.local.solution0}.	Then there exist $\delta>0$ and $T>0$ such that \eqref{3.be} admits a unique analytic solution $f:\R^d\times\T^d\times[0,T)\to\R$ with
		\[f(x,p,0)=f_0(x,p):=\frac{1}{\eta+e^{-\lambda_0^0(x)-\lambda_1^0(x)\epsilon(p)}}+g_0(x,p),\]
		if $g_0:\R^d\times\T^d\to\R$ is analytic with $0\leq f_0(x,p)\leq \eta^{-1}$ and satisfies
		
		\begin{equation*}%\label{eq.g0.1}
		\int_{\T^d}|\del_x^a\del_p^bg_0(x,p)|dp\leq C_0 n_{f_0}(x)(1-\eta n_{f_0}(x))a!b!\nu^{-|a+b|}
		\end{equation*}\label{eq.g0.2}
		as well as
		
		\begin{equation*}
		|\del_x^an_{g_0}(x)|+|\del_x^aE_{g_0}(x)|\leq \delta n_{F_{\lambda^0}}(x)(1-\eta n_{F_{\lambda^0}}(x))a!\nu^{-|a|}
		\end{equation*}
		for all $x\in \R^d$ and $a,b\in\mathbb N_0^d$ with $a+b\neq0$. 
		
		Moreover, there exist $\tilde C,\tilde{\nu}>0$ and $\tilde{T}\in(0,T)$ such that
		
		\begin{equation*}
		\int_{\T^d}|\partial_x^a\partial_p^b f(x,p,t)|dp\leq \tilde Ca!b!\tilde\nu^{-|a+b|}
		\end{equation*}
		for all $x\in\R^d,\ t\in[0,\tilde T]$.
	\end{theorem}
	
	\begin{remark}
		The solution is well-posed in the following sense: There exist $\tilde{\nu}>0$ and $\tilde C>0$, $\tilde{T}\in(0,T)$ such that two solutions $f^1,f^2$ of \eqref{3.be} fulfill
		
		\begin{multline*}
		\int_{\T^d}|\partial_x^a\partial_p^b(f^2(x,p,t)-f^1(x,p,t))|dp\\\leq  \tilde Ca!b!\tilde\nu^{-|a+b|}\sum_{\substack{\alpha,\beta\in\mathbb N_0^d\\\alpha+\beta>0}} \frac{\nu^{|\alpha+\beta|}}{\alpha!\beta!}\int_{\T^d}|\del_x^a\del_p^b(g^2_0(x,p)-g^1_0(x,p)|dp
		\end{multline*}
		for all $x\in\R^d,\ t\in[0,\tilde T]$, where
		\[f^i_0(x,p):=f^i(x,p,0)\quad\mbox{and}\quad g^i_0(x,p):=f^i_0(x,p)-\frac{1}{\eta+e^{-\lambda_0^0(x)-\lambda_1^0(x)\epsilon(p)}}\]
		satisfy the same conditions as $ f_0$ and $g_0$ from Theorem \ref{thm.BGK.local.solution} for $i=1,2$.  
	\end{remark}
	
	\section{Analytic norms}\label{analytic-norms}
	Our strategy to solve \eqref{3.be} will be applying a fixed-point argument. Therefore, we require suitable functions spaces: we use the following analytic norms, which are similar to those from \cite{MoVi11}.
	\begin{definition}
		Let $\nu>0$, $d\in\mathbb N$. We define
		
		\begin{equation*}
		\NN{f}_{C^{\nu}}:= \sum_{a,b\in\mathbb N_0^d} \frac{\nu^{|a+b|}}{a!b!}\|\del_x^a\del_p^b f\|_{W^{1,\infty}_xW^{1,1}_p}
		\end{equation*}
		for $f:\R^d\times\T^d\to \R^k$ being analytic, where we use the notation
		
		\begin{equation*}
		\|f\|_{W^{1,\infty}_xW^{1,1}_p}:=\sum_{\substack{a,b\in\mathbb N_0^d\\|a+b|\leq 1}}\|\del_x^a\del_p^bf\|_{L^{\infty}_xL^{1}_p}\quad\mbox{and}\quad \|f\|_{L^{\infty}_xL^{1}_p}:=\sup_{x\in \R^d}\int_{\T^d}|f(x,p)|dp.
		\end{equation*}
		Moreover, we define the semi-norm
		
		\begin{equation*}
		\NN{Df}_{C^{\nu}}:= \sum_{\substack{a,b\in\mathbb N_0^d\\|a+b|= 1}}\NN{\del_x^a\del_p^af}_{C^{\nu}}
		\end{equation*}
		and we set
		\[
		\NN{\velo}_{C^{\nu,\infty}}:= \max_{i=1,\ldots,d}\sum_{b\in\mathbb N_0^d}\frac{\nu^{|b|}}{b!}\NN{\del^{b} \velo_i}_{W^{1,\infty}(\T^d)}.
		\]
		for $\velo=(\velo_1,\ldots,\velo_d):\T^d\to\R^d$.
	\end{definition}
	Comparing these norms to the analytic norms 
	
	\begin{align*}
	\N{f}_{C^{\nu}}&:= \sum_{a,b\in\mathbb N_0^d} \frac{\nu^{|a+b|}}{a!b!}\|\del_x^a\del_p^b f\|_{L^{\infty}_xL^{1}_p}
	\quad\mbox{and}\quad
	\N{\velo}_{C^{\nu,\infty}}:= \sum_{b\in\mathbb N_0^d}\frac{\nu^{|b|}}{b!}\NN{\del^{b} \velo}_{L^{\infty}(\T^d)}
	\end{align*} 
	from \cite{MoVi11}, we have the trivial estimate $\N{\cdot}_{C^{\nu}}\leq \NN{\cdot}_{C^{\nu}}$. For the inverse estimate, we can only compare $\N{\cdot}_{C^{\mu}}$ with $\NN{\cdot}_{C^{\nu}}$ if $\mu>\nu$ as the following lemma suggests. As we will see later on, the norm $\NN{\cdot}_{C^{\nu}}$ is suited better for treating semiconductor Boltzmann-Dirac-Benney type equations. The idea is to do the analysis with our tailor-made norms $\NN{\cdot}_{C^{\nu}}$. We only use the more ``standard" analytic norms  $\N{\cdot}_{C^{\mu}}$ afterward for the statements by using the following comparison estimate.
	
	\begin{lemma}\label{lem: nnorm leq C norm 1}
		Let $\mu>\nu>0$ and $d\in\mathbb N$. Then there exists $C_{\mu,\nu}>0$ such that 
		
		\[
		\NN{f}_{C^{\nu}}\leq C_{\mu,\nu}\N{f}_{C^{\mu}}
		\]
		for all analytic $f:\R^d\times\T^d\to \R$. %$=1+\frac{2}{\nu}\sup_{a\geq0}a\big(\frac{\nu}{\mu}\big)^a=1+2/(\nu \operatorname e \log(\mu/\nu))$.
	\end{lemma}
	\begin{proof}
		It suffices to show that we have $\NN{\del f}_{C^{\nu}}\leq C\N{f}_{C^{\mu}_x}$ for $\del\in\{\del_x,\del_p\}$ for some $C>0$. Let $\del=\del_x$ and compute
		
		\begin{align*}
		\NN{\del_x f}_{C^{\nu}}
		&=
		\sum_{i,j\in\mathbb N_0 }\frac{\nu^{i+j}}{i!j!}\|\del_x^{i+1}\del_p^jf\|_{L^\infty_xL^1_p}
		=\frac{1}{\nu}
		\sum_{i,j\in\mathbb N_0 }i\frac{\nu^{i+j}}{i!j!}\int_{\T^d}\|\del_x^{i}\del_p^jf\|_{L^\infty_xL^1_p}\\
		&\leq\frac{1}{\nu}\sup_{a\in\mathbb N}a\frac{\nu^a}{\mu^a}
		\sum_{i,j\in\mathbb N_0 }\frac{\mu^{i+j}}{i!j!}\int_{\T^d}\|\del_x^{i}\del_p^jf\|_{L^\infty_xL^1_p}=	C\N{\del_x f}_{C^{\mu}_x}	
		\end{align*}
		for $C= \sup_{a\in\mathbb N}a\frac{\nu^{a-1}}{\mu^a}<\infty$. The estimate for $\del=\del_p$ can be proved similarly.
	\end{proof}
	The equation \eqref{2.be} consists of terms which involve product. Therefore, the following algebraic properties are particularly useful for treating equation \eqref{3.be}.
	%\begin{remark}
	%	Let $f:\R^d\times \T^d\to \R$ and $n:\R^d\to\R$ be both analytic, where $n$ is independent from $p$. Then we have that $\N{fn}_{C^{\nu}}\leq \N{f}_{C^{\nu}}\N{n}_{C^{\nu}}$ for $\nu\geq0$ and $x\in \R^d$, which can be proved similarly as in  \cite{MoVi11}.
	%\end{remark}
	%For the analysis of the semiconductor Boltzmann equation for ultracold atoms, we require slightly different analytic norms.
	\begin{lemma}\label{lem: analytic: del1 f mal  del2 g Abschaetzung}
		Let $f:\R^d\times\T^d\to\R$, $n:\R^d\to\R$ and $\velo:\T^d\to\R^d$ be analytic.	Let $\nu\geq0$. Then it holds
		
		\begin{align*}
		\NN{fn}_{C^{\nu}}\leq \NN{f}_{C^{\nu}} \NN{n}_{C^{\nu}}\quad\mbox{and}\quad\NN{\velo\cdot\nabla_x f}_{C^{\nu}} \leq\NN{\velo}_{C^{\nu,\infty}}\NN{Df}_{C^{\nu}}
		\intertext{as well as}
		\NN{\nabla_xn\cdot\nabla_p f}_{C^{\nu}} \leq\NN{n}_{C^{\nu}}\NN{Df}_{C^{\nu}} +\NN{Dn}_{C^{\nu}}\NN{f}_{C^{\nu}}.
		\end{align*}
	\end{lemma}
	\begin{proof}
		First, we try to rewrite the norm $\NN{\cdot}_{C^{\nu}}$ in such a way that we can use the results of \cite[section 4]{MoVi11}. Then we can easily show using the Leibniz rule that  $\N{fn}_{C^{\nu}}\leq \N{f}_{C^{\nu}}\N{n}_{C^{\nu}}$  and $\N{ f\velo}_{C^{\nu}}\leq \N{f}_{C^{\nu}}\N{\velo}_{C^{\nu,\infty}}$ (see \cite[section 4]{MoVi11}).
		Using this and the chain rule, we have 
		
		\begin{align*}
		\NN{fn}_{C^{\nu}}=\!\!\sum_{|a+b|\leq 1}\N{\del_x^a\del_x^b(fn)}_{C^{\nu}}\leq \!\!\sum_{|a_i+b_i|\leq 1}\!\!\N{\del_x^{a_1}\del_x^{b_1}f}_{C^{\nu}}\N{\del_x^{a_2}\del_x^{b_2}n}_{C^{\nu}}=\NN{f}_{C^{\nu}}\NN{n}_{C^{\nu}}.
		\end{align*}
		Likewise,
		
		\begin{align*}
		\NN{\velo\cdot\nabla f}_{C^{\nu}} 
		&
		\leq \sum_{i=1}^d\sum_{a,b\in \mathbb N_0^d,|a+b|\leq 1}\big(\N{\velo_i}_{C^{\nu,\infty}}\N{\del^{a}_x\del^b_p \del_{x_i}f}_{C^{\nu}} +\NN{\del_p^{b}\velo_i}_{C^{\nu,\infty}}\N{\del^a_x \del_{x_i}f}_{C^{\nu}}\big)
		\\&\leq\NN{\velo}_{C^{\nu,\infty}}\sum_{i=1}^d\NN{\del_{x_i}f}_{C^{\nu}}
		\intertext{and}
		\NN{\nabla_xn\cdot\nabla_p f}_{C^{\nu}} &\leq 
		\sum_{i=1}^d\sum_{a,b\in \mathbb N_0^d,|a+b|\leq 1}\big(\NN{\del_{x_i}n}_{C^{\nu}}\N{\del^{a}_x\del^b_p \del_{p_i}f}_{C^{\nu}} 
		\\&\hspace{5cm}+\N{\del_x^{a}\del_{x_i}n}_{C^{\nu}}\N{\del^b_p \del_{p_i}f}_{C^{\nu}}\big)
		\\
		&\leq\sum_{i=1}^d\NN{n}_{C^{\nu}}\NN{\del_{p_i}f}_{C^{\nu}} +\sum_{i=1}^d\NN{\del_{x_i}n}_{C^{\nu}}\NN{f}_{C^{\nu}}.\qedhere
		\end{align*}
	\end{proof}

	In \cite{MoVi11}, Mouhot and Villani unleashed the full potential of these analytic norms by varying the index $\nu$ over time. Motivated by their results, we define the following norm and derive the proceeding lemma.
	\begin{definition}\label{Def Norm lambda mu T}
		For $\nu,T>0$, $\mu\in[0,\nu/T)$, we define 
		
		\begin{equation*}
		\NN{f}_{\nu,\mu}:=\sup_{0\leq t<T}\left(\NN{f(t)}_{C^{\nu-\mu t}}+\mu\int_0^t\NN{Df(s)}_{C^{\nu-\mu s}}ds\right)
		\end{equation*}
		for $f:\R^d\times\T^d\times[0,T)\to\R$ being analytic in $(x,p)$ and continuous in $t$ writing $f(t)=f(\cdot,\cdot,t)$.% We moreover put $\NN{f}_{\lambda_0,\mu,T}=\sup_x \NN{f}_{\lambda_0,\mu,T,x}$.
	\end{definition}
	
	\begin{lemma}\label{lem: d_t nnormClambda  f estimate}
		For $\nu,T>0$, $\mu\in[0,\nu/T)$ and $f:\R^d\times\T^d\times(0,T)\to \R$  be analytic in $x,p$ and  continuously differentiable in $t$. Then
		
		\begin{equation*}
		%	\NN{f(t)}_{C^{\nu-\mu t}}+\mu\int_s^t\NN{Df(\tau)}_{C^{\nu-\mu\tau}}d\tau
		\NN{f}_{\nu,\mu}
		\leq \NN{f(0)}_{C^{\nu}}+\int_0^T\NN{\del_tf(\cdot,\cdot,t)}_{C^{\nu-\mu t}}d t.
		\end{equation*}
	\end{lemma}
	\begin{proof}
		Let $0<t\leq T$. Throughout this proof, we write $f(t):=f(\cdot,\cdot,t)$. Without loss of generality, we assume that
		
		\begin{equation*}
		\tau\mapsto \NN{\del_tf(\tau)}_{C^{\nu-\mu \tau}}\in L^1(0,t),
		\end{equation*}
		because otherwise, the assertion is trivial. Setting
		\[P_{f,N}(\lambda,t):=\sum_{|a|,|b|\leq N}\frac{\lambda^{|a+b|}}{a!b!}\NN{\del_x^{i+a}\del_p^{j+b}f(t)}_{W^{1,\infty}_xW^{1,1}_p}\quad\ \ \]
		and
		\[Q_N(\lambda,t):=\sum_{|i+j|=1}\sum_{\substack{|a|,|b|\leq N\\|a+b|<2N}}\frac{\lambda^{|a+b|}}{a!b!}\NN{\del_x^{i+a}\del_p^{j+b}f(t)}_{W^{1,\infty}_xW^{1,1}_p},\]
		we have $P_{f,N}(\lambda,t)\to \NN{f(t)}_{C^\lambda}$ and $Q_N(\lambda,t)\to \NN{Df(t)}_{C^\lambda}$ as $N\to \infty$.  Let $i,j,a,b\in\mathbb N_0^d$ and $0<s<t$. Then
		
		\begin{align*}
		&\N{\NN{\partial_x^{i+a}\partial_p^{j+b}f(t)}_{W^{1,\infty}_xW^{1,1}_p}-\NN{\del_x^{i+a}\del_p^{j+b}f(x,\cdot,s)}_{W^{1,\infty}_xW^{1,1}_p}}
		\\&\hspace{2cm}\leq \NN{\del_x^{i+a}\del_p^{j+b}f(x,\cdot,t)-\del_x^{i+a}\del_p^{j+b}f(x,\cdot,s)}_{W^{1,\infty}_xW^{1,1}_p}
		\\&\hspace{2cm}\leq \sup_{s\leq \tau\leq t}\NN{\del_x^{i+a}\del_p^{j+b}\del_tf(x,\cdot,\tau)}_{W^{1,\infty}_xW^{1,1}_p}(t-s)
		\end{align*}
		implies
		\[\N{P_{f,N}(\lambda,t)-P_{f,N}(\lambda,s)}\leq  \sup_{s\leq \tau\leq t}P_{\del_t f,N}(\lambda,\tau)(t-s).\]
		Next, let $\lambda=\lambda(t)=\nu-\mu t$. Using the  estimate
		
		\begin{align*}
		\frac{\lambda(t)^a-\lambda(s)^a}{a!}&=
		\frac{(\lambda(s)+\mu(s-t))^a-\lambda(s)^a}{a!}
		\\	&=\mu(s-t)\sum_{j=0}^{a-1}\frac{\lambda(s)^j\lambda^{a-1-j}}{j!(a-1-j)!(a-j)}
		\begin{cases}
		\leq \mu(s-t)\frac{\lambda(s)^{a-1}}{(a-1)!}\\
		\geq \mu(s-t)\frac{\lambda(t)^{a-1}}{(a-1)!},
		\end{cases}
		\end{align*}
		we can derive that 
		
		\begin{align*}
		\N{P_{f,N}(\nu-\mu t,t)-P_{f,N}(\nu-\mu s,s)}&\leq  \sup_{s\leq \tau\leq t}P_{\del_t f,N}(\nu-\mu t,\tau)(t-s) \\&\hspace{2cm}+ \mu\sup_{s\leq \tau\leq t}Q_N(\nu-\mu\tau,s)(t-s).
		\end{align*}
		Thus, $P_{f,N}(\nu-\mu t,t)$ is Lipschitz continuous w.r.t.~$t$ and belongs to $W^{1,\infty}((0,T))$ with
		\[\frac{d}{dt}P_{f,N}(\nu-\mu t,t)\leq P_{\del_t f,N}(\nu-\mu t,t)-\mu Q_N(\nu-\mu t,t),\]
		since $P_{f,N}$, $P_{\del_tf,N}$ and $Q_N$ are continuous.
		
		Since $P_{\del_t f,N}(\nu-\mu\tau,\tau)\leq\NN{\del_tf(\tau)}_{C^{\nu-\mu\tau}_x}\in L^1(0,T)$, the dominated convergence theorem implies that \[\int_0^TP_{\del_t f,N}(\nu-\mu\tau,\tau)d\tau\to \int_0^T\NN{\del_tf(\tau)}_{C^{\nu-\mu\tau}}d\tau\quad\mbox{as}\quad N\to\infty.\]
		Moreover, we can utilize the monotone convergence theorem in order to obtain that $\int_0^TQ_N(\nu-\mu\tau,\tau)d\tau \to \int_0^T\NN{Df(\tau)}_{C^{\nu-\mu\tau}}d\tau$. Thus, we summarize
		
		\begin{align*}
		&\NN{f(t)}_{C^{\nu-\mu t}}+\int_0^t\left(\mu\NN{Df(\tau)}_{C^{\nu-\mu\tau}}-\NN{\del_tf(\tau)}_{C^{\nu-\mu\tau}}\right)d\tau
		\\&\hspace{2cm} \stackrel{N\to\infty}\leftarrow
		P_{f,N}(\nu(t),t)+\int_0^t\left(\mu Q_N(\nu-\mu\tau,\tau)-P_{\del_t f,N}(\nu-\mu\tau,\tau)\right) d\tau
		\\&\hspace{2cm}\ \ \leq 
		P_{f,N}(\nu,0)
		\leq \NN{f(0)}_{C^{\nu}}
		\end{align*}
		finishing the proof.
	\end{proof}

	\begin{remark}\label{lem: Norm f_0<R}
		Let $T,\nu>0$ and let $f_0:\R^d\times\T^d\to\R$ be analytic. Then  
		
		\begin{equation*}
		\NN{f_0}_{\nu,\mu}
		%=\NN{f_0}_{C^{\nu}_x}+\mu\int_0^T\NN{Df_0}_{C^{\nu-\mu t}_x}dt
		=\NN{f_0}_{C^{\nu}} \quad\mbox{for every }\mu\in(0,\nu/T).
		\end{equation*}
		%	This is a direct consequence of Lemma \ref{lem: d_t nnormClambda  f estimate} since
		%	\[
		%	\NN{f_0}_{C^{\nu-\mu t}}+\mu\int_0^t\NN{Df_0}_{C^{\nu-\mu s}}ds\leq \NN{f_0}_{C^{\nu}}.
		%	\]
	\end{remark}
	\section{Local well-posedness in analytic norms}\label{local}
	In this section, we analyze the semiconductor Boltzmann equation \eqref{1.be} for ultracold atoms (setting $V_f:=-Un_f$ for $U\in \R$) in combination with a relaxation time approximation with fixed equilibrium. We consider
	
	\begin{equation}\label{2.be}
	\pa_t f + u(p)\cdot\na_x f - U\na_x n_f\cdot\na_p f = \gamma n_f(1-\eta n_f)(F-f),
	\end{equation}
	with $f(x,p,0)= f_0(x,p)$ for some given $F=F(x,p,t)$ and $\gamma\geq0$.

	\begin{theorem}\label{thm: boltzmann local: hT solution}
		Let $C,R,\nu>0$ and $f_0:\R^d\times\T^d\to \R$ and  $F:\R^d\times\T^d\times[0,T')\to \R$ be analytic such that
		\begin{equation}\label{eins}
		\NN{f_0}_{C^{\nu}}<R.
		\end{equation}
		Then if $\mu>0$ is sufficiently large, $T\in (0,\nu/\mu)$ and $F:\R^d\times\T^d\times[0,T')\to \R$ is analytic such that
		\begin{equation}\label{zwei}
		\NN{F(t)}_{C^{\nu-\mu t}}\leq C
		\end{equation}
		for all $0\leq t\leq T$, then equation \eqref{2.be}  admits a unique analytic solution $f:\R^d\times\T^d\times[0,T)\to\R$ with $\NN{f}_{\nu,\mu}\leq R$ and $f(x,p,0)= f_0(x,p)$ for all $x\in\R^d$ and $p\in \T^d$. 
		
		Moreover, let $\Psi: (f_0,F)\mapsto f$ be defined by the unique solution of
		
		\begin{equation*}
		\del_t f+\velo\cdot\nabla_x f - U\nabla_xn_f\cdot\nabla_p f =\gamma n_f(1-\eta n_f)(F-f)
		\end{equation*}
		with $f(x,p,0)=f_0(x,p)$. If $\mu>0$ is sufficiently large, the mapping $\Psi$ is Lipschitz continuous, i.e.,
		
		\begin{equation*}
		\NN{\Psi(f_0,F)-\Psi(g_0,G)}_{\nu,\mu}\leq 2\NN{(f_0,F)-(g_0,G)},
		\end{equation*}
		where
		\begin{align*}
		\NN{(f_0,F)}&:= \NN{f_0}_{C^\nu}+\mu^{-\frac12}\sup_{0\leq t<T}\NN{F}_{C^{\nu-\mu t}}.
		\end{align*} 
		for $f_0,g_0$ and $F,G$ satisfying \eqref{eins} and \eqref{zwei}, respectively.
	\end{theorem}
	\begin{remark}
		A sufficient condition for $\mu$ is given by \[\mu\geq \frac{C}{R-\NN{f_0}_{C^{\nu}}}+C\sup_{0\leq t\leq T}\NN{F}_{C^{\nu-\mu t}}\] for some $C>0$ independent from $f_0$ and $F$. 
	\end{remark}
	The key idea for the proof relies on the contraction mapping principle/Banach's fixed-point theorem. We define the mapping $\Phi$
	
	\begin{equation}\label{eq: VDhom ht: iterative}
	\Phi(f):=f_0-\int_0^t\big(\velo\cdot\nabla_x f - U\nabla_xn_f\cdot\nabla_p f -\gamma n_f(1-\eta n_f)(F-f)\big)dt.
	%\partial_t f_{j+1} + \velo\cdot\nabla_x f_j - U\nabla_xn_j\cdot\nabla_p f_j =\gamma n_j(1-\eta n_j)(F-f_j),
	\end{equation}
	for $f$ being analytic in $(x,p)$ and continuous in time. 
	%with $f_{j+1}(x,p,0)= f_0(x,p)$. 
	In order to prove that $\Phi$ admits a (unique) fixed-point,
	% using Banach's fixed-point theorem,
	we require the next lemmas.
	\begin{lemma}\label{lem.proof.thm.boltzmann local: hT solution1}
		Let $\NN{f_0}_{C^{\nu}}\leq R$. For every sufficient large $\mu$, there exists a $T\in(0,\nu/\mu)$ such that $\NN{f}_{\nu,\mu}\leq R$ implies $\NN{\Phi(f)}_{\nu,\mu}\leq R$. Here, a sufficient condition for $\mu$ is given by $\mu\geq C/(R-\NN{f_0}_{C^{\nu}})$ for some $C>0$ independent from $f_0$.
	\end{lemma}
	\begin{proof}
		First we fix $ \mu>0$ and $T\in(0,\nu/\mu)$. According to Lemma \ref{lem: d_t nnormClambda  f estimate}, we have
		
		\begin{align*}
		&\|\Phi(f)\|_{\nu,\mu}-\NN{f_0}_{C^{\nu}}\leq  \int_0^T\NN{\del_t \Phi(f)}_{C^{\nu-\mu t}}dt
		\\&\qquad\leq \int_0^T\big(\NN{\velo\cdot\na_x f}_{C^\lambda}+\N{U}\NN{\na_xn_f\cdot\na_p f}_{C^\lambda}+\gamma\NN{ n_f(1-\eta n_f)(f-F)}_{C^\lambda}\big)dt.
		\end{align*}
		Using the submultiplicativity obtained by Lemma \ref{lem: analytic: del1 f mal  del2 g Abschaetzung} and $\NN{n_f}_{C^\lambda}\leq \NN{f}_{C^\lambda}$ as well as $\NN{Dn_f}_{C^\lambda}\leq \NN{Df}_{C^\lambda}$, we obtain
		
		\begin{multline}\label{eq.inproof.73}
		\int_0^T\NN{\del_t\Phi(f)}_{C^{\nu-\mu t}}d\tau
		\leq 
		\frac{1}{\mu}\left(\NN{\velo}_{C^{\nu,\infty}}+\N{U}\NN{f}_{\nu,\mu}\right)\NN{f}_{\nu,\mu}
		\\ +
		\gamma T\NN{f}_{\nu,\mu}\left(1+\eta\NN{f}_{\nu,\mu}\right)\left(\NN{f}_{\nu,\mu}+\sup_{0\leq t\leq T}\NN{F}_{C^{\nu-\mu t}}\right).
		\end{multline}
		Thus, assuming $\NN{f}_{\nu,\mu}\leq R$ entails that
		
		\begin{align*}
		\NN{\Phi(f)}_{\nu,\mu}
		&\leq \NN{f_0}_{C^\nu}+\frac{1}{\mu}\left(\NN{\velo}_{C^{\nu,\infty}}+\N{U}R\right)R
		\\&\hspace{2cm}+
		T\gamma R\left(1+\eta R\right)\left(R+\sup_{0\leq t\leq T}\NN{F}_{C^{\nu-\mu t}}\right)
		\end{align*}
		Let $\mu_R:=2\frac{(\NN{\velo}_{C^{\nu,\infty}}+\N{U}R)R}{R-\NN{f_0}_{C^\nu}}>0 $. Then for all $\mu\geq \mu_R$, we have
		
		\begin{align*}
		\NN{\Phi(f)}_{\nu,\mu}
		&\leq \NN{f_0}_{C^\nu}+\frac{1}{2}(R-\NN{f_0}_{C^\nu})
		\\&\hspace{2cm}+
		T\gamma R\left(1+\eta R\right)\left(R+\sup_{0\leq t\leq T}\NN{F}_{C^{\nu-\mu t}}\right)
		\leq R
		\end{align*}
		if $T\in(0,\nu/\mu)$ satisfies
		
		\begin{align}\label{inproof_muT_Bedingung2}
		T\gamma R\left(1+\eta R\right)\left(R+\sup_{0\leq t\leq T}\NN{F}_{C^{\nu-\mu t}}\right)\leq \frac12\left(R-\NN{f_0}_{C^\nu}\right).
		\end{align}
		Therefore, for every sufficient large $\mu$, i.e.~$\mu\geq\mu_R$, every $T<\nu/\mu$ satisfies condition \eqref{inproof_muT_Bedingung2}. Thus, $\NN{\Phi(f)}_{\nu,\mu}\leq R$.
	\end{proof}
	\begin{lemma}\label{lem: Norm f_j+1-f_j<f_j-1-f_j}
		Let $R>\NN{f_0}_{C^{\nu}}$. For every sufficient large $\mu$, there exists a $T\in(0,\nu/\mu)$ such that
		$\NN{f_1}_{\nu,\mu},\NN{f_2}_{\nu,\mu}\leq R$ imply 
		\[\NN{\Phi(f_1)-\Phi(f_2)}_{\nu,\mu}\leq \frac12 \NN{f_1-f_2}_{\nu,\mu}.\]
		Here, a sufficient condition for $\mu$ is given by \[\mu\geq \frac{C}{R-\NN{f_0}_{C^{\nu}}}+C\sup_{0\leq t\leq T}\NN{F}_{C^{\nu-\mu t}}\] for some $C>0$ independent from $f_0$ and $F$.  
	\end{lemma}
	\begin{proof}
		The difference $g:=\Phi(f_2)-\Phi(f_1)$ is given by
		
		\begin{equation*}
		g(t)=\int_0^t\big(-\velo\cdot\na_xg+U\na_xn_{f_1}\cdot\na_p g+U\na_xn_g\cdot\na_p f_{2}
		+Q(f_2)-Q(f_{1})\big)ds
		\end{equation*}
		where $n_{f_j}=\int_{\T^d} f_j dp$, $n_g= \int_{\T^d} g dp$ and 
		
		\begin{equation*}
		Q(f_j):%=-\gamma\rho_j(1-\eta\rho_j)(f_j-\fnEhT(\rho_j,\E_j))
		=\gamma n_j(1-\eta n_j)\left(F-f_j\right)\quad\mbox{for }j=1,2.
		\end{equation*}
		Since $Q(f_j)$ is affine in $f_j$ and quadratic in $n_{f_j}$, we use the submultiplicativity properties of the norm $\NN{\cdot}_{C^{\nu-\mu t}}$ from Lemma \ref{lem: analytic: del1 f mal  del2 g Abschaetzung} to ensure that
		
		\begin{align*}
		%	     \NN\del_t g_{i+1}}_{C^{\nu-\mu t}_x}
		\NN{ Q(f_2)- Q(f_{1})}_{C^{\nu-\mu t}}
		&\leq 
		\gamma (2+3\eta R)\left(\NN{F}_{C^{\nu-\mu t}}+ R\right)\NN{g}_{C^{\nu-\mu t}}
		\end{align*}
		for  $\NN{f_{1}}_{C^{\nu-\mu t}},\NN{f_{2}}_{C^{\nu-\mu t}}\leq R$. We derive similarly to the proof of Lemma \ref{lem.proof.thm.boltzmann local: hT solution1} that
		
		\begin{multline*}
		\NN{\del_t (\Phi(f_2)-\Phi(f_1))}_{C^{\nu-\mu t}}
		\\\leq \big(\NN{\velo}_{C^{\nu,\infty}}+\N{U}\big(\NN{f_1}_{C^{\nu-\mu t}}+\NN{f_2}_{C^{\nu-\mu t}_x}\big)\big)\NN{Dg}_{C^{\nu-\mu t}}
		\\+
		\big(\N{U}\NN{Df_2}_{C^{\nu-\mu t}}+\N{U}\NN{Df_{1}}_{C^{\nu-\mu t}}\\+\gamma (2+3\eta R)\big(R+\NN{F}_{C^{\nu-\mu t}}\big)\big)\NN{g}_{C^{\nu-\mu t}}.
		\end{multline*}
		By Lemma \ref{lem: d_t nnormClambda  f estimate}, we obtain for all $\NN{f_j}_{\nu,\mu}\leq R$
		
		\begin{align*}
		\NN{\Phi(f_2)-\Phi(f_1)}_{\nu,\mu}&\leq \int_0^T  \NN{\del_t (\Phi(f_2)-\Phi(f_1))}_{C^{\nu-\mu t}_x} dt
		\\&\leq \frac1\mu\left(\NN{\velo}_{C^{\nu,\infty}}+4\N{U}R\right)\NN{g}_{\nu,\mu}
		\\&\hspace{1cm}	
		+T\gamma (2+3\eta R)\left(R+\sup_{0\leq t\leq T}\NN{F}_{C^{\nu-\mu t}}\right)\NN{g}_{\nu,\mu}\\&\leq\frac {C_{\nu,R}}\mu\NN{g}_{\nu,\mu}=\frac {C_{\nu,R}}\mu\NN{f_2-f_1}_{\nu,\mu},%<\frac12\NN{g_j}_{\lambda_0,\mu,T}
		\end{align*}
		where $C_{\nu,R}=\NN{\velo}_{C^{\nu,\infty}}+4\N{U}R+\nu\gamma (2+3\eta R)(R+\sup_{0\leq t\leq T}\NN{F}_{C^{\nu-\mu t}})$ using $T<\nu/\mu$. Finally, we obtain the assertion by assuming that $\mu\geq 2C_{\nu,R}$.
	\end{proof}
	\begin{proof}[Proof of Theorem \ref{thm: boltzmann local: hT solution}]
		Let $X$ consist of all functions $f:\R^d\times\T^d\times[0,T)\to \R$ being analytic in $x,p$ and continuous in $t$ such that $\|f\|_{\nu,\mu}\leq R$. Combining the previous two lemmata, we directly obtain that $\Phi:X\to X$ defined by \eqref{eq: VDhom ht: iterative} is a contraction requiring that $\mu$ is sufficiently large and $T\in(0,\nu/\mu)$. Thus, Banach's fixed-point theorem implies that equation \eqref{2.be} admits a unique mild solution in the space $X$. Using a bootstrap argument yields that $f$ is also analytic in $t$ and satisfies equation \eqref{2.be} classically.
		
		For the second part of the assertion, let 
		$f=\Psi(f_0,F)$, $g=\Psi(g_0,G)$.  There exists a $\tilde{\mu}>0$ such that for $T\in[0,\nu/\tilde{\mu})$ the functions $f,g$ are both defined on $[0,T)$ and satisfy $\NN{f}_{\tilde{\mu},\nu},\NN{g}_{\tilde{\mu},\nu}\leq R$. Defining $h=f-g, h_0=f_0-g_0$ as well as $H=F-G$, we have	
		
		\begin{equation*}
		\del_t h+\velo\cdot\nabla_x h - U\nabla_xn_h\cdot\nabla_p f- U\nabla_xn_g\cdot\nabla_p h =Q(f,F)-Q(g,G)
		\end{equation*}
		with $h(x,p,0)=h_0(x,p)$, where $Q(f,F)= \gamma n_f(1-\eta n_f)(F-f)$. Similar to the proof of Lemma \ref{lem: Norm f_j+1-f_j<f_j-1-f_j}, we estimate
		
		\begin{align*}
		\NN{h}_{\nu,\mu}&\leq \NN{h_0}_{C^\nu_x}+\int_0^T  \NN{\del_t h}_{C^{\nu-\mu t}} dt
		\\&\leq \NN{h_0}_{C^\nu}+\frac1\mu\left(\NN{\velo}_{C^{\nu,\infty}}+4\N{U}R\right)\NN{h}_{\nu,\mu}
		\\&\hspace{2cm}	
		+\frac{\nu}{\mu}\gamma (2+3\eta R)\left(R+\sup_{0\leq t\leq T}\NN{F}_{C^{\nu-\mu t}}\right)\NN{h}_{\nu,\mu}
		\\&\hspace{2cm} +\frac{\nu}{\mu}\gamma \sup_{0\leq t\leq T}\NN{n_f(1-\eta n_f)}_{C^{\nu-\mu t}}\NN{H}_{C^{\nu-\mu t}}
		%		\\&\hspace{2cm} +T\gamma \sup_{0\leq t\leq T}\NN{n_g(1-\eta n_g)}_{C^{\nu-\mu t}_x}\NN{H}_{C^{\nu-\mu t}_x}
		\end{align*}
		for $\mu\geq\tilde{\mu}$. Using $\NN{F(t)}_{C^{\nu-\mu t}}\leq C$ and choosing $\mu>0$ sufficiently large, we obtain 
		%		 and \[\sup_{0\leq t\leq T}\NN{n_g(1-\eta n_g)}_{C^{\nu-\mu t}_x}\leq R(1+\eta R),\]
		
		\begin{align}
		\NN{h}_{\nu,\mu}&\leq \frac12	\NN{h}_{\nu,\mu}+\NN{h_0}_{C^\nu} +\frac{\gamma\nu}{\mu}\sup_{0\leq t\leq T}\NN{n_f(1-\eta n_f)}_{C^{\nu-\mu t}}\NN{H}_{C^{\nu-\mu t}}
		\end{align}
		implying
		
		\begin{align}\label{eq.strong.Lipschitz}
		\NN{f-g}_{\nu,\mu}&\leq 
		2\NN{f_0-g_0}_{C^\nu} +\frac{2\gamma\nu}{\mu}\sup_{0\leq t\leq T}\NN{n_f(1-\eta n_f)}_{C^{\nu-\mu t}}\NN{F-G}_{C^{\nu-\mu t}}
		%		\\&\leq2\NN{f_0-g_0}_{C^\nu} +\frac{2}{\sqrt{\mu}}\sup_{0\leq t\leq T}\NN{F-G}_{C^{\nu-\mu t}}.\notag
		\end{align}
		for $\mu>\tilde{\mu}$ being sufficiently large.	
		Moreover, we can again use the submultiplicative property of the norm $\NN{\cdot}_{C^{\nu-\mu t}}$ and the fact that 
		
		\[\|n_f\|_{C^{\nu-\mu t}}\leq \|f\|_{C^{\nu-\mu t}}\leq \|f\|_{\nu,\mu}\leq R\] to obtain that
		
		\begin{align*}
		\frac{\gamma\nu}{\mu}\NN{n_f(1-\eta n_f)}_{C^{\nu-\mu t}}
		&\leq \frac{\gamma\nu}{\mu} R\left(1+\eta R\right)
		\leq \frac{1}{\sqrt{\mu}} 
		\end{align*}
		if $\mu>\tilde \mu$ is sufficiently large. This finishes the proof.
	\end{proof}
	\section{BGK-type collision operator}
	In this section, we focus on the semiconductor Boltzmann-Dirac-Benney equation 
	
	\begin{equation}\label{3.be.again}\tag{\ref{3.be}}
	\pa_t f + u(p)\cdot\na_x f - U\na_x n_f\cdot\na_p f = \gamma n_f(1-\eta n_f)(\F_f-f)
	\end{equation}
	with $f(x,p,0)=f_0(x,p)$ for given $U\neq0$ and $\gamma\geq 0$.
	
	It can also be understood as a version of Eq. \eqref{2.be} with a self-consistent equilibrium distribution function $ F= \F_f(x,p,t) = \big(\eta + \exp(-\lambda_0(x,t)-\lambda_1(x,t)\epsilon(p))\big)^{-1}$, for $x\in\R^d,\ p\in\T^d$ and $t>0$. Here, $\lambda_0,\lambda_1$ shall be chosen in such a way that
	
	\begin{equation}
	n_f(x,t):=n_{\F_f}(x,t)\quad\mbox{and}\quad E_f(x,t)=E_{\F_f}(x,t),
	\end{equation}
	where $n_f(x,t):=\int_{\T^d}f(x,p,t)dp$ and $E_f(x,t)=\int_{\T^d}\epsilon(p)f(x,p,t)dp$. This is well-defined according to \cite{Bra17} section 5.1. 
	%Thus, we change the equation \eqref{2.be} to

	%Thus, replacing $F$ with $\F_f$ we arrive at equation \eqref{3.be.nochmal}, namely
	%
	%\begin{equation}\label{3.be.again}
	%\pa_t f + u(p)\cdot\na_x f - U\na_x n_f\cdot\na_p f = \gamma n_f(1-\eta n_f)(\F_f-f)
	%\end{equation}
	\begin{theorem}\label{thm.eins}
		Let $\eta,\nu,R>0$ and $\alpha>0$. There exist $\delta>0$ and $\mu_0>0$ such that the following is true:
		
		Let $\bar n,\bar E\in\R$ and  $f_0:\R^d\times\T^d\to [2\alpha,\eta^{-1}-2\alpha]$ be analytic such that
		
		\begin{equation}
		\NN{f_0}_{C^{\nu}}\leq \frac{R}2\quad\mbox{and}\quad
		\NN{n_{f_0}-\bar n}_{C^{\nu}}+\NN{E_{f_0}-\bar E}_{C^{\nu}}\leq \frac12\delta.
		\end{equation}
		Then if $\mu\geq\mu_0$, $T\in (0,\nu/\mu)$, equation \eqref{2.be} admits a unique analytic solution $f:\R^d\times\T^d\times[0,T)\to[\alpha,\eta^{-1}-\alpha]$ satisfying $f\big|_{t=0}= f_0$ and 
		\[
		\NN{f(t)}_{C^{\nu-\mu t}}\leq R\quad\mbox{and}\quad
		\NN{n_{f}(t)-\bar n}_{C^{\nu-\mu t}}+\NN{E_{f}(t)-\bar E}_{C^{\nu-\mu t}}\leq \delta
		\]
		for all $0\leq t<T$. Moreover, let $f,g$ be  the unique solution of \eqref{3.be.again} for with $f(x,p,0)=f_0(x,p)$ and $g(x,p,0)=g_0(x,p)$, where $f_0$ and $g_0$ satisfy both the hypothesis of this theorem. Then there exists a $C>0$ such that
		
		\begin{equation*}
		\NN{f(t)-g(t)}_{C^{\nu-\mu t}}\leq C\NN{f_0-g_0}_{C^\nu}
		\end{equation*}
		for all for all $0\leq t<T$.
	\end{theorem}
	In order to prove that Eq.~\eqref{3.be.again} admits a local, analytic solution, we basically require Theorem \ref{thm: boltzmann local: hT solution} %, Proposition \ref{prop: boltzmann local: hT Lipschitz}
	and the following Lipschitz estimate from Proposition \ref{prop.Ff.is.Lipschitz}. 
	
	Let $\Psi:(f_0,F)\mapsto f$ be the mapping as in Theorem \ref{thm: boltzmann local: hT solution} defined by the solution of 
	
	\begin{equation*}
	\del_t f+\velo\cdot\nabla_x f - U\nabla_xn_f\cdot\nabla_p f =\gamma n_f(1-\eta n_f)(F-f)
	\end{equation*}
	with $f(x,p,0)=f_0(x,p)$. With this, we define the mapping 
	
	\begin{equation*}
	\Theta(g):=\Psi(f_0,\F_g).
	\end{equation*}
	Therefore, every fixed-point of $\Theta$ is a classical solution of \eqref{3.be.again}. At first, we need to show that $\Theta$ is well-defined.
	\begin{proposition}\label{lem.Ff.is.Lipschitz}
		Let $\eta,\nu,R>0$ and $\alpha>0$. There exists an $C,\delta>0$ such that the following is true.
		
		Let $f,g:\R^d\times\T^d\to[\alpha,\eta^{-1}-\alpha]$ be analytic satisfying $\NN{f}_{C^{\nu}},\NN{g}_{C^{\nu}}\leq R$ and
		
		\begin{equation*}
		\NN{n_h-\bar n}_{C^{\nu}}+\NN{E_h-\bar E}_{C^{\nu}}\leq \delta\quad\mbox{for }h\in\{f,g\}
		\end{equation*}
		and some $\bar{n},\bar E\in \R$, it holds
		
		\begin{equation*}
		\NN{\F_f}_{C^\nu},\NN{\F_g}_{C^\nu}\leq C
		\end{equation*}
		and
		
		\begin{equation*}
		\NN{\F_f-\F_g}_{C^\nu}\leq C\NN{f-g}_{C^\nu}.
		\end{equation*}
	\end{proposition}
	\begin{proof}
		See appendix.
	\end{proof}
	Using this proposition, we can define the metric space $Y$ on which $\Theta$ is a contraction.
	\begin{definition}
		For $R,\nu,\eta,\alpha>0$, let $\delta>0$ be as in Proposition \ref{lem.Ff.is.Lipschitz}. 
		%	which only depends on $a$ according to Remark \ref{rem.alpha.depends.only.on.a}. Given $R,\nu>0$,
		% let $\delta>0$ be as in Proposition \ref{prop.Ff.is.Lipschitz}. 
		Moreover, let $\bar n,\bar E\in \R$, $\mu>0$ and $T\in(0,\nu/\mu)$.
		\\
		%	We define $X$ as the complete metric space of all analytic functions $f:\R^d\times\T^d\times[0,T)\to \R$ with $\NNN{f}_{\nu,\mu}\leq R$. 
		%%	Note that this is a complete metric space if the metric is induced by the norm $\NN{\cdot}_{\nu,\mu}$.
		We define $Y$  as the space of all analytic functions $f:\R^d\times\T^d\times[0,T)\to [\alpha,\eta^{-1}-\alpha]$ satisfying
		\begin{enumerate}
			\item $\NN{f}_{\nu,\mu}\leq R$,
			\item
			$\NN{n_{f}(\cdot,t)-\bar n}_{C^{\nu-\mu t}}+\NN{E_{f}(\cdot,t)-\bar E}_{C^{\nu-\mu t}}\leq \delta $ 
		\end{enumerate}
		for all $t\in[0,T)$. 
		% $\NNN{f}_{\nu,\mu}\leq R$ such that the mapping	$(n_f,E_f):\R^d\times[0,t)\to \overline{\mathcal U_{a,\alpha}}$ defined by \[(n_f,E_f):=\int_{\T^d}(1,\epsilon(p))f(\cdot,p,\cdot)dp\] is well-defined and fulfills	
		%	\[\NN{n_{f}(\cdot,s)}_{\dot C^{\nu-\mu s}_x}+\NN{E_{f}(\cdot,s)}_{\dot C^{\nu-\mu s}_x}\leq \delta n_{f}(x,s)(1-\eta n_{f}(x,s))\quad\mbox{for }x\in\R^d,s\in[0,t).\]
		Thus, $Y$ is a complete if the metric is induced by the norm $\NN{\cdot}_{\nu,\mu}$.
	\end{definition}
	As we plan to apply the Banach fixed-point theorem, we need to show that $\Theta$ is a contraction, i.e., the image of $\Theta$ is included in $Y$ and $\Theta$ is Lipschitz continuous with Lipschitz constant $L<1$. 
	\begin{lemma}
		Let $\bar n,\bar E\in \R$, $f_0:\R^d\times\T^d\to [2\alpha,\eta^{-1}-2\alpha]$ be analytic such that
		
		\begin{equation*}
		\NN{f_0}_{C^{\nu}}\leq \frac R2\quad\mbox{and}\quad
		\NN{n_{f_0}-\bar n}_{C^{\nu}},\NN{E_{f_0}-\bar E}_{C^{\nu}}\leq \frac12\delta.
		\end{equation*}
		If $\mu>0$ is sufficiently large and $g\in Z$, then $\Theta(g)\in Z$.
	\end{lemma}
	\begin{proof}
		By definition, we have 
		\[\Theta(g)=\Psi(f_0,\F_g).\]
		For $g\in Z$, we know from Proposition \ref{lem.Ff.is.Lipschitz} that $\NN{\F_g}_{C^{\nu-\mu t}}\leq C$ for some $C>0$ and all $t\in[0,T)$. Hence, $f:=\Theta(g)$ is well-defined for sufficiently large $\mu>0$ and $\NN{f}_{\nu,\mu}\leq R$. Clearly, by continuity, if $\mu$ sufficiently large and thus $T>0$ sufficiently small, then the image of $f$ belongs to $[\alpha,\eta^{-1}-\alpha]$. 
		
		Therefore, it remains to show that
		$\NN{n_{f}(\cdot,t)}_{C^{\nu-\mu t}}+\NN{E_{f}(\cdot,t)}_{C^{\nu-\mu t}}\leq \delta $ for all $t\in[0,T)$. 	Using Lemma \ref{lem: d_t nnormClambda  f estimate} entails that
		
		\begin{align*}
		\NN{n_{f}(\cdot,t)-\bar n}_{C^{\nu-\mu t}}
		-
		\NN{n_{f_0}-\bar n}_{C^{\nu}}\leq\int_0^T\NN{\del_t n_f}_{C^{\nu-\mu t}}d t
		\leq 
		\int_0^T\NN{\del_t f}_{C^{\nu-\mu t}}d t.
		\end{align*}
		and likewise

		\begin{align*}
		\NN{E_{f}(\cdot,t)-\bar E}_{C^{\nu-\mu t}}
		-
		\NN{E_{f_0}-\bar E}_{C^{\nu}}
		\leq 
		\|\epsilon\|_{L^\infty}\int_0^T\NN{\del_t f}_{C^{\nu-\mu t}}d t.
		\end{align*}
		As in the proof of Lemma \ref{lem.proof.thm.boltzmann local: hT solution1} (see inequality \ref{eq.inproof.73}), we can show that
		
		\begin{align*}
		\int_0^T\NN{\del_tf}_{C^{\nu-\mu t}}d t
		\leq 
		\frac{1}{\mu}\left(\NN{\velo}_{C^{\nu,\infty}}+\N{U}R\right)\NN{f}_{\nu,\mu}
		+
		\gamma T R\left(1+\eta R\right)\left(R+C\right)
		\leq \frac{\tilde C}\mu
		\end{align*}
		for some $\tilde C>0$. Thus,
		
		\begin{multline*}
		\NN{n_{f}(\cdot,t)-\bar n}_{C^{\nu-\mu t}}+\NN{E_{f}(\cdot,t)-\bar E}_{C^{\nu-\mu t}}
		\\	\leq
		\NN{n_{f_0}-\bar n}_{C^{\nu}}+\NN{E_{f_0}-\bar E}_{C^{\nu}}
		%(1+\|\epsilon\|_{C^\nu})\int_0^T\NN{\del_t n_f}_{C^{\nu-\mu t}}d t	
		+\frac{\bar C}\mu\leq \frac\delta 2 +\frac{\bar C}\mu\leq \delta
		\end{multline*}
		for all $t\in[0,T)$ and some $\bar C>0$ if $\mu\geq \tilde 2C/\delta$, which proves the assertion.
	\end{proof}
	\begin{lemma}
		Let $\bar n,\bar E\in \R$, $f_0:\R^d\times\T^d\to [2\alpha,\eta^{-1}-2\alpha]$ be analytic such that
		
		\begin{equation*}
		\NN{f_0}_{C^{\nu}}\leq \frac R2\quad\mbox{and}\quad
		\NN{n_{f_0}-\bar n}_{C^{\nu}},\NN{E_{f_0}-\bar E}_{C^{\nu}}\leq \frac12\delta.
		\end{equation*}
		If $\mu>0$ is sufficiently large, then for $f,g\in Y$ it holds
		\[\NN{\Theta(f)-\Theta(g)}_{\nu,\mu}\leq \frac12\NN{f-g}_{\nu,\mu}.\]
	\end{lemma}
	\begin{proof}
		According to the previous Lemma, we can apply Theorem \ref{thm: boltzmann local: hT solution} entailing  for sufficiently large $\mu>0$ that
		
		\begin{align*}
		\NN{\Theta(f)-\Theta(g)}_{\nu,\mu}\leq 2\mu^{-\frac12}\sup_{t}\NN{\F_f-\F_g}_{C^{\nu-\mu t}}.
		\end{align*}
		Then the second statement of Proposition \ref{lem.Ff.is.Lipschitz} yields that 
		
		\begin{align*}
		\NN{\Theta(f)-\Theta(g)}_{\nu,\mu}\leq C\mu^{-\frac12}\sup_{t}\NN{f-g}_{C^{\nu-\mu t}}\leq C\mu^{-\frac12}\NN{f-g}_{\nu,\mu} 
		\end{align*}
		for some $C>0$. This implies the assertion for sufficiently large $\mu$ satisfying $\mu\geq 4C^2$.
	\end{proof}
	
	\begin{proof}[Proof of Theorem \ref{thm.eins}]
		The contraction mapping theorem ensures that $\Psi$ has a unique fixed-point implying that equation \eqref{3.be.again} admits a unique solution. Finally, the Lipschitz estimate is a direct consequence of Theorem \ref{thm: boltzmann local: hT solution}.
	\end{proof}
	With Theorem \ref{thm.eins} we can now easily prove the following weaker version of Theorem  \ref{thm.BGK.local.solution.bitorus}.
	\begin{theorem}\label{thm.BGK.local.solution.bitorus.easy}
		Let $\eta>0$, $\gamma\geq0,$ $U\neq0$ and $f_0:\T^d\times\T^d\to(0,\eta^{-1})$ be analytic such that
		
		\begin{align*}
		n_{f_0}(x)=\int_{\T^d}f_0(x,p)dp=const. \quad\mbox{and}\quad E_{f_0}=\int_{\T^d}\epsilon(p)f_0(x,p)dp=const.
		\end{align*}
		w.r.t.\ $x\in\T^d$. Then there exists a time $T>0$ such that \eqref{3.be} admits a unique analytic solution $f:\T^d\times\T^d\times[0,T)\to\R$ with $f(x,p,0)=f_0(x,p)$.
	\end{theorem}
	\begin{proof}
		Since $f_0$ is analytic and hence continuous, there exists a $\alpha>0$ such that $2\alpha<f_0<\eta^{-1}-2\alpha$. The key difference to Theorem \ref{thm.eins} is that now the spacial domain is essentially restricted to a compact set $\T^d$, which can be extended periodically to $\R^d$. Any analytic function $f_0$ on a compact domain has a minimal radius $r$ of convergence, i.e.~a number $r>0$ such that for all $(x,p)$ the series 
		
		\begin{equation*}
		\sum_{i,j\in\mathbb N_0^d}\frac{\del_x^i\del_p^j f_0(x,p)}{i!j!}x^ip^j
		\end{equation*}
		converges absolutely for $|x|+|p|\leq r$. This implies that
		
		\begin{equation*}
		M(x,p):=\sum_{i,j\in\mathbb N_0}\frac{(r/2)^{|i+j|}}{i!j!}\N{\del_x^i\del_p^j f_0(x,p)}<\infty
		\end{equation*}
		for all $(x,p)\in(\T^d)^2$. Now, choose $(x_1,p_1),\ldots,(x_N,p_N)\in(\T^d)^2$ such that
		\[\bigcup_{i=1}^N B_{r/4}(x_i,p_i)\supset(\T^d)^2\]
		and define
		\[K:= \max_{i=1,\ldots,N}M(x_i,p_i).\]
		Let $\nu:=r/8$. Then for every $(x,p)\in(\T^d)^2$ there exists an $i\in\{1,\ldots ,N\}$ such that $\N{(x-x_i,p-p_i)}\leq \nu$ and
		
		\begin{align*}
		\sum_{k,j\in\mathbb N_0^d}\frac{(2\nu)^{|k+j|}}{k!j!}\N{\del_x^k\del_p^j f_0(x,p)}
		\leq 
		\sum_{k,j\in\mathbb N_0^d}\frac{(4\nu)^{|k+j|}}{k!j!}\N{\del_x^k\del_p^j f_0(x_i,p_i)}\leq K.
		\end{align*}
		This directly implies that $\N{f_0}_{C^{2\nu}}<\infty$. Moreover, by Lemma \ref{lem: nnorm leq C norm 1}, we can see that also $R:=\NN{f_0}_{C^{\nu}}+1$ is finite as $\nu<2\nu$. As $\bar n:= n_{f_0}$ and $\bar E:=E_{f_0}$ are constant, we have shown all the hypothesis of Theorem \ref{thm.eins} and finally obtain a analytic solution on a small time interval.
	\end{proof}
	For a full proof of Theorem  \ref{thm.BGK.local.solution.bitorus}, we refer to section \ref{sec.space.local} and section \ref{sec.BGK.space.local}, in which we refine the presented technique using that the collision operator is local in space. The next section is devoted to an application of Theorem \ref{thm.eins} showing the ill-posedness of equation \eqref{3.be}.

	\section{On the ill-posedness of the semiconductor Boltzmann-Dirac-Benney equation}
	This section is motivated by the ill-posedness result of \cite{HaNg15}  and \cite{BaNo12} for the Vlasov-Dirac-Benney equation. Similar to \cite{BaNo12}, we linearize the equation around an equilibrium. Let $\bar{\lambda}=(\bar \lambda_0,\bar \lambda_1)\in \R^2$. Then 
	\[F_{\bar{\lambda}}:\T^d\to\R, \quad  p\mapsto\frac{1}{\eta+e^{-\bar \lambda_0-\bar \lambda_1\epsilon(p)}}\]
	is a stationary analytic solution of \eqref{3.be}, which is constant in $x$. 
	\subsection{Linearized equation}
	Now let us formally linearize the left-hand side of \eqref{3.be} around $F_{\bar{\lambda}}$ and consider
	
	\begin{align}\label{4.be.linear}
	&\pa_t g + u(p)\cdot\na_x g - U\na_x n_g\cdot\na F_{\bar \lambda}(p) 
	%=0
	%\\&= \gamma \del_f(n_f(1-\eta n_f)(\F_f-f))(g)
	%\\&=\gamma \del_n(n(1-\eta n))(n_g)(\F_f-f))|_{f=F_\lambda}+ n_f(1-\eta n_f)\del(\F_f-f))(g) |_{f=F_{\lambda}}
	%\\&=0+ n_f(1-\eta n_f)\del_f(\F_f-f))(g) |_{f=F_{\lambda}}
	=\gamma n_{F_{\bar \lambda}}(1-\eta n_{F_{\bar \lambda}})\left(G(\bar \lambda;p)\cdot\binom{n_g}{E_g}-g\right)
	\end{align}
	with $g(x,p,0)=g_0(x,p)$ 
	and $G(\bar \lambda;p):= \del_{(n_f,E_f)}\F_f(p)|_{f=F_{\bar \lambda}}$. Recall that  $u(p)=\nabla_p\epsilon(p)$ with 
	
	\begin{equation*}
	\epsilon(p) = -2\epsilon_0\sum_{i=1}^d \cos(2\pi p_i), \quad p\in\T^d,
	\end{equation*}
	for some $\epsilon_0>0$.
	\begin{remark}
		The definition of $G(\bar \lambda;p)$ has to be understood according to Definition \ref{def.FrhoE}: $\F_f$ can be written  by $\F_f=\FrhoE(n_f,E_f;p)$ for some analytic $\FrhoE:V\subset\R^2\times\T^d\to[0,\eta^{-1}]$.  By Lemma \ref{lem.Berechnung.von.G} from the appendix, it holds
		\[
		G(\bar\lambda;p)=\frac{F_{\bar\lambda}(p)(1-\eta F_{\bar\lambda}(p))}{\int_{\T^d}\epsilon^2d\mu\int_{\T^d}1d\mu-\big(\int_{\T^d}\epsilon d\mu\big)^2}\int_{\T^d}\binom{-\epsilon(p')}{1}(\epsilon(p)-\epsilon(p'))d\mu_{p'},
		\]
		where $d\mu_p:= F_{\bar\lambda}(p)(1-\eta F_{\bar\lambda}(p))dp$.
		
		In the following, we will denote the components of $G$ as $G_1,G_2$ and write $p=(p_1,\ldots,p_d)$, $x=(x_1,\ldots,x_d)$ and $u(p)=(u_1(p),\ldots u_d(p))$.
	\end{remark}
	\begin{lemma}\label{lem.Def.of.g.cxi}
		For $\bar\lambda\in\R^2$ we abbreviate $\gamma_{\bar{\lambda}}:=  \gamma n_{F_{\bar \lambda}}(1-\eta n_{F_{\bar \lambda}})$. Assume that there exists a bounded set $K\subset(\R\setminus\{0\})^2$ with $\overline K\subset \R\setminus\{0\}\times \R$  such that the eigenvalues of
		
		\begin{equation*}
		B=B(\alpha,\beta):=\int_{\T^d}\frac{(U\velo_1(p)\partial_{p_1}F_{\bar \lambda}(p)+\beta G_1(\bar{\lambda};p),\beta G_2(\bar{\lambda};p))}{\velo_1(p)^2+\alpha^2}\binom{1}{\epsilon(p)}dp
		\end{equation*}
		are $0$ and $1$ for $(\alpha,\beta)\in K$. Let $(\hat n_{\alpha,\beta},\hat E_{\alpha,\beta})$ denote the eigenvector to the eigenvalue $1$ and define
		
		\begin{equation*}
		A_{\alpha,\beta}(p):=\frac{1}{\velo_1(p)-i\alpha}\left(U\partial_{p_1}F_{\bar\lambda}(p)\hat n_{\alpha,\beta} -i\frac{\beta}{\alpha} G(\bar{\lambda};p)\cdot\binom{\hat n_{\alpha,\beta}}{\hat {E}_{\alpha,\beta}}\right).
		\end{equation*}
		Then
		\[
		g_{\alpha,\beta}(x,p):=A_{\alpha,\beta}(p)e^{i\gamma n_{F_{\bar{\lambda}}}(1-\eta n_{F_{\bar{\lambda}}})\frac{\alpha}{\beta}x_1}
		\]
		is a solution of
		
		\begin{equation*}
		\velo(p)\cdot\na_x g_{\alpha,\beta} - U\na_x n_{g_{\alpha,\beta}}\cdot\na F_{\bar \lambda}(p) 
		=\gamma n_{F_{\bar \lambda}}(1-\eta n_{F_{\bar \lambda}})\left(G(\bar \lambda;p)\cdot\binom{n_{g_{\alpha,\beta}}}{E_{g_{\alpha,\beta}}}-\frac{\alpha^2}{\beta}g_{\alpha,\beta}\right).
		\end{equation*}
		Moreover, let $N\in\mathbb N$. There exists $C_N>0$ such that
		
		\begin{align*}
		%	C \velo_1(p)^2\leq
		%	 \N{g_{\alpha,\beta}(x,p)}\leq
		\sup_{(\alpha,\beta)\in K}\NN{g_{\alpha,\beta}}_{W^{N,\infty}(\R^d\times\T^d)}\leq C_N(1+\N{\beta}^{-N})
		\end{align*}
		for all $(x,p)\in\R^d\times\T^d$.
		%	\begin{align*}
		%%		\leq \int_{\T^d}\N{\frac{\velo_1(p)}{\velo_1(p)-ic}}F_{\bar{\lambda}}(p)(1-\eta F_{\bar \lambda}(p))dp\leq 
		%		%\int_{\T^d}\frac{\velo_1(p)^2}{\velo_1(p)^2+c^2}F_{\bar{\lambda}}(p)(1-\eta F_{\bar \lambda}(p))dp.
		%	\NN{g_{c,\varphi}(x,p)}_{L^1_p(\T^d)}\leq \frac{1}{4\eta}	\quad\mbox{and}\quad  \NN{g_{c,\varphi}(0,p)}_{L^1_p(\T^d)}=
		%		\frac{1}{U\bar \lambda_1}.
		%	\end{align*}
		In addition, there exists a $\nu_0>0$ and such that
		
		\begin{align*}%\label{eq.lem.g.c.varphi}
		\NN{g_{\alpha,\beta}}_{C^\nu}\leq C_{\nu_0}(1+|\beta|^{-1}) e^{\frac{c\nu}{\N\beta}}
		\end{align*}
		and for all $\nu\leq \nu_0$ and some $c,C_{\nu_0}>0$ being independent from $\alpha,\beta$.
		%	Let $c:=\alpha-\frac{\beta}{\alpha}$ and $\xi :=\frac{\alpha\gamma_{\bar{\lambda}}}{\beta}$ if $\frac{\alpha\gamma_{\bar{\lambda}}}{\beta}\neq0$ or $\xi$ be arbitrary if $\gamma_{\bar\lambda}=0$. Then
		%	\begin{equation*}
		%	g_{c,\xi}(x,p):= A_{\alpha,\beta}(p)e^{i\xi x_1}
		%	\end{equation*}
		%	is a solution of
		%	\begin{equation*}
		%	c\xi g_{c,\xi} + \velo(p)\cdot\na_x g_{c,\xi} - U\na_x n_{g_{c,\xi}}\cdot\na F_{\bar \lambda}(p) 
		%	=\gamma n_{F_{\bar \lambda}}(1-\eta n_{F_{\bar \lambda}})\left(G(\bar \lambda;p)\cdot\binom{n_{g_{c,\xi}}}{E_{g_{c,\xi}}}-g_{c,\xi}\right).
		%	\end{equation*}
	\end{lemma}
	\begin{proof}
		Note that $G$ is symmetric an $\partial_{p_1}$ is anti-symmetric, i.e., $G(\bar{\lambda};-p)=G(\bar{\lambda};p)$ and $\partial_{p_1}F_{\bar{\lambda}}(-p)=-\partial_{p_1}F(p)$, which is a consequence of $\velo(-p)=-\velo(p)$ as well as $\F_{\bar{\lambda}}(-p)=\F_{\bar{\lambda}}(p)$ for $p\in \T^d$.  Therefore, since the denominator is even, we may add an odd function to the denominator without changing the integral. Thus, we can divide the integrand by $\velo_1(p)+i\alpha$ and obtain
		
		\begin{equation*}
		B(\alpha,\beta)=\int_{\T^d}\frac{(U\partial_{p_1}F_{\bar \lambda}(p)-i\frac{\beta}{\alpha} G_1(\bar{\lambda};p),-i\frac{\beta}{\alpha} G_2(\bar{\lambda};p))}{\velo_1(p)-i\alpha}\binom{1}{\epsilon(p)}dp.
		\end{equation*}
		Since $(\hat n_A,\hat E_A)$ is the eigenvector to the eigenvalue $1$ of $B$, we infer $\int_{\T^d}A_{\alpha,\beta}(p)dp=\hat n_{\alpha,\beta}$ and $\int_{\T^d}\epsilon(p)A_{\alpha,\beta}(p)dp=\hat E_{\alpha,\beta}$. Finally, we directly compute 
		
		\begin{align*}
		\gamma &n_{F_{\bar{\lambda}}}(1-\eta n_{F_{\bar{\lambda}}})\frac{\alpha^2}{\beta} g_{\alpha,\beta} + \velo(p)\cdot\na_x g_{\alpha,\beta}
		\\&=i\gamma n_{F_{\bar{\lambda}}}(1-\eta n_{F_{\bar{\lambda}}})\frac{\alpha}{\beta}\left(-i\alpha+\velo_1(p)\right)A_{\alpha,\beta}(p)e^{i\gamma n_{F_{\bar{\lambda}}}(1-\eta n_{F_{\bar{\lambda}}})\frac{\alpha}{\beta}x_1}
		\\&=i\gamma n_{F_{\bar{\lambda}}}(1-\eta n_{F_{\bar{\lambda}}})\frac{\alpha}{\beta}\left(U\partial_{p_1}\F_{\bar\lambda}(p)\hat n_{\alpha,\beta} -i\frac\beta\alpha G(\bar{\lambda};p)\cdot\binom{\hat n_{\alpha,\beta}}{\hat E_{\alpha,\beta}}\right)e^{i\gamma n_{F_{\bar{\lambda}}}(1-\eta n_{F_{\bar{\lambda}}})\frac{\alpha}{\beta} x_1}
		\\&=U\na_x n_{g_{\alpha,\beta}}\cdot\na F_{\bar \lambda}(p) 
		+\gamma n_{F_{\bar{\lambda}}}(1-\eta n_{F_{\bar{\lambda}}})G(\bar \lambda;p)\cdot\binom{n_{g_{\alpha,\beta}}}{E_{g_{\alpha,\beta}}}.
		\end{align*}
		Since $p\mapsto g_{\alpha,\beta}(0,p)$ is analytic on $\T^d$ and $\T^d$ is compact and $\overline K\subset \R\setminus\{0\}\times \R$ is compact, there exists a $\nu>0$ such that
		
		\begin{align*}
		C_\nu:=\sup_{(\alpha,\beta)\in K}\sum_{j\in\mathbb N_0^d} \frac{\nu^j}{j!}\NN{\del^j_p g_{\alpha,\beta}(0,p)}_{W^{1,1}_p(\T^d)}<\infty.
		\end{align*}
		Thus,
		
		\begin{align*}
		\NN{g_{\alpha,\beta}}_{C^\nu}
		&=
		\sum_{j,l\in\mathbb N_0^d} \frac{\nu^{|j+l|}}{j!l!}\NN{\del^l_x\del^j_p g_{\alpha,\beta}(x,p)}_{W^{1,\infty}_xW^{1,1}_p(\T^d)}
		\\&=
		\sum_{j\in\mathbb N_0^d} \frac{\nu^{|j|}}{j!}\NN{\del^j_p g_{\alpha,\beta}(0,p)}_{W^{1,1}_p(\T^d)}\sum_{k=0}^1\sum_{l=0}^\infty \frac{\nu^{l}}{l!}\N{\del^{l+k}_{x_1} e^{-i\varphi x_1}}
		\\&=C_\nu (1+|\varphi|)e^{\nu\N{\varphi}}
		\end{align*}
		setting $\varphi:=\gamma n_{F_{\bar{\lambda}}}(1-\eta n_{F_{\bar{\lambda}}})\frac{\alpha}{\beta}$ for all $(\alpha,\beta)\in K$. If we want to estimate only a finite number of derivatives, we see that for all $N>0$ there exists a $C_N>0$ such that $\sup_{(\alpha,\beta)\in K}\sum_{|j|\leq N}\NN{\del^j_pg_{\alpha,\beta}(0,p)}_{L^\infty_p(\T^d)}\leq C_N$ since $g_{\alpha,\beta}$ is smooth. This yields
		
		\begin{align*}
		\sum_{\substack{j,l\in\mathbb N_0^d\\  |j|,|l|\leq N}}\NN{\del^l_x\del^j_pg_{\alpha,\beta}}_{L^\infty(\R^d\times \T^d)}
		&=
		\sum_{|j|\leq N}\NN{\del^j_pg_{\alpha,\beta}(0,p)}_{L^\infty_p( \T^d)}\sum_{|l|\leq N}\NN{\del^l_x e^{-i\phi x_1}}_{L^\infty(\R^d)}\\&
		\leq C_N\sum_{l=0}^{N}\N{\varphi}^l
		\leq NC_N(1+\N{\varphi}^N).\qedhere
		\end{align*}
		%	???????Moreover, we have that 
		%	$\N{\Re(h(0,p))},\N{\Im(h(0,p))}\geq C\velo_1(p)^2$ for some $C>0$ independent from $\varphi$ for all $p\in\T^d$.
		%	\begin{align*}
		%	c\xi& g_{c,\xi} + \velo(p)\cdot\na_x g_{c,\xi} +\gamma_{\bar{\lambda}}g_{c,\xi} 
		%	\\&=i\xi\left(-i\left(c+\frac{\gamma_{\bar{\lambda}}}{\xi}\right)+\velo_1(p)\right)A_{\alpha,\beta}(p)e^{i\xi x_1}
		%	\\&=i\xi\left(U\partial_{p_1}\F_{\bar\lambda}(p)\hat n_A -i\frac{\gamma_{\bar{\lambda}}}{\xi}G(\bar{\lambda};p)\cdot\binom{\hat n_A}{\hat E_A}\right)e^{i\xi x_1}
		%	\\&=U\na_x n_{g_{c,\xi}}\cdot\na F_{\bar \lambda}(p) 
		%	+\gamma_{\bar \lambda}G(\bar \lambda;p)\cdot\binom{n_{g_{c,\xi}}}{E_{g_{c,\xi}}}. \qedhere
		%	\end{align*}
	\end{proof}
	In order to prove that the hypothesis of the previous lemma can be fulfilled, we start with an easier case, where $\beta=0$. Then the condition simplifies to
	\begin{equation}\label{eq.Bed.for.lambda,c}
	1=U \bar\lambda_1\int_{\T^d}\frac{\velo_1(p)^2}{\velo_1(p)^2+\alpha_0^2}F_{\bar{\lambda}}(p)(1-\eta F_{\bar \lambda}(p))dp
	\end{equation}
	for some $\alpha_0\neq 0$.
	\begin{lemma}\label{lem.Bed.for.lambda,c}
		Let $U\neq0$.	Then there exist $\bar \lambda\in\R^2$ and an $\alpha_0>0$ such that \eqref{eq.Bed.for.lambda,c} is satisfied. In addition, the solution $\alpha_0$ of \eqref{eq.Bed.for.lambda,c} is unique (up to its sign) for fixed $\bar{\lambda}$.
		%	\begin{align*}%\label{eq.Bed.for.lambda,c}
		%	1=U \bar\lambda_1\int_{\T^d}\frac{\velo_1(p)^2}{\velo_1(p)^2+c^2}F_{\bar{\lambda}}(p)(1-\eta F_{\bar \lambda}(p))dp.
		%	\end{align*}
	\end{lemma}
	
	\begin{proof}
		At first, we define
		\[
		\kappa(\lambda):=\lambda_1\int_{\T^d}F_{\lambda}(p)(1-\eta F_{\lambda}(p))dp.
		\]
		According to \cite{Bra17} section 5.3, it holds that $\sup_{\lambda\in\R^2}\kappa(\lambda)=\infty$ and by symmetry $\inf_{\lambda\in\R^2}\kappa(\lambda)=-\infty$. Thus, there exists $\bar \lambda\in\R^2$ such that
		\[U \bar\lambda_1\int_{\T^d}F_{\bar{\lambda}}(p)(1-\eta F_{\bar \lambda}(p))dp=U\kappa(\bar \lambda)>1.\]
		Finally, by
		\[
		1<U\kappa(\bar \lambda)\stackrel{c\to0}{\leftarrow}U \bar\lambda_1\int_{\T^d}\frac{\velo_1(p)^2}{\velo_1(p)^2+c^2}F_{\bar{\lambda}}(p)(1-\eta F_{\bar \lambda}(p))dp\stackrel{c\to\infty}{\to}0,
		\]
		the intermediate value theorem yields the first assertion. The uniqueness is a consequence of the monotonicity of $U \bar\lambda_1\int_{\T^d}\frac{\velo_1(p)^2}{\velo_1(p)^2+c^2}F_{\bar{\lambda}}(p)(1-\eta F_{\bar \lambda}(p))dp$ w.r.t.~$c$.
	\end{proof}\begin{remark}
	We used in the proof that Equation \eqref{eq.Bed.for.lambda,c} admits a solution if 
	\[
	1<U \bar\lambda_1\int_{\T^d}F_{\bar{\lambda}}(p)(1-\eta F_{\bar \lambda}(p))dp
	\]
	is satisfied.
\end{remark}
Now, we go back to the general case, where $\beta\neq0$. 
\begin{lemma}\label{lem.Loesung.der.EW.Gleichung}
	Let $\bar{\lambda}$ and $\alpha_0$ be as in Lemma \ref{lem.Bed.for.lambda,c}. There exist an open Interval $0\in I\subset \R$ and a function $\alpha:I\to\R$ with $\alpha(0)=\alpha_0$ such that $B(\alpha(\beta),\beta)$ possesses the eigenvalues $0$ and $1$ for all $\beta\in I$.
\end{lemma}
\begin{proof}
	According to Lemma \ref{lem.Bed.for.lambda,c}, we know that $1$ is an eigenvalue of $B(\alpha_0,0)$ which is equivalent to $\det(B(\alpha_0,0)-\mathrm{Id})=0$. Since \[\phi:\R\setminus\{0\}\times\R\to \R,\quad (a,b)\mapsto \det(B(a,b)-\mathrm{Id})\]
	is smooth, there exists an $\alpha:I\ni0\to \R$ with $\alpha(0)=\alpha_0$ if the derivative of $\phi$ has full rank at $(a,b)=(\alpha_0,0)$. In order to show this, we only need to look at the derivative w.r.t.~$a$:
	
	\begin{align*}
	\partial_a\phi&= \partial_a((B_{11}-1)(B_{22}-1)-B_{12}B_{21})
	\\&=\partial_a B_{11}(B_{22}-1)+(B_{11}-1)\partial_aB_{22}-\partial_aB_{12}B_{21}-B_{12}\partial_aB_{21} 
	\end{align*}
	%	According to Lemma \ref{lem.Bed.for.lambda,c}, we know that $B_{11}(\alpha_0,0)=1$. Moreover, we see directly from the definition of $B$ that $B_{12}(\alpha_0,0)=\partial_a B_{12}(\alpha_0,0)=B_{22}(\alpha_0,0)=0$. Hence, 
	
	\begin{align*}
	\partial_a\phi(\alpha_0,0)=-\partial_aB_{11}(\alpha_0,0)
	&=
	2\alpha_0U\int_{\T^d}\frac{\velo_1(p)^2F_{\bar \lambda}(p)(1-\eta F_{\bar \lambda}(p))}{(\velo_1(p)^2+\alpha_0^2)^2}dp\neq0.
	\end{align*}
	Thus, the derivative of $\phi$ has at $(\alpha_0,0)$ full rank and therefore the zero-set of $\phi$ is locally a one-dimensional manifold at $(\alpha_0,0)$. According to Lemma \ref{lem.Bed.for.lambda,c}, $\phi(a,0)=0$ has only one positive solution at $a=\alpha_0$. Finally, the fact that $B$ has rank $1$ implies directly the trivial eigenvalue and finishes the proof.
\end{proof}
\begin{proposition}\label{lem.g.c.varphi}
	Let $\alpha_0>0,\bar{\lambda}\in\R^2$ be a solution of \begin{equation*}%\label{eq.Bed.for.lambda,c}
	1=U \bar\lambda_1\int_{\T^d}\frac{\velo_1(p)^2}{\velo_1(p)^2+\alpha_0^2}F_{\bar{\lambda}}(p)(1-\eta F_{\bar \lambda}(p))dp
	\end{equation*}
	(see Lemma \ref{lem.Bed.for.lambda,c}). Then there exist an open interval $I\ni 0$ and function $\omega:I\setminus\{0\}\to\R$ and analytic $g_\beta:\T^d\times\T^d\times[0,\infty)\to\R$ solutions of \eqref{4.be.linear}  for $\beta\in I\setminus\{0\}$	such that the following holds:
	\begin{itemize}
		\item $\beta\mapsto \beta\omega(\beta)$ can be extended on $I$ to a positive continuous function.
		\item $g_\beta(x,p,t)=g_\beta(x,p,0)e^{\omega(\beta)t}$.
		%Moreover,
		\item There exists a $\nu_0>0$ and such that
		
		\begin{align}\label{eq.lem.g.c.varphi}
		\NN{g_{\beta}}_{C^\nu}\leq C_{\nu_0}(1+|\beta|^{-1}) e^{\frac{c\nu}{\N\beta}+\omega(\beta)t}\quad\mbox{for all }\beta\in I\setminus\{0\}
		\end{align}
		and for all $\nu\leq \nu_0$ and some $c,C_{\nu_0}>0$ being independent from $x$ and $\beta$.
		% In addition,  
		\item There exists $C_N>0$ such that
		
		\begin{align}\label{eq.lem.g.c.varphi2}
		%	C \velo_1(p)^2\leq
		%	 \N{g_{\alpha,\beta}(x,p)}\leq
		\NN{g_{\beta}(\cdot,\cdot,t)}_{W^{N,\infty}(\R^d\times\T^d)}\leq C_N(1+\N{\beta}^{-N})e^{\omega(\beta)t}\quad\mbox{for all }\beta\in I\setminus\{0\}.
		\end{align}
		
	\end{itemize}
\end{proposition}
\begin{proof}
	Let 
	%	$\alpha_0>0,\bar{\lambda}\in\R^2$ be a solution of \begin{equation*}%\label{eq.Bed.for.lambda,c}
	%	1=U \bar\lambda_1\int_{\T^d}\frac{\velo_1(p)^2}{\velo_1(p)^2+\alpha_0^2}F_{\bar{\lambda}}(p)(1-\eta F_{\bar \lambda}(p))dp
	%	\end{equation*}
	%	(see Lemma \ref{lem.Bed.for.lambda,c}). Then let
	$\alpha:I\to \R$ with $\alpha(0)=\alpha_0$ be given by Lemma \ref{lem.Loesung.der.EW.Gleichung}. For $\beta\in I\setminus\{0\}$, we define
	
	\begin{align*}
	g_{\beta}(x,p,t):=\Re (g_{\alpha(\beta),\beta}(x,p))e^{\omega(\beta)t}\quad\mbox{and}\quad \omega(\beta):=\gamma n_{F_{\bar{\lambda}}}(1-\eta n_{F_{\bar{\lambda}}})\frac{\alpha(\beta)^2-\beta}{\beta}
	\end{align*} 
	and where $g_{\alpha,\beta}$ is given by Lemma \ref{lem.Def.of.g.cxi}. Then $g_{\beta}$ 	is a solution of \eqref{4.be.linear} fulfilling $g_\beta(x,p,0)=\Re (g_{\alpha(\beta),\beta}(x,p))$ for all $\beta\in I\setminus\{0\}$. The remaining parts are a direct consequence of Lemma \ref{lem.Def.of.g.cxi}.
	%	for all $(x,p)\in\R^d\times\T^d$ and all $\beta\in I$.
	%	is a solution of \eqref{4.be.linear} with $g(x,p,0)=g_{c,\varphi}(x,p)$ and $\gamma=0$. Let $N\in\mathbb N$. There exist $C,C_N>0$ such that
	%	\begin{align*}
	%	C \velo_1(p)^2\leq \N{g_{c,\varphi}(x,p)}\leq \NN{g_{c,\varphi}}_{W^{N,\infty}(\R^d\times\T^d)}\leq C_N(1+\N{\varphi}^N)
	%	\end{align*}
	%	for all $(x,p)\in\R^d\times\T^d$ and $\varphi\in\R$.
	%%	\begin{align*}
	%%%		\leq \int_{\T^d}\N{\frac{\velo_1(p)}{\velo_1(p)-ic}}F_{\bar{\lambda}}(p)(1-\eta F_{\bar \lambda}(p))dp\leq 
	%%		%\int_{\T^d}\frac{\velo_1(p)^2}{\velo_1(p)^2+c^2}F_{\bar{\lambda}}(p)(1-\eta F_{\bar \lambda}(p))dp.
	%%	\NN{g_{c,\varphi}(x,p)}_{L^1_p(\T^d)}\leq \frac{1}{4\eta}	\quad\mbox{and}\quad  \NN{g_{c,\varphi}(0,p)}_{L^1_p(\T^d)}=
	%%		\frac{1}{U\bar \lambda_1}.
	%%	\end{align*}
	%	In addition, there exists a $\nu_0>0$ and such that
	%	\begin{align}\label{eq.lem.g.c.varphi}
	%		\N{g_{c,\varphi}}_{C^\nu_x}\leq C_{\nu_0} e^{\nu\N\varphi}\quad\mbox{for all }\varphi \in \R\mbox{ and all }x\in \R^d
	%	\end{align}
	%	and for all $\nu\leq \nu_0$, where $C_{\nu_0}>0$ is a constant independent from $\varphi$ and $x$.
\end{proof}

\subsection{Nonlinear equation}
Fix $\bar{\lambda}$ and $\alpha_0$ such that \eqref{eq.Bed.for.lambda,c} is fulfilled (see Lemma \ref{lem.Bed.for.lambda,c}). We now choose $\nu>0$ such that $\NN{F_{\bar \lambda}}_{C^{\nu}}<\infty$. Let $g_\beta$ be as in Proposition \ref{lem.g.c.varphi} and let $c>0$ be given such that \eqref{eq.lem.g.c.varphi} is fulfilled. We set

\begin{align*}
%g^\delta(x,p) =  g_{c,\log(1/\delta)/\nu }(x,p)\quad\mbox{and}\quad 
f_0^\beta(x,p)=F_{\bar \lambda}(p)+\beta e^{-\frac{c\nu}{\N{\beta}}} g_{\beta}(x,p,0).
\end{align*}
Then $\|f_0^\beta\|_{C^\nu}$ is uniformly bounded w.r.t.~$\beta>0$. Since $F_{\bar \lambda}(p)\in [b,\eta^{-1}-b]$ for some $b>0$ and all $p\in\T^d$ and $g_{\beta}(x,p,0)$ is uniformly bounded w.r.t.~$x,p$ and $\beta$, we can apply Theorem \ref{thm: boltzmann local: hT solution}: there exists a $\beta_0>0$ and a $T>0$ such that

\begin{equation*}%\label{3.be}
\pa_t f^\beta + u(p)\cdot\na_x f^\beta - U\na_x n_{f^\beta}\cdot\na_p f^\beta = \gamma n_{f^\beta}(1-\eta n_{f^\beta})(\F_{f^\beta}-f^\beta)
\end{equation*}
has a unique analytic solution $f^\beta:\R^d\times\T^d\times[0,T)\to \R$ for each $\beta\in(-\beta_0,\beta_0)\setminus\{0\}$ with 

\begin{align*}
f^\beta(x,p,0)=f_0^\beta(x,p)%=F_{\bar \lambda}(p)+e^{-\frac{c\nu}{\N{\beta}}} g_{c,\varphi}(x,p).
\end{align*}
By shrinking $T>0$, the theorem moreover implies that there exist $\tilde{\nu}\in(0,\nu)$ and $\tilde C>0$ such that $\NN{f^\beta(t)}_{C^{\tilde \nu}}\leq \tilde C$ for all $\beta\in(-\beta_0,\beta_0)\setminus\{0\}$ and $t\in[0,T)$.
Define $h^\beta$ by the equation

\begin{equation*}
f^\beta(x,p,t)=F_{\bar \lambda}(p)+\beta e^{-\frac{c\nu}{\N{\beta}} }\left(g_{\beta}(x,p,t)+h^\beta(x,p,t)\right).
\end{equation*}
Then $h^\beta$ solves 

\begin{multline*}
\del_t h^\beta +\velo(p)\cdot\na_x h^\beta-U\na_p f^\beta\cdot\na_x n_{h^\beta}-U\beta e^{-\frac{c\nu}{\N{\beta}}}\na_x n_{g_{\beta}}\cdot \na_ph^\beta\\=Q^\beta+ U\beta e^{-\frac{c\nu}{\N{\beta}}} \na_x n_{g_{c,\varphi}}\cdot \na_pg_{c,\varphi}
\end{multline*}
with 
\[
Q^\beta:=\gamma n_{f^\beta}(1-\eta n_{f^\beta})\frac{\F_{f^\beta}-f^\beta}\beta e^{\frac{c\nu}{\N{\beta}}}-\gamma n_{F_{\bar{\lambda}}}(1-\eta n_{F_{\bar{\lambda}}})\left(G(\bar{\lambda})\cdot\binom{n_{g_{\beta}}}{E_{g_\beta}}-g_\beta\right)
\]
and $h(x,p,0)=0$. Note that $c$ is the constant provided by Proposition \ref{lem.g.c.varphi}.

\begin{lemma}\label{lem.h.is.small}
	There exist $C,\tau>0$ such that
	
	\begin{align*}
	\NN{h^\beta(\cdot,\cdot,t)}_{L^\infty(\R^d\times\T^d)}\leq Ct\quad\mbox{for }0\leq t\leq \tau\mbox{ and all }\N\beta<\beta_0.
	\end{align*}
\end{lemma}
\begin{proof}
	Recall the norms
	
	\begin{align*}
	\N{f}_{C^{\nu}}&:= \sum_{a,b\in\mathbb N_0^d} \frac{\nu^{|a+b|}}{a!b!}\|\del_x^a\del_p^b f\|_{L^{\infty}_xL^{1}_p}
	\quad\mbox{and}\quad
	\N{\velo}_{C^{\nu,\infty}}:= \sum_{b\in\mathbb N_0^d}\frac{\nu^{|b|}}{b!}\NN{\del^{b} \velo}_{L^{\infty}(\T^d)}
	\end{align*} 
	% First of all, we recall that $\N{f^\varphi(\cdot,\cdot,0)}_{C^{\tilde \nu}_x}\leq \tilde C$ for all $\phi\geq \phi_0$, $x\in\R^d$ and $t\in[0,T)$. According to Lemma \ref{lem: nnorm leq C norm 1}, for all $\mu<\tilde\nu$, there exists a $C_1\geq \tilde C$ such that
	% \[
	% \N{\del_pf^\varphi(\cdot,\cdot,0)}_{C^{\tilde \nu}_x}\leq \frac{C_1}{\tilde C}\N{f^\varphi(\cdot,\cdot,0)}_{C^{\tilde \nu}_x}\leq C_1
	% \]
	% for all $\phi\geq \phi_0$, $x\in\R^d$ and $t\in[0,T)$. Redoing this argument iteratively, for all $\mu<\tilde{\nu}$ and all $N\in\mathbb N$, there exists a constant $C_N>0$ such that
	%  \[
	%  \sum_{i=0}^{N}\N{\del_p^if^\varphi(\cdot,\cdot,0)}_{C^{\tilde \nu}_x}\leq C_N
	%  \]
	%  for all $\phi\geq \phi_0$, $x\in\R^d$ and $t\in[0,T)$. Using the Sobolev imbedding $C^0(\T^d)\subset W^{N,1}(\T^d)$ for large $N$, we have that
	%  \[\sup_{\varphi\geq\varphi_0}\sup_{t\in[0,T)}\sup_{x\in\R^d}\N{f}_{C^{\mu,\infty}}<\infty\]
	%  for all $\mu<\tilde \nu$. 
	and let $\mu\in(0,\nu/2)$ and $M>0$. 
	%  such that
	%\[...\]
	Similar to the proof of Lemma \ref{lem: d_t nnormClambda  f estimate}, we see that
	
	\begin{align*}
	\del_t \N{h^\beta}_{C^{\mu-Mt}}&\leq\N{\del_th^\beta}_{C^{\mu-Mt}} - M\N{\del_xh^\beta}_{C^{\mu-Mt}}-M\N{\del_ph^\beta}_{C^{\mu-Mt}}
	\\&\leq \N{Q^\beta}_{C^{\mu-Mt}}+\left(\N{\velo}_{C^{\mu-Mt,\infty}} +\N{U\na_p f^\beta}_{C^\mu} -M\right) \N{\del_xh^\beta}_{C^{\mu-Mt}}
	\\&\hspace{2cm}+\left( e^{-\frac{c\nu}{\N{\beta}}}\N{\beta U\na_x n_{g_{\beta}}}_{C^\mu}-M\right)\N{\del_ph}_{C^{\mu-Mt}}
	\\&\hspace{2cm} + e^{-\frac{c\nu}{\N{\beta}}} \N {\beta U\na_x n_{g_{\beta}}}_{C^{\mu-Mt}} \N{ \na_pg_{\beta}}_{C^{\mu-Mt}}.
	\end{align*}
	Using Proposition \ref{lem.g.c.varphi}, we note that there exists a constant $C_0>0$ independent from $\beta$ such that
	
	\begin{align*}
	\N{\beta U\na_x n_{g_{\beta}}}_{C^\mu}\leq C_0 \exp\left({\frac{\beta\omega(\beta)t+c\mu}{\N\beta}}\right).
	\end{align*}
	Thus, for $t\leq\tau:= \min\{c(\nu/2-\mu)/(\max_{\N\beta\leq \beta_0}\beta\omega(\beta)),T\}$ we have that
	
	\begin{align*}
	\N{U\na_x n_{g_{\beta}}}_{C^\mu}\leq C_0 \exp\left(\frac{c\nu}{2\N{\beta}}\right).
	\end{align*}
	Similarly to the proof of Lemma \ref{lem.Def.of.g.cxi} Proposition \ref{lem.g.c.varphi}, we can show that 	$\N{\na_p g_{\beta}}_{C^\mu}\leq C_{1} e^{\frac{c\mu}{\N\beta}}$ for some $C_1>0$ which does not depend on $\beta$. Choosing now 
	\[
	M:=C_0+\sup_{\N{\beta}<\beta_0}\sup_{t\in[0,\tau)}\left(\N{U\na_p f^\varphi}_{C^\mu}+\N{\velo}_{C^{\mu,\infty}}\right)<\infty.
	\]
	We note that $M$ is finite due to the choice of $\tilde{\tau}$ and Lemma \ref{lem: nnorm leq C norm 2}, because $ \NN{f^\varphi}_{C^\nu}$ is uniformly bounded and $\mu<\nu$. This choice of $M$ implies that
	
	\begin{align*}
	\del_t \N{h^\beta}_{C^{\mu-Mt}_x}&\leq \N{Q^\beta}_{C^{\mu-Mt}_x}+C_0C_1
	\end{align*}
	for $0\leq t\leq \min\{\tau,\mu/2M\}$. In order to show that the first term on the r.h.s.~is also bounded for small $t$, we define $H_s:=F_{\bar \lambda}(p)+s\left(g_{\beta}(x,p,t)+h^\beta(x,p,t)\right)$ and \[\phi:s\mapsto \gamma n_{H_s}(1-\eta n_{H_s})(\F_{H_s}-H_s)/s.
	\]
	Then, we have
	
	\begin{align*}
	Q^\beta=\phi(\beta^{-1}e^{-\frac{c\nu}{\N{\beta}}})-\lim_{s\to0}\phi(s^{-1}e^{-\frac{c\nu}{s}})=\int_0^{\N\beta}\phi'(s^{-1}e^{-\frac{c\nu}{s}})e^{-\frac{c\nu}{s}}(c\nu-s)\frac{ds}{s^3}
	\end{align*}
	with
	
	\begin{align*}
	\phi'(s)&=\gamma (1-2\eta n_{H_s})\partial_s n_{H_s}\frac{\F_{H_s}-H_s}{s}
	\\&\hspace{2cm}-\gamma n_{H_s} (1-\eta n_{H_s})\frac{\F_{H_s}-s\partial_{H_s}\F_{H_s}\partial_s H_s-(H_s-s\partial_s H_s)}{s^2}
	\\&=
	\gamma (1-2\eta n_{H_s})(n_{g_\beta}+n_{h^\beta})\frac{\F_{H_s}-H_s}{s}
	\\&\hspace{2cm}-\gamma n_{H_s} (1-\eta n_{H_s})\frac{\F_{H_s}-\F_{H_0}-s\partial_{s}\F_{H_s}}{s^2}
	\end{align*}
	Since $f\mapsto \F_f$ is twice differentiable (see appendix) and $H_s$ is linear in $s$, one can prove that $\N{\phi'(s)}_{C^{\mu-Mt}}$ is uniformly bounded for small $s>0$. Thus,
	
	\begin{align*}
	\N{Q^\beta}_{C^{\mu-Mt}}\leq 
	%	\sup_{0\leq s\leq \N{\beta}}
	%		 \N{\phi'(s)}_{C^{\mu-Mt}}\int_0^{|\beta|}e^{-\frac{c\nu}{s}}(c\nu-s)\frac{ds}{s^3}=
	\sup_{0\leq s\leq \N{\beta}}\N{\phi'(s)}_{C^{\mu-Mt}}\N\beta^{-1}e^{-\frac{c\nu}{\N{\beta}}}
	\leq C_2
	\end{align*}
	for $0<\beta\leq \beta_0$ some $C_2>0$ depending only $\beta_0$.	 Therefore, 
	\[
	|h^\beta|_{C^{\mu-Mt}}\leq (C_2+C_0C_1)t
	\]
	for $t\leq \min\{\tau,\mu/2M\}$. Finally, we can use $\NN{\cdot}_{L^\infty(\R^d\times\T^d)}\leq C\N{\cdot}_{C^{\mu-Mt}}$ for some $C>0$ and all $0\leq t\leq \mu/2M$ in order to finish the proof.
\end{proof}

\begin{remark}
	For every $\delta>0$ there exists a constant $C_\delta>0$ such that 
	\[
	\NN{g_{\beta}(\cdot,\cdot,t)}_{L^1(B_\delta(p,x))}\geq 2C_\delta e^{\omega(\beta)t}\quad\mbox{for all }\N{\beta}\leq \beta_0
	\]
	and $(x,p)\in \R^d\times\T^d$ and small $t$. 
	Thus, by Lemma \ref{lem.h.is.small}, there exists a $\tau_\delta>0$ such that 
	
	\begin{align}\label{eq.f.beta-F.geq.C}
	\NN{f^\beta(\cdot,t)-F_{\bar{\lambda}}}_{L^1(B_\delta(x,p))}\geq C_\delta e^{\omega(\beta)t-\frac{c\nu}{\N{\beta}}}
	\end{align}
	for all $t<\tau_\delta$, $x\in \R^d$, $p\in \T^d$ and $\N{\beta}<\beta_0$, where $ \omega$ is given by Proposition \ref{lem.g.c.varphi} satisfying  
	\[\beta\omega(\beta)\geq \tilde c\quad\mbox{for some }\tilde c>0\mbox{ and all }|\beta|<\beta_0.\] 
	%	(\beta):=\gamma n_{F_{\bar{\lambda}}}(1-\eta n_{F_{\bar{\lambda}}})\frac{\alpha(\beta)^2-\beta}{\beta}$.
\end{remark}
\begin{proof}
	The first part is clear due to the definition of $g_\beta$. The second assertion is then a consequence of
	Lemma \ref{lem.h.is.small}, which guarantees  for sufficiently small $t>0$ that \[\NN{h^{\beta}(\cdot,\cdot,0)}_{L^1(B_\delta(p,x))}\leq C t.\qedhere\]%_\delta e^{\omega(\beta)t}$.
\end{proof}
\begin{proof}[Proof of Theorem \ref{thm.illposed}]
	Let $\theta>0$, $\delta>0$ and $k\in\mathbb N_0$. If we combine \eqref{eq.lem.g.c.varphi2} with \eqref{eq.f.beta-F.geq.C}, we see that there exists a constant $C_{\delta,k,\nu}>0$ such that
	
	\begin{multline*}
	\NN{f^\beta(\cdot,t)-F_{\bar{\lambda}}}_{L^1(B_\delta(x,p))}\\\geq \N{\beta}^{-1}C_{\delta,k,\nu}\left(\inf_{\N\beta\leq\beta_0} \frac{e^{\omega(\beta)t-\frac{c\nu}{\N{\beta}}}}{\N{\beta}(1+\N{\beta}^{-k})^\theta}\right)\NN{f^\beta(\cdot,0)-F_{\bar{\lambda}}}_{W^{k,\infty}(\R^d\times\T^d)}^\theta
	\end{multline*}
	We recall $\omega$ from Proposition \ref{lem.g.c.varphi}
	%	the definition of $\omega(\beta)=\gamma n_{F_{\bar{\lambda}}}(1-\eta n_{F_{\bar{\lambda}}})\frac{\alpha(\beta)^2-\beta}{\beta}$
	and see that 
	\[\inf_{\N\beta\leq\beta_0} \frac{e^{\omega(\beta)t-\frac{c\nu}{\N{\beta}}}}{\N{\beta}(1+\N{\beta}^{-k})^\theta}>0
	%	\]
	\qquad
	\mbox{if} 
	%	\[
	\qquad
	t>\tau_{\min}(\nu):=\frac{c\nu}{\inf_{|\beta|<\beta_0}\beta\omega(\beta)}
	%	\frac{c\nu}{\gamma n_{F_{\bar{\lambda}}}(1-\eta n_{F_{\bar{\lambda}}})(\inf_{|\beta|<\beta_0}\alpha(\beta)^2-\beta)}
	\]
	assuming that $\beta_0$ is sufficiently small such that $\beta\omega(\beta)$ is positive for all $\N\beta\leq \beta_0$. Since the parameter $\nu>0$ was arbitrary, we may choose $\tau_{\min}(\nu)<\delta/2$. Therefore, we just have proved that for any $\delta>0$ and $k\in\mathbb N$, there exist a $C_{\delta,k,\theta}>0$ and a $\tau>\delta$ such that
	
	\begin{align*}
	\NN{f^\beta(\cdot,\cdot,t)-F_{\bar{\lambda}}}_{L^1(B_\delta(x,p))}\geq \frac{C_{\delta,k,\theta}}{\N{\beta}}\NN{f^\beta(\cdot,\cdot,0)-F_{\bar{\lambda}}}_{W^{k,\infty}(\R^d\times\T^d)}^\theta\quad\mbox{for all }\N{\beta}\leq\beta_0
	\end{align*}
	and for all $x\in\R^d$, $p\in\T^d$ and $t\in(\delta,\tau)$. This implies the assertion of the theorem as $\beta\to0$.
\end{proof}

\section{Space local method}\label{sec.space.local}
In order to improve the existence results we have obtained so far, we need to make use of the fact that the collision operator of the semiconductor-Boltzmann-Dirac-Benny equation is local in space. Therefore, we are now focusing on a space local version of the method presented in sections \ref{analytic-norms} and \ref{local}. For this we replace the analytic norms $\N\cdot_{C^\nu}$ to space-local semi-norms, i.e.\ we define for every point $x$ in the physical space a semi-norm $\NN{f}_{C^{\nu}_x}$ that only consists of all the derivatives of $f$ evaluated at the point $x$.
\begin{definition}
	Let $\nu>0$, $d\in\mathbb N$ and fix $x\in \R^d$. We consider the space-local semi-norms
	
	\begin{equation*}
	\NN{f}_{C^{\nu}_x}:= \sum_{\substack{i,j\in\mathbb N_0^d\\|i+j|\leq 1}}\sum_{a,b\in\mathbb N_0^d} \frac{\nu^{|a+b|}}{a!b!}\int_{\R^d}|\del_x^{a+i}\del_p^{b+j} f(x,p)|dp
	\end{equation*} and
	
	\begin{equation*}
	\NN{Df}_{C^{\nu}_x}:= \sum_{\substack{a,b\in\mathbb N_0^d\\|a+b|= 1}}\NN{\del_x^a\del_p^af}_{C^{\nu}_x}
	\end{equation*}
	as well as
	
	\begin{equation*}
	%		\N{f}_{\dot C^{\nu}_x}:= \N{f}_{C^{\nu}_x}-\N{f}_{C^{0}_x}\quad\mbox{and}\quad
	\NN{f}_{\dot C^{\nu}_x}:= \NN{f}_{C^{\nu}_x}-\NN{f}_{C^{0}_x}.
	\end{equation*} 
	for $f:\R^d\times\T^d\to \R^k$ being analytic. 
	
	Let $\nu,T>0$, $\mu\in[0,\nu/T)$. Using the semi-norms from above, we define 
	
	\begin{equation*}
	\NNN{f}_{\nu,\mu}:=\sup_{x\in\R^d}\NNN{f}_{\nu,\mu,x},\ \NNN{f}_{\nu,\mu,x}:=\sup_{0\leq t<T}\left(\!\NN{f(t)}_{C^{\nu-\mu t}_x}+\mu\int_0^t\NN{Df(s)}_{C^{\nu-\mu s}_x}ds\!\right)
	\end{equation*}
	for $f:\R^d\times\T^d\times[0,T)\to\R$ being analytic in $(x,p)$ and continuous in $t$ writing $f(t)=f(\cdot,\cdot,t)$.	
\end{definition}
Note that we can prove the following version of Lemma \ref{lem: nnorm leq C norm 1} for these semi-norms. The proof is similar to that of Lemma \ref{lem: nnorm leq C norm 1} and will be omitted.
\begin{lemma}\label{lem: nnorm leq C norm 2}
	Let $\mu_2>\mu_1>0$ and $d\in\mathbb N$. Then there exists a constant $C=C_{\mu_1,\mu_2}>0$ such that for all analytic $f:\R^d\times\T^d\to \R$ and all $x\in\R^d$, it holds
	\[
	\NN{f}_{\dot C^{\nu}_x}\leq \nu C_{\mu_1,\mu_2}\N{f}_{\dot C^{\mu_2}_x}
	\]
	for all $\nu\in[0,\mu_1]$, where
	\[
	\N{f}_{\dot C^{\nu}_x}:=\sum_{0\neq(a,b)\in\mathbb N_0^{2d}} \frac{\nu^{|a+b|}}{a!b!}\int_{\R^d}|\del_x^{a}\del_p^{b} f(x,p)|dp.
	\]
\end{lemma}

With the same arguments as in the previous section, one can prove the following counterpart to Theorem \ref{thm: boltzmann local: hT solution}.
\begin{theorem}\label{thm: boltzmann local: hT solution2}
	Let $C,R,\nu>0$ and $f_0:\R^d\times\T^d\to \R$ be analytic such that
	\begin{equation}\label{eins}
	\sup_{x\in\R^d}\NN{f_0}_{C^{\nu}_x}<R.
	\end{equation}
	Then if $\mu>0$ is sufficiently large, $T\in (0,\nu/\mu)$ and $F:\R^d\times\T^d\times[0,T)\to \R$ is analytic such that
	\begin{equation}\label{zwei}
	\NN{F(t)}_{C^{\nu-\mu t}_x}\leq C
	\end{equation}
	for all $0\leq t\leq T$ and $x\in\R^d$, then the equation 
	\begin{equation}\label{2.be.nochmal}
	\del_t f+\velo\cdot\nabla_x f - U\nabla_xn_f\cdot\nabla_p f =\gamma n_f(1-\eta n_f)(F-f)
	\end{equation}
	admits a unique analytic solution $f:\R^d\times\T^d\times[0,T)\to\R$ with $\NNN{f}_{\nu,\mu,x}\leq R$ and $f(x,p,0)= f_0(x,p)$ for all $x\in\R^d$ and $p\in \T^d$. 
	
	Moreover, let $\Psi: (f_0,F)\mapsto f$ be defined by the unique solution of \eqref{2.be.nochmal}
	with $f(x,p,0)=f_0(x,p)$. If $\mu>0$ is sufficiently large, the mapping $\Psi$ is Lipschitz continuous, i.e., for all $x\in \R^d$
	
	\begin{equation*}
	\NNN{\Psi(f_0,F)-\Psi(g_0,G)}_{\nu,\mu,x}\leq 2\NNN{(f_0,F)-(g_0,G)}_x,
	\end{equation*}
	where
	\begin{align*}
	\NNN{(f_0,F)}_x&:= \NN{f_0}_{C^\nu_x}+\mu^{-\frac12}\sup_{0\leq t<T}\NN{F}_{C^{\nu-\mu t}_x}.
	\end{align*} 
	for $f_0,g_0$ and $F,G$ satisfying \eqref{eins} and \eqref{zwei}, respectively.
\end{theorem}
Similarly as in estimate \eqref{eq.strong.Lipschitz} in the proof of Theorem \ref{thm: boltzmann local: hT solution}, we can improve the Lipschitz estimate.
\begin{lemma} Let $f:=\Psi(f_0,F)$ and $g:=\Psi(g_0,G)$. We have
	\begin{equation}\label{eq.strong.Lipschitz2}
	\NNN{f-g}_{\nu,\mu,x}\leq
	2\NN{f_0-g_0}_{C^\nu_x} +\frac{4\gamma\nu}{\mu}\sup_{0\leq t\leq T}\NN{n_f(1-\eta n_f)}_{C^\nu_x}\NN{F-G}_{C^{\nu-\mu t}_x}.
	\end{equation}
	if $\mu>0$ is sufficiently large.
\end{lemma}

\section{BGK-type collision operator - space local method}\label{sec.BGK.space.local}
In this section, we consider again equation

\begin{equation}\label{3.be.nochmal}\tag{\ref{3.be}}
\pa_t f + u(p)\cdot\na_x f - U\na_x n_f\cdot\na_p f = \gamma n_f(1-\eta n_f)(\F_f-f)
\end{equation}
with $f(x,p,0)=f_0(x,p)$ for given $U\neq0$ and $\gamma\geq 0$.
% \eqref{3.be.again}, namely
%
%\begin{equation}\label{3.be.nochmal}
%\pa_t f + u(p)\cdot\na_x f - U\na_x n_f\cdot\na_p f = \gamma n_f(1-\eta n_f)(\F_f-f)
%\end{equation}
%for given $U\neq0$ and $\gamma,\eta\geq0$. 
As before, we use the self-consistent equilibrium distribution function \[ \F_f(x,p,t) = \big(\eta + \exp(-\lambda_0(x,t)-\lambda_1(x,t)\epsilon(p))\big)^{-1}\ \mbox{for }x\in\R^d,\ p\in\T^d\mbox{ and }t>0,\] where  $\lambda_0,\lambda_1$ satisfy

\begin{equation}
n_f(x,t):=n_{\F_f}(x,t)\quad\mbox{and}\quad E_f(x,t)=E_{\F_f}(x,t),
\end{equation}
for $n_f(x,t):=\int_{\T^d}f(x,p,t)dp$ and $E_f(x,t)=\int_{\T^d}\epsilon(p)f(x,p,t)dp$. 

The main goal is to improve the existence result from Theorem \ref{thm.eins} using the space local semi-norms. Similar as before, the key ingredient will Theorem \ref{thm: boltzmann local: hT solution2} and the Lipschitz estimate \eqref{eq.strong.Lipschitz2}.

\begin{definition}
	%		Let $\nu>0$ and fix $x\in \R^d$. We consider the space-local semi-norms
	
	%	\begin{equation*}
	%\N{f}_{\dot C^{\nu}_x}:= \N{f}_{C^{\nu}_x}-\N{f}_{C^{0}_x}	\quad\mbox{and}\quad\NN{f}_{\dot C^{\nu}_x}:= \NN{f}_{C^{\nu}_x}-\NN{f}_{C^{0}_x}.
	%	\end{equation*} 
	Let $a\geq 1$ and $\delta>0$. We define
	%	\begin{align*}
	%	Z_a:=\left\{p\mapsto\frac{1}{\eta+e^{-\lambda_0-\lambda_1\epsilon(p)}}: \lambda_0,\lambda_1\in\R\mbox{ with } \N{\lambda_1}\leq \log a\right\}
	%	\end{align*}
	%	and
	%	\begin{align*}
	%	M_a:=\left\{\int_{\T^d}(1,\epsilon(p))f(p)dp: f\in Z_a \right\}.
	%	\end{align*}
	
	\begin{align*}
	M_a:=\left\{\int_{\T^d}\frac{(1,\epsilon(p))dp}{\eta+e^{-\lambda_0-\lambda_1\epsilon(p)}}: \lambda_0,\lambda_1\in\R\mbox{ with } \N{\lambda_1}\leq \log a\right\}\subset \mathbb R^2
	\end{align*}
	and
	
	\begin{align*}
	\mathcal U_{a,\delta}
	&:=
	%	\left\{y\in \R^2: \exists (m_0,m_1)\in M_a \mbox{ with }|(m_0,m_1)-y|<\delta m_0(1-\eta m_0)\right\}
	%	\\&=
	\bigcup_{(m_0,m_1)\in M_a}B_{\delta m_0(1-\eta m_0)}(m_0,m_1)\supset M_a,
	\end{align*}
	where $B_\theta(y)$ denotes the ball in $\mathbb R^2$ centered at $y$ with radius $\theta$.
\end{definition}
\begin{proposition}\label{prop.BGK.contraction}
	Let $\eta,\nu_0,R>0$, $\gamma\geq0$, $a\geq1$. Then there exist $\alpha,\beta,\mu>0$ such that the following holds:
	
	Let $f_0:\R^d\times\T^d\to \R$ is analytic such that $\NN{f_0}_{C^\nu_x}\leq R/2$ for some $\nu\in(0,\nu_0)$ and \[\binom{n_{f_0}(x)}{E_{f_0}(x)}:=\int_{\T^d}\binom{1}{\epsilon(p)}f_0(x,p)dp\in \mathcal U_{a,\alpha/2}\]  is well-defined for all $x\in\R^d$. 
	%			\begin{equation*}
	%			\NN{n_{f_0}}_{\dot C^\nu_x}+\NN{E_{f_0}}_{\dot C^\nu_x}\leq \frac\delta2 n_{f_0}(x)(1-\eta n_{f_0}(x))%\quad\mbox{for }x\in\R^d
	%			\end{equation*}
	%			for all $x\in\R^d$. 
	Moreover, suppose that
	
	\begin{equation*}
	\NN{{f_0}}_{\dot C^{\nu}_x}\leq \beta n_{f_0}(x)(1-\eta n_{f_0}(x))\quad\mbox{for }x\in\R^d.
	\end{equation*}
	Then equation \eqref{3.be.nochmal}	with $f\big|_{t=0}=f_0$ admits an analytic solution $f:\R^d\times\T^d\times[0,T)\to\R$ with $\NNN{f}_{\nu,\mu}\leq R$ for $T<\nu/\mu$.
\end{proposition}
The theorem will also be proved using the Banach fixed-point theorem. In order to define the right metric space, we require some properties of the equilibrium distribution.

\begin{proposition}\label{prop.Ff.is.Lipschitz}
	Let $\eta>0$, $a\geq1$ and $R,\nu>0$. Then there exist $\alpha>0$ such that for all $f,g:\R^d\times\T^d\to\R$ being analytic with $\NN{f}_{C^\nu_x},\NN{g}_{C^\nu_x}\leq R$ and
	
	\begin{equation*}
	\NN{n_h}_{\dot C^\nu_x}+\NN{E_h}_{\dot C^\nu_x}\leq \alpha n_h(x)(1-\eta n_h(x))\quad\mbox{for }h\in\{f,g\}
	\end{equation*}
	and $(n_h(x),E_h(x))\in \mathcal U_{a,\alpha}$ for $h\in\{f,g\}$, it holds 
	
	\begin{equation*}
	\NN{\F_f}_{C^\nu_x},\NN{\F_g}_{C^\nu_x}\leq C
	\end{equation*}
	and
	
	\begin{equation*}
	\NN{\F_f-\F_g}_{C^\nu_x}\leq C\NN{f-g}_{C^\nu_x}.
	\end{equation*}
	for some $C>0$ and all $x\in \R^d$.
\end{proposition}
\begin{proof}
	See appendix.
\end{proof}
\begin{remark}\label{rem.alpha.depends.only.on.a}
	According to the proof in the appendix, the parameter $\alpha$ only depends on $a$. More precisely, it can be written as $\alpha=1/(2B_a)$ for  $B_a$ from Lemma \ref{Ableitung.FrhoE}.
\end{remark}
\begin{definition}\label{def.Z^t}
	For $R,\nu,\eta>0$, $a\geq1$ let $\alpha>0$ be as in Proposition \ref{prop.Ff.is.Lipschitz}. 
	%	which only depends on $a$ according to Remark \ref{rem.alpha.depends.only.on.a}. Given $R,\nu>0$,
	% let $\delta>0$ be as in Proposition \ref{prop.Ff.is.Lipschitz}. 
	Moreover, let $\mu>0$ and $T\in(0,\nu/\mu)$. We assume that
	
	\begin{align*}
	\NN{f_0}_{\dot C^{\nu}_x}\leq \beta n_{f_0}(x)(1-\eta n_{f_0}(x))
	\end{align*}
	\\
	%	We define $X$ as the complete metric space of all analytic functions $f:\R^d\times\T^d\times[0,T)\to \R$ with $\NNN{f}_{\nu,\mu}\leq R$. 
	%%	Note that this is a complete metric space if the metric is induced by the norm $\NN{\cdot}_{\nu,\mu}$.
	Let $Z$ space of all analytic functions $f:\R^d\times\T^d\times[0,T)\to [0,\eta^{-1}]$ satisfying
	\begin{enumerate}
		\item $\NNN{f}_{\nu,\mu}\leq R$,
		\item
		$\NN{n_{f}(\cdot,t)}_{\dot C^{\nu-\mu t}_x}+\NN{E_{f}(\cdot,t)}_{\dot C^{\nu-\mu t}_x}\leq \alpha n_{f}(x,t)(1-\eta n_{f}(x,s))$ and
		
		\item $(n_f(x,t),E_f(x,t))\in \overline{\mathcal U_{a,\alpha}}$ 
	\end{enumerate}
	for all $x\in\R^d$ and $t\in[0,T)$. 
	% $\NNN{f}_{\nu,\mu}\leq R$ such that the mapping	$(n_f,E_f):\R^d\times[0,t)\to \overline{\mathcal U_{a,\alpha}}$ defined by \[(n_f,E_f):=\int_{\T^d}(1,\epsilon(p))f(\cdot,p,\cdot)dp\] is well-defined and fulfills	
	%	\[\NN{n_{f}(\cdot,s)}_{\dot C^{\nu-\mu s}_x}+\NN{E_{f}(\cdot,s)}_{\dot C^{\nu-\mu s}_x}\leq \delta n_{f}(x,s)(1-\eta n_{f}(x,s))\quad\mbox{for }x\in\R^d,s\in[0,t).\]
	Thus, $Z$ is a complete if the metric is induced by the norm $\NN{\cdot}_{\nu,\mu}$. 
	
\end{definition}	
Let $f_0:\R^d\times\T^d\to \R$ be analytic such that $\NN{f_0}_{C^\nu_x}\leq R/2$ and \[(n_{f_0}(x),E_{f_0}(x)):=\int_{\T^d}(1,(\epsilon(p)))f_0(x,p)dp\in \mathcal U_{a,\alpha/2}\]  is well-defined for all $x\in\R^d$. 
%			\begin{equation*}
%			\NN{n_{f_0}}_{\dot C^\nu_x}+\NN{E_{f_0}}_{\dot C^\nu_x}\leq \frac\delta2 n_{f_0}(x)(1-\eta n_{f_0}(x))%\quad\mbox{for }x\in\R^d
%			\end{equation*}
%			for all $x\in\R^d$. 
Moreover, suppose that

\begin{equation*}
\NN{{f_0}}_{\dot C^{\nu}_x}\leq \beta n_{f_0}(x)(1-\eta n_{f_0}(x))\quad\mbox{for }x\in\R^d
\end{equation*}
for some small $\beta>0$. For sufficiently large $\mu>0$, we define
the mapping 
\[
\Theta:Z\ni g\mapsto f,
\]
where $f$ is the solution of
\[\pa_t f + u(p)\cdot\na_x f - U\na_x n_f\cdot\na_p f = \gamma n_f(1-\eta n_f)(\F_g-f)\]
with $f\big|_{t=0}=f_0$. This is well-defined for large $\mu>0$ according to Theorem \ref{thm: boltzmann local: hT solution2} and Proposition \ref{prop.Ff.is.Lipschitz}.
As we plan to apply the Banach fixed-point theorem, we need to show that $\Theta$ is a contraction, i.e., the image of $\Theta$ is included in $Z$ and $\Theta$ is Lipschitz continuous with Lipschitz constant $L<1$. We start with the Lipschitz estimate, which is in this case the easier assertion.
\begin{lemma}
	Let $\mu>0$ be sufficiently large. Then for $f,g\in Z$ it holds
	\[\NNN{\Theta(f)-\Theta(g)}_{\nu,\mu}\leq \frac12\NNN{f-g}_{\nu,\mu}.\]
\end{lemma}
\begin{proof}
	Using $\Psi$ from Theorem \ref{thm: boltzmann local: hT solution2}, we can rewrite $\Theta$ as 
	\[\Theta(f)=\Psi(f_0,\F_f).\]
	For $f\in Z$, we know from Proposition \ref{prop.Ff.is.Lipschitz} that $\NN{\F_f}_{C^{\nu-\mu t}_x}\leq C$ for some $C>0$ and all $x\in \R^d$ and $t\in[0,T)$. Thus, Theorem \ref{thm: boltzmann local: hT solution2} entails that for sufficiently large $\mu>0$, 
	
	\begin{align*}
	\NNN{\Theta(f)-\Theta(g)}_{\nu,\mu}\leq 2\mu^{-\frac12}\sup_{t,x}\NN{\F_f-\F_g}_{C^{\nu-\mu t}_x}.
	\end{align*}
	Then the second statement of Proposition \ref{prop.Ff.is.Lipschitz} yields that 
	
	\begin{align*}
	\NNN{\Theta(f)-\Theta(g)}_{\nu,\mu}\leq C\mu^{-\frac12}\sup_{t,x}\NN{f-g}_{C^{\nu-\mu t}_x}\leq C\mu^{-\frac12}\NNN{f-g}_{\nu,\mu} 
	\end{align*}
	for some $C>0$. This implies the assertion for sufficiently large $\mu$ satisfying $\mu\geq 4C^2$.
\end{proof}
\begin{lemma}
	Let $\mu>0$ be sufficiently large, $(1+\nu^2)\beta>0$ sufficiently small and $g\in Z$. Then $\Theta(g)\in Z$.
\end{lemma}
\begin{proof}
	Let $g\in Z$  and define $f:=\Theta(g)$.
	\\
	{\bf Claim 1}: $\NNN{f}_{\nu,\mu}\leq R$ if $\mu$ is sufficiently large.
	\\
	This is a direct consequence of Theorem \ref{thm: boltzmann local: hT solution2} combined with Proposition \ref{prop.Ff.is.Lipschitz}.
	\\
	{\bf Claim 2}: We have	
	\[\NN{n_{f}(\cdot,t)}_{\dot C^{\nu-\mu t}_x}+\NN{E_{f}(\cdot,t)}_{\dot C^{\nu-\mu t}_x}\leq \alpha n_{f}(x,t)(1-\eta n_{f}(x,t))\quad\mbox{for }x\in\R^d,t\in[0,T).\]
	%			Since $(x,p,t)\mapsto f_0(x,p)\in Z^\tau$ and $\Theta(\F_{f_0})\neq f_0????$, we obtain similarly as before that
	%			\[\NN{f-f_0}_{X^\tau}\leq \tau L_{T'}C_{T'}\NN{g-f_0}_{X^\tau}.\]
	%			Using that $\NN{n_f-n_{f_0}}_{\dot C^{\nu-\mu t}_x}+\NN{E_f-E_{f_0}}_{\dot C^{\nu-\mu t}_x}\leq (1+\NN{\eps}_{L^\infty(\T^d)})\NN{f-{f_0}}_{X^\tau}$ for $t\in[0,\tau)$ and $x\in\R^d$, there exists a constant $C>0$ depending on $T'$ such that
	%			\begin{align*}
	%			\NN{n_f-n_{f_0}}_{\dot C^{\nu-\mu t}_x}+\NN{E_f-E_{f_0}}_{\dot C^{\nu-\mu t}_x}
	%			&\leq C\tau \NN{g-f_0}_{X^\tau}
	%			\\&\leq C\tau \left(\NN{g}_{X^\tau}+\sup_{x\in \R^d}\NN{f_0}_{C^\nu_x}\right)
	%			\end{align*}
	%		
	Fix $x\in \R^d$ and define
	\[h_0:\R^d\times\T^d\to \R, \ (y,p)\mapsto f_0(x,p)+\del_xf_0(x,p)y\]
	as well as
	
	\begin{align*}
	h(y,p,t)=h_1(p,t)+h_2(p,t)y,
	\end{align*}
	where $h_1,h_2$ solve
	
	\begin{align*}
	\del_t h_1+h_2\velo(p)=U\na_xn_{f_0}(x)\cdot\na_p h_1 \quad\mbox{and}\quad \del_t h_2=U\na_x n_{f_0}(x)\cdot\na_p h_2.
	\end{align*}
	with $h_1(p,0)=f_0(x,p)$ and $h_2(p,0)=\del_xf(x,p)$. Then it holds 
	
	\begin{align*}
	\del_t h+\velo\cdot\na_y h&= \del_t h_1+ \del_t h_2 y+h_2\velo
	%\velo\cdot\nabla_y h- U\nabla_yn_h\cdot\nabla_p h 
	\\&=U\na_xn_{f_0}(x)\cdot\na_p (h_1+h_2y)=U\nabla_yn_h\cdot\nabla_p h.
	\end{align*}
	Note that the equations for $h_1$ and $h_2$ are linear transport equation. We thus can solve them explicitly, e.g.
	\[h_2(p,t)=\del_x f_0(x,p-tU\na_xn_{f_0}(x)).\]
	With this, we can easily compute the density $n_{h_1}=\int_{\T^d}h_1(p,\cdot)dp$ by
	\[
	n_{h_1}(t)=n_{f_0}(x)-\int_0^t\int_{\T^d}\del_xf_0(x,p-sU\na_xn_{f_0}(x))\velo(p)dpds
	\]
	and estimate
	
	\begin{align}\label{eq.in.proof.n_h_1}
	\N{n_{h_1}(t)-n_{f_0}(x)}\leq t \NN{\velo}_{L^\infty(\T^d)}\NN{\del_x f^0(x,p)}_{L^1_p(\T^d)}.
	\end{align}
	
	Next, we infer from the Lipschitz estimate \eqref{eq.strong.Lipschitz2} that 		\begin{equation}
	\NNN{f-h}_{\nu,\mu,x}\leq
	2\NN{f_0-h_0}_{C^\nu_x} +\frac{4\gamma\nu}{\mu}\sup_{0\leq t\leq T}\NN{n_h(1-\eta n_h)}_{C^\nu_x}\NN{\F_g}_{C^{\nu-\mu t}_x}.
	\end{equation}
	for sufficiently large $\mu>0$. 
	%		\begin{align*}
	%		\NN{f-h}_{\nu,\mu,t,x}&\leq \frac{1}{1-\alpha}\NN{f_0-h_0}_{C^\nu_x}
	%		\\&\hspace{2cm} + t\frac{\gamma}{1-\alpha}\sup_{0\leq t\leq T}\NN{n_h(1-\eta n_h)}_{C^{\nu-\mu t}_x}\NN{\F_g-h}_{C^{\nu-\mu t}_x}.
	%		\end{align*}
	%		for $t\leq\tau$.
	At first, we note that $\NN{\F_g}_{C^{\nu-\mu t}_x}$ and $\NN{h}_{C^{\nu-\mu t}_x}$ are uniformly bounded. Then, we see by the definition of $h$ that we can estimate the r.h.s.~using that $\NN{f_0-h_0}_{C^\nu_x}\leq\NN{f_0}_{\dot C^\nu_x}$ and obtain
	
	\begin{align*}
	\NNN{f-h}_{\nu,\mu,x}&\leq 2\NN{f_0}_{\dot C^\nu_x}+ \frac {C\nu}\mu\sup_{0\leq t\leq T}\N{n_{h_1}(t)(1-\eta n_{h_1}(t))}.
	\end{align*}
	for some $C>0$ independent from $\nu$.  Moreover, it holds 
	
	\begin{align*}
	\sup_{0\leq t\leq T}\N{n_{h_1}(1-\eta n_{h_1})}&\leq \N{n_{f_0}(x)(1-\eta n_{f_0}(x))}+CT\NN{\del_x f^0(x,p)}_{L^1_p(\T^d)}
	\\&\leq
	\N{n_{f_0}(x)(1-\eta n_{f_0}(x))}+\frac{C\nu}{\mu}\NN{f^0}_{\dot C^\nu_x}
	\end{align*}
	because $T<\nu/\mu$.
	Thus, there exists a constant $C>0$ independent from $\nu$ such that for all $t\leq \tau_0$, we have
	
	\begin{align*}
	\NNN{f-h}_{\nu,\mu,x}&\leq \left(2+\frac{C\nu^2}{\mu^2}\right)\NN{f_0}_{\dot C^\nu_x}+ \frac {C\nu}\mu \N{n_{f_0}(x)(1-\eta n_{f_0}(x))}.
	\end{align*}
	%		By enlarging that constant, we obtain from Lemma \ref{lem: nnorm leq C norm 2} that  
	%		
	%		\begin{align*}
	%		\NNN{f-h}_{\nu,\mu,x}&\leq C\N{f_0}_{\dot C^{\nu_0}_x}+ tC\N{n_{f_0}(x)(1-\eta n_{f_0}(x))}.
	%		\end{align*}
	Note that $h$ is affine in $y$, hence $\del_y^i h=0$ for $|i|\geq2$ and
	
	\begin{align*}
	\sum_{\substack{|i|,|j|=0,1, a,b\in \mathbb N_0^d\\ |i+a|\geq2}}\frac{(\nu-\mu s)^{a+b}}{a!b!}\NN{\del_x^{i+a}\del_p^{j+b}f(x,p,t)}_{L^1_p(\T^d)}\leq \NNN{f-h}_{\nu,\mu,x}
	\end{align*}
	for $0\leq t< T$. In particular,
	
	\begin{multline*}
	\NN{n_f}_{\dot C^{\nu-\mu t}_x}+\NN{E_f}_{\dot C^{\nu-\mu t}_x}\leq (1+\NN{\epsilon}_{L^\infty(\T^d)})\NNN{f-h}_{\nu,\mu,x}\\+\nu\left(\N{\del_xn_f(x,t)}+\N{\del_xE_f(x,t)}\right).
	\end{multline*}
	Moreover, we can estimate the latter two terms by
	
	\begin{align*}
	\N{\del_xn_f(x,s)}&\leq \left[\N{\del_y(n_f(y,s)-n_{h}(y,s))}+\N{\del_yn_{h}(y,s)}\right]_{y=x}
	\\&\leq \NNN{f-h}_{\nu,\mu,x}+\N{n_{h_2}(s)}=\NNN{f-h}_{\nu,\mu,x}+\N{\del_xf_0(x,p)}_{L^1_p(\T^d)}
	\end{align*}
	and likewise,
	
	\begin{align*}
	\N{\del_xE_f(x,s)}\leq\NN{\epsilon}_{L^\infty}\left(\NNN{f-h}_{\nu,\mu,x}+\N{\del_xf_0(x,p)}_{L^1_p(\T^d)}\right).
	\end{align*}
	Since $\nu\N{\del_xf_0(x,p)}_{L^1_p(\T^d)}\leq \NN{f_0}_{\dot C^{\nu}_x}$, there exist a constant $C>0$ such that for sufficiently large $\mu>0$ it holds
	
	\begin{equation}\label{in.proof.constant.C.for.beta}
	\NN{n_f}_{\dot C^{\nu-\mu t}_x}+\NN{E_f}_{\dot C^{\nu-\mu t}_x}\leq C(1+\nu^2)\NN{f_0}_{\dot C^{\nu}_x}+ \frac C\mu n_{f_0}(x)(1-\eta n_{f_0}(x))
	\end{equation}
	for all $0\leq t< T$. 
	By the hypothesis, we have
	
	\begin{align*}
	\NN{n_f}_{\dot C^{\nu-\mu s}_x}+\NN{E_f}_{\dot C^{\nu-\mu s}_x}&\leq C\left(\beta(1+\nu^2)+\frac1\mu\right)n_{f_0}(x)(1-\eta n_{f_0}(x))
	\\&\leq\frac\alpha2 n_{f_0}(x)(1-\eta n_{f_0}(x)) 
	\end{align*}
	for $0\leq t <T$ if $(1+\nu^2)\beta\leq \alpha/(4C)$ and $\mu>0$ is sufficiently large. 
	
	However, we still need to ``replace" $n_{f_0}$ by the density of the solution $f$ in the estimate. In order to show that $n_f$ and $n_{f_0}$ are closely related, we use the equation for $n_f$ and derive similarly as above that
	
	\begin{align*}
	\N{\del_t (n_f(1-\eta n_f))} &\leq  \NN{v}_{L^\infty(\T^d)} \int_{\T^d} \N{\na f} dp
	\\&\leq
	\NN{\velo}_{L^\infty}\left(\NNN{f-h}_{\nu,\mu,x}+\N{\del_xf_0(x,p)}_{L^1_p(\T^d)}\right)
	\end{align*}
	which entails
	
	\begin{align}\label{eq.in.proof nf1-nf.Anfang}
	n_f(x,t)(1-\eta n_f(x,t))\geq (1-Ct)n_{f_0}(x)(1-\eta n_{f_0}(x))
	\end{align}
	for some $C>0$ and all $t\leq T<\nu/\mu$ if $\mu$ is sufficiently large. Thus, we even have
	
	\begin{align*}
	n_f(x,t)(1-\eta n_f(x,t))\geq \frac12n_{f_0}(x)(1-\eta n_{f_0}(x))  
	\end{align*}
	if $\mu>0$ is sufficiently large implying
	
	\begin{align*}
	\NN{n_f}_{\dot C^{\nu-\mu s}_x}+\NN{E_f}_{\dot C^{\nu-\mu s}_x}\leq \delta n_{f}(x,s)(1-\eta n_{f}(x,s)) .
	\end{align*}
	This proves the claim. Finally, there is only one assertion left:
	\\
	{\bf Claim 3:}  $(n_f(x,t),E_f(x,t))\in \overline{\mathcal{U}_{a,\alpha}}$ for all $x\in \R^d$ and $0\leq t<T$ if $\mu>0$ is sufficiently large.
	\\
	Recall that
	\[	\mathcal U_{a,\alpha}
	=
	%	\left\{y\in \R^2: \exists (m_0,m_1)\in M_a \mbox{ with }|(m_0,m_1)-y|<\delta m_0(1-\eta m_0)\right\}
	%	\\&=
	\bigcup_{(m_0,m_1)\in M_a}B_{\delta m_0(1-\eta m_0)}(m_0,m_1).\]
	Similar to \eqref{eq.in.proof nf1-nf.Anfang}, we obtain that
	
	\begin{align*}
	\N{n_f(x,t)-n_0(x)}+\N{E_f(x,t)-E_0(x)}\leq Ctn_{f_0}(x)(1-\eta n_{f_0}(x))
	\end{align*}
	for some $C>0$ independent from $x$ and small $t>0$. Since by assumption $(n_{f_0}(x),E_{f_0}(x))\in \mathcal U_{a,\alpha/2}$, there exits $(m_0,m_1)\in M_a$ such that
	\[|(n_{f_0}(x),E_{f_0}(x))-(m_0,m_1)|<\frac\alpha2 m_0(1-\eta m_0).\]
	Thus, we compute that
	
	\begin{align*}
	|(n_{f}(x,t),E_{f}(x,t))-(m_0,m_1)|&\leq |(n_{f_0}(x),E_{f_0}(x))-y|+ Ctn_{f_0}(x)(1-\eta n_{f_0}(x))
	\\&<\frac\alpha2 m_0(1-\eta m_0)+ Ctn_{f_0}(x)(1-\eta n_{f_0}(x)).
	\end{align*}
	Hence, for sufficiently large $\mu>0$, it holds 
	
	\begin{align*}
	|(n_{f}(x,t),E_{f}(x,t))-(m_0,m_1)|<\delta m_0(1-\eta m_0).
	\end{align*}
	for all $0\leq t<T<\nu/\mu$.
	%		Concluding, we have proved that there exists a $\tau=\tau_4>0$ and a $\beta=\delta/(4C)$ with $C>0$ as in \eqref{in.proof.constant.C.for.beta} such that $\Theta:Z^\tau\to Z^\tau$ is a contraction.
\end{proof}
\begin{proof}[Proof of Proposition \ref{prop.BGK.contraction}]
	Combining the previous two lemmata, we see that $\Theta$ is a contraction and admits a unique fixed-point being the solution of  equation \eqref{3.be.nochmal}.
\end{proof}
%\begin{lemma}\label{lemma.Theta:Z^t}
%	Let $\eta>0$, $a\geq1$, $R,\nu,\mu,T,\delta>0$ be as in Definition \ref{def.Z^t}. If $f_0:\R^d\times\T^d\to \R$ is analytic such that $\sup_{x\in\R^d}\NN{f_0}_{C^\nu_x}
%	\leq R/2$, then 
%\end{lemma}
\begin{proof}[Proof of Theorem \ref{thm.BGK.local.solution.bitorus}]
	The proof is similar to the proof of Theorem \ref{thm.BGK.local.solution.bitorus.easy}. This time we want to apply Proposition \ref{prop.BGK.contraction} and thus have to show that the initial data satisfy its hypothesis.
	
	Again we use the fact that $f_0$ is continuous to guarantee that there exists an $\theta>0$ such that $\theta<f_0<\eta^{-1}(1-\theta)$. Likewise to the proof of  Theorem \ref{thm.BGK.local.solution.bitorus.easy}, we can use the analyticity of $f_0$ to show that $R:=2\sup_{x\in \T^d}\NN{f_0}_{C^{\nu_0}_x}+1<\infty$ for sufficiently small $\nu_0>0$. Moreover, it is easy to check that 
	
	\[\NN{f_0}_{C^{\nu}_x}\leq \frac R2\quad\mbox{and}\quad \NN{f_0}_{C^{\nu}_x}\leq R\frac{\nu}{2\nu_0}\]
	holds for all $0<\nu\leq\nu_0$ and $x\in \T^d$.	Using the bounds for $f_0$, we see that
	
	\[n_{f_0}(1-\eta n_{f_0})\geq \theta^2\]
	and thus given $\beta$, we have
	
	\[\NN{f_0}_{C^{\nu}_x}\leq \frac R2\quad\mbox{and}\quad\NN{f_0}_{C^{\nu}_x}\leq \beta n_{f_0}(x)(1-\eta n_{f_0}(x))\]
	for all $x\in \T^d$ if $\nu\leq \min\{\nu_0,\beta\nu_0/(K\theta^2)\}$.
	
	The next step is to show the hypothesis on the macroscopic densities of $f_0$. We claim that 
	\[(n_{f_0}(x),E_{f_0}(x)):=\int_{\T^d}(1,(\epsilon(p)))f_0(x,p)dp\in \mathcal U_{a,\alpha/2}\]
	for some $a\geq1$ and given $\alpha>0$.
	According to \cite{Bra17} section 5.1 and $\theta<f_0<\eta^{-1}(1-\theta)$, there exists $\lambda^0=(\lambda_0^0,\lambda_1^0):\T^d\to\R^2$ analytic and bounded such that
	
	\begin{align*}
	n(\lambda^0)=n_{f_0}\quad\mbox{and}\quad E(\lambda^0)=E_{f_0}.
	\end{align*}
	Thus,
	
	\begin{align*}
	(n_{f_0}(x),E_{f_0}(x)\in \left\{\int_{\T^d}\frac{(1,\epsilon(p))dp}{\eta+e^{-\lambda_0-\lambda_1\epsilon(p)}}: \lambda_0,\lambda_1\in\R\mbox{ with } \N{\lambda_1}\leq \log a\right\}\subset\mathcal U_{a,\alpha}
	\end{align*}
	for $a:=\exp(\|\lambda^0_1\|_{L^\infty})$ and all $\alpha>0$. Finally, we can apply Proposition \ref{prop.BGK.contraction} and obtain the assertion.
\end{proof}

\begin{proof}[Proof of Theorem \ref{thm.BGK.local.solution}]
	The idea is to adjust the parameter such that the hypothesis of Proposition \ref{prop.BGK.contraction} are fulfilled. At first, we see by Lemma \ref{lem.NN.Frho.nE} and Proposition \ref{lem: nnorm leq C norm 2} that there exists an $R>0$ and a $\nu_0>0$ such that $\sup_{x\in\R^d}\NN{f_0}_{C^{\nu_0}_x}\leq R/2$ and
	
	\begin{align*}
	\NN{f_0}_{\dot C^{\nu_0}_x}\leq C\nu_1n_{f_0}(x)(1-\eta n_{f_0}(x))
	\end{align*}
	holds for all $x\in\R^d$ and all $\nu_1\leq\nu_0$. Now, we set $\alpha_0=\alpha/2$ and $\nu_1:=\beta/{C}$, where $\alpha$ and $\beta$ are given by Proposition \ref{prop.BGK.contraction}. Then Proposition \ref{prop.BGK.contraction} guarantees a unique analytic solution $f$ on a short time interval. The well-posedness is then a direct consequence of Theorem \ref{thm: boltzmann local: hT solution2}. Finally, using Lemma \ref{lem: nnorm leq C norm 2}, we obtain the well-posedness also in the desired norm with a larger constant.
\end{proof}
Finally, we note that Theorem \ref{thm.BGK.local.solution0} is actually a corollary of Theorem \ref{thm.BGK.local.solution}.

\appendix
\section{Proof of Proposition \ref{lem.Ff.is.Lipschitz}}
\begin{definition}\label{def.FrhoE}
	Let  $\lambda_0=\lambda_0(n,E)$ and $\lambda_1=\lambda_1(n,E)$ be functions of the densities $n,E$ given by
	
	\begin{align*}%\label{intro:rhoEvonlambda}
	\binom{n}{E}=\int_{\T^d}\binom{1}{\epsilon(p)}\frac{dp}{\eta+e^{-\lambda_0(n,E)-\lambda_1(n,E)\epsilon(p)}}.
	\end{align*}
	We define 
	
	\begin{equation*}
	\FrhoE(n,E;p):=\frac{1}{\eta+e^{-\lambda_0(n,E)-\lambda_1(n,E)\epsilon(p)}}
	\end{equation*}
	for $(n,E)\in \{\int_{\T^d}(1,\epsilon(p))g(p)dp: g\in L^1(\T^d;(0,\eta^{-1}))\}$ and $p\in \T^d$.
	
\end{definition}
Our goal is to estimate the norm of $\F_f$ be means of $f$. Due to the preceding Definition, we can rewrite $\F_f$ as a composition by
\[
\F_f(x,p)= \FrhoE(n_f(x),E_f(x);p),\]
where $n_f(x):=\int_{\T^d}f(x,p)dp$ and $E_f(x)=\int_{\T^d}\epsilon(p)f(x,p)dp$. 

Thus, we can easily compute the first derivative of $\F_f$ w.r.t.~$x$ as $\partial_x\F_f=\partial_n\FrhoE(n_f,E_f) \partial_xn_f+\partial_E\FrhoE(n_f,E_f)\partial_xE_f$ by using the chain rule and the following Lemma.
\begin{lemma}\label{lem.Berechnung.von.G}
	\[
	\frac{\partial\FrhoE}{\partial(n,E)}(n,E;p)=\frac{\FrhoE(n,E;p)(1-\eta\FrhoE(n,E;p))}{\int_{\T^d}\epsilon^2d\mu\int_{\T^d}1d\mu-\big(\int_{\T^d}\epsilon d\mu\big)^2}\int_{\T^d}\binom{-\epsilon(p')}{1}(\epsilon(p)-\epsilon(p'))d\mu_{p'},
	\]
	where $d\mu_p:=\FrhoE(n,E;p)(1-\eta\FrhoE(n,E;p))dp$.
\end{lemma}
\begin{proof}
	For $\lambda=(\lambda_0,\lambda_1)\in \R^2$, let us denote  $\F(\lambda;p):=1/(\eta+e^{-\lambda_0-\lambda_1\epsilon(p)})$  and	
	
	\begin{align*}
	n(\lambda):=\int_{\T^d}\frac{dp}{\eta+e^{-\lambda_0-\lambda_1\epsilon(p)}},\quad E(\lambda):=\int_{\T^d}\frac{\epsilon(p)dp}{\eta+e^{-\lambda_0-\lambda_1\epsilon(p)}}.
	\end{align*}
	We have
	\[
	\frac{\partial(n,E)}{\partial\lambda}(\lambda)=\int_{\T^d}\begin{pmatrix}
	1&\epsilon(p)\\\epsilon(p)&\epsilon(p)^2 
	\end{pmatrix}
	d\mu_p,
	\]
	where $d\mu_p:=\F(\lambda;p)(1-\eta\F(\lambda;p))dp$. Since $d\mu_p$ is a positive measure, we can use the Cauchy-Schwarz inequality to see that $\frac{\partial(n,E)}{\partial\lambda}(\lambda)$ is invertible. Finally, we easily compute 
	\[
	\frac{\partial\lambda}{\partial(n,E)}=\frac{1}{\int_{\T^d}\epsilon^2d\mu\int_{\T^d}1d\mu-\big(\int_{\T^d}\epsilon d\mu\big)^2}\int_{\T^d}\begin{pmatrix}
	\epsilon^2&-\epsilon\\-\epsilon&1
	\end{pmatrix}d\mu
	\]
	by the inverse function theorem and the chain rule ensures the assertion.
\end{proof}
Note that our main techniques are based on the analytic norms	
\begin{equation*}
\NN{f}_{C^{\nu}}:= \sum_{a,b\in\mathbb N_0^d} \frac{\nu^{|a+b|}}{a!b!}\|\del_x^a\del_p^b f\|_{W^{1,\infty}_xW^{1,1}_p}
\end{equation*}
for $f:\R^d\times\T^d\to \R^k$ being analytic, where we use the notation

\begin{equation*}
\|f\|_{W^{1,\infty}_xW^{1,1}_p}:=\sum_{\substack{a,b\in\mathbb N_0^d\\|a+b|\leq 1}}\|\del_x^a\del_p^bf\|_{L^{\infty}_xL^{1}_p}\quad\mbox{and}\quad \|f\|_{L^{\infty}_xL^{1}_p}:=\sup_{x\in \R^d}\int_{\T^d}|f(x,p)|dp.
\end{equation*}

This motivates Proposition \ref{lem.Ff.is.Lipschitz}, which we restate for the reader's convenience.
\begin{proposition}\label{lem.Ff.is.Lipschitz.nochmal}
	Let $\eta,\nu,R>0$ and $\alpha>0$. There exists an $C,\delta>0$ such that the following is true.
	
	Let $f,g:\R^d\times\T^d\to[\alpha,\eta^{-1}-\alpha]$ be analytic satisfying $\NN{f}_{C^{\nu}},\NN{g}_{C^{\nu}}\leq R$ and
	
	\begin{equation*}
	\NN{n_h-\bar n}_{C^{\nu}}+\NN{E_h-\bar E}_{C^{\nu}}\leq \delta\quad\mbox{for }h\in\{f,g\}
	\end{equation*}
	and some $\bar{n},\bar E\in \R$, it holds
	
	\begin{equation*}
	\NN{\F_f}_{C^\nu},\NN{\F_g}_{C^\nu}\leq C
	\end{equation*}
	and
	
	\begin{equation*}
	\NN{\F_f-\F_g}_{C^\nu}\leq C\NN{f-g}_{C^\nu}.
	\end{equation*}
\end{proposition}
The main steps to prove this proposition is again to consider $\F_f$ as the composition	\begin{equation*}
\FrhoE(n,E;p):=\frac{1}{\eta+e^{-\lambda_0(n,E)-\lambda_1(n,E)\epsilon(p)}}
\end{equation*}
for $(n,E)\in \{\int_{\T^d}(1,\epsilon(p))g(p)dp: g\in L^1(\T^d;(0,\eta^{-1}))\}$ and $p\in \T^d$.

In the analytic norm $\NN{\cdot}_{C^\nu}$ involves all derivatives. As a first step we consider the derivatives of $\FrhoE$. Using the inverse mapping theorem, we can easily see that $\lambda_0,\lambda_1$ are analytic in  their domain $\{\int_{\T^d}(1,\epsilon(p))g(p)dp: g\in L^1(\T^d;(0,\eta^{-1}))\}$. This proves the following.
\begin{lemma}
	$\FrhoE$ is analytic on $\{\int_{\T^d}(1,\epsilon(p))g(p)dp: g\in L^1(\T^d;(0,\eta^{-1}))\}\times\T^d$. In particular, for all 
	\[(\bar{n},\bar E)\in \left\{\int_{\T^d}(1,\epsilon(p))g(p)dp: g\in L^1(\T^d;(0,\eta^{-1}))\right\}\]
	there exist constants $A\geq0$ such that
	
	\begin{align}\label{eq.42}
	\N{\partial_{(n,E)}^i\partial_p^j\FrhoE(\bar n,\bar E;p)} 
	&\leq
	i!j!A^{|i|+|k|}
	\end{align}
	for all $p\in \T^d$ for $i\in\mathbb N_0^2$ and $j\in\mathbb N_0^d$.
\end{lemma}
\begin{corollary}\label{cor.FrhoE.C.delta}
	Let $(\bar{n},\bar E)\in M:=\left\{\int_{\T^d}(1,\epsilon(p))g(p)dp: g\in L^1(\T^d;(0,\eta^{-1}))\right\}$. Then there exist a $\delta>0$ and an open neighborhood $\mathcal U\subset M$ of $(\bar{n},\bar E)$ such that
	\begin{equation*}
	\NN{\FrhoE}_{C^{\delta}(\mathcal U)}:= \sum_{|i|+|j|\leq 1}\sum_{a\in \mathbb N_0^2,b\in\mathbb N_0^d} \frac{\delta^{|a|+|b|}}{a!b!}\sup_{(n,E)\in\mathcal U}\int_{\T^d}|\del_{(n,E)}^{a+i}\del_p^{b+j} \FrhoE(n, E;p)|dp
	\end{equation*}
	is finite
\end{corollary}
\begin{proof}
	According to the previous lemma, there exists an $A>0$ such that the estimate \eqref{eq.42} is satisfied. Using the Taylor formula for $\FrhoE$ w.r.t.\ $(n,E)$ makes sure that
	\begin{align}
	\N{\partial_{(n,E)}^i\partial_p^j\FrhoE(n,E;p)} 
	&\leq
	i!j!(2A)^{|i|+|k|}
	\end{align}
	holds true in a neighborhood $\mathcal U\subset M$ of $(\bar{n},\bar E)$. Thus, summing up all derivatives with the right weight, we can show that $\NN{\FrhoE}_{C^{\nu_0}(\mathcal U)}$ for $\delta<1/(2A)$.
\end{proof}
The last ingredient for the proof of Proposition \ref{lem.Ff.is.Lipschitz.nochmal} is a formula for the analytic norms of composition of functions which is in fact a corollary of the Fa\`a di Bruno formula. It was firstly derived by \cite{MoVi11}. Note that Mouhot and Villani \cite{MoVi11} also state a version for $d>1$. However, in their proof, they use only the one dimensional Fa\`{a} di Bruno formula such that they leave the multidimensional case to the reader. For $d\geq1$, we also refer to \cite{Bra17} Lemma 4.2.5, where the definition of the norm $\N{\cdot}_{C^\nu}$ slightly differs from our case and involves full derivatives. The same techniques can still be used for this case.
\begin{lemma}\label{lemma.Verkettung2}
	Let $x\in V\subset \R^k$ open and let $g:\R^d\times\T^d\to V$, $\phi:V\to \R$ be analytic. Then
	
	\begin{equation*}
	\N{\phi\circ g}_{C^\nu}\leq \N{\phi}_{C^\mu} \qquad\text{ for }\mu=\N{g-v}_{C^\nu}\
	\end{equation*}
	for $\nu>0$ and all $v\in V$, where	
	
	\begin{equation*}
	\N{g}_{C^{\nu}}:= \sum_{a,b\in\mathbb N_0^d} \frac{\nu^{|a+b|}}{a!b!}\|\del_x^a\del_p^b g\|_{L^{\infty}_xL^{1}_p}
	\end{equation*}
	for $f:\R^d\times\T^d\to \R^k$ being analytic.
\end{lemma}
\begin{corollary}
	Given $\bar n,\bar E\in\R$ and $\nu>0$. Let $\delta>0$ and $\mathcal U$ be as in Corollary \ref{cor.FrhoE.C.delta}. Then there exists a $C>0$ such that for all $(n,E):\R^d\to\mathcal \R^2$ being analytic such that
	
	\begin{equation*}
	\NN{n-\bar n}_{C^{\nu}}+\NN{E-\bar E}_{C^{\nu}}\leq \delta,
	\end{equation*}
	it holds

	\begin{equation*}
	\NN{\FrhoE(n,E)}_{C^\nu}\leq C.
	\end{equation*}
\end{corollary}
\begin{proof}
	Using the analytic norms from Lemma \ref{lemma.Verkettung2}, we can write
	
	\begin{align*}
	\NN{\FrhoE(n,E)}_{C^\nu}=\N{\FrhoE(n,E)}_{C^\nu}+\sum_{i=1}^d\N{\del_{x_i}\FrhoE(n,E)}_{C^\nu}+\sum_{i=1}^d\N{\del_{p_i}\FrhoE(n,E)}_{C^\nu}.
	\end{align*}
	By Lemma \ref{lemma.Verkettung2}, we obtain
	
	\begin{align*}
	\N{\FrhoE(n,E)}_{C^\nu}\leq \N{\FrhoE}_{C^\delta(\mathcal U)},\quad\mbox{because }\N{n-\bar n}_{C^{\nu}}+\N{E-\bar E}_{C^{\nu}}\leq \delta.
	\end{align*}
	By assumption $\N{\FrhoE}_{C^\delta(\mathcal U)}<\infty$ and thus, $\N{\FrhoE(n,E)}_{C^\nu}$ is bounded. We can do the same trick for the other terms. Her we only need to use the chain rule and the submultiplicativity of $|\cdot|_{C^\nu}$ to split the terms into
	
	\begin{align*}
	\N{\del_{x_i}\FrhoE(n,E)}_{C^\nu_x}=\N{\del_{(n,E)}\FrhoE(n,E)\del_{x_i} (n,E)}_{ C^\nu}&\leq \N{\del_{(n,E)}\FrhoE(n,E)}_{C^\nu}\N{\del_{x_i} (n,E)}_{ C^\nu}
	\\&\leq \N{\del_{(n,E)}\FrhoE(n,E)}_{C^\nu}\NN{(n,E)}_{ C^\nu}
	\end{align*}
	with slightly abuse of notation. Note that a version of Corollary \ref{cor.FrhoE.C.delta} for $\del_{(n,E)}\FrhoE(n,E)$ holds true. This can be shown in the same manner as for Corollary \ref{cor.FrhoE.C.delta}.			
\end{proof}
Without loss of generality, we can assume that $\mathcal U$ is convex. We can apply the same arguments for $\del_{(n,E)}\FrhoE(n,E)$ and obtain by 
\[\FrhoE(n_1,E_1)-\FrhoE(n_0,E_0)=\binom{n_1-n_0}{E_1-E_0}\cdot\int_0^1\del_{(n,E)}\FrhoE(n_1t+(1-t)n_0,E_1t+(1-t)E_0)dt.\]
This leads to the following statement.
\begin{corollary}
	Given $\bar n,\bar E\in\R$ and $\nu>0$. Let $\delta>0$ and $\mathcal U$ be as in Corollary \ref{cor.FrhoE.C.delta}. Let $\mathcal U'\subset U$ be convex. Then there exists a $C>0$ such that for all $(n_i,E_i):\R^d\to\mathcal U'$, $i=0,1$, being analytic such that
	
	\begin{equation*}
	\NN{n_i-\bar n}_{C^{\nu}}+\NN{E_i-\bar E}_{C^{\nu}}\leq \delta,\quad i=0,1,
	\end{equation*}
	it holds

	\begin{equation*}
	\NN{\FrhoE(n_1,E_1)-\FrhoE(n_0,E_0)}_{C^\nu}\leq C\left(\NN{n_1-n_0}_{C^\nu}+\NN{E_1-E_0}_{C^\nu}\right).
	\end{equation*}
\end{corollary}
\begin{proof}[Proof of Proposition \ref{lem.Ff.is.Lipschitz.nochmal}]
	The assertion is basically a direct consequence of the foregoing corollaries. The only the difference is that we do not want to assume explicitly that $(n,E)(\R^d)\subset \mathcal U$. We can neglect this hypothesis by choosing $\delta$ sufficiently small such that there exist a ball $B_{\delta}(\bar n,\bar E)\subset \mathcal U$ with radius $\delta$. Then
	
	\begin{align*}
	\NN{n-\bar n}_{L^\infty}+\NN{E-\bar E}_{L^\infty}\leq 	\NN{n-\bar n}_{C^{\nu}}+\NN{E-\bar E}_{C^{\nu}}\leq \delta.
	\end{align*}
	implies that $(n(x),E(x))\in B_{\delta}(\bar n,\bar E)\subset \mathcal U$ for all $x\in \R^d$.
\end{proof}
\section{Proof of Proposition \ref{prop.Ff.is.Lipschitz}}
In this section we are going to prove Proposition \ref{prop.Ff.is.Lipschitz}, which we restate for the sake of convenience.
Let 

\begin{equation*}
\NN{f}_{C^{\nu}_x}:= \sum_{\substack{i,j\in\mathbb N_0^d\\|i+j|\leq 1}}|\del_x^{i}\del_p^{j} f|_{C^{\nu}_x},\ \mbox{where }\N{f}_{C^{\nu}_x}:=\sum_{a,b\in\mathbb N_0^d} \frac{\nu^{|a+b|}}{a!b!}\int_{\R^d}|\del_x^{a}\del_p^{b} f(x,p)|dp
\end{equation*}
and

\begin{equation*}
\N{f}_{\dot C^{\nu}_x}:= \N{f}_{C^{\nu}_x}-\N{f}_{C^{0}_x}	\quad\mbox{and}\quad\NN{f}_{\dot C^{\nu}_x}:= \NN{f}_{C^{\nu}_x}-\NN{f}_{C^{0}_x}.
\end{equation*} 
\begin{proposition}\label{prop.Ff.is.Lipschitz.nochmal}
	Let $\eta>0$, $a\geq1$ and $R,\nu>0$. Then there exist $\alpha,\delta>0$ such that for all $f,g:\R^d\times\T^d\to\R$ being analytic with $\NN{f}_{C^\nu_x},\NN{g}_{C^\nu_x}\leq R$ and
	
	\begin{equation*}
	\NN{n_h}_{\dot C^\nu_x}+\NN{E_h}_{\dot C^\nu_x}\leq \delta n_h(x)(1-\eta n_h(x))\quad\mbox{for }h\in\{f,g\}
	\end{equation*}
	and $(n_h(x),E_h(x))\in \mathcal U_{a,\alpha}$ for $h\in\{f,g\}$, it holds 
	
	\begin{equation*}
	\NN{\F_f}_{C^\nu_x},\NN{\F_g}_{C^\nu_x}\leq C
	\end{equation*}
	and
	
	\begin{equation*}
	\NN{\F_f-\F_g}_{C^\nu_x}\leq C\NN{f-g}_{C^\nu_x}.
	\end{equation*}
	for some $C>0$ and all $x\in \R^d$.
\end{proposition} 
%\begin{definition}\label{def.FrhoE}
%	Let  $\lambda_0=\lambda_0(n,E)$ and $\lambda_1=\lambda_1(n,E)$ be functions of the densities $n,E$ given by
%
%	\begin{align*}%\label{intro:rhoEvonlambda}
%	\binom{n}{E}=\int_{\T^d}\binom{1}{\epsilon(p)}\frac{dp}{\eta+e^{-\lambda_0(n,E)-\lambda_1(n,E)\epsilon(p)}}.
%	\end{align*}
%	We define 
%
%	\begin{equation*}
%	\FrhoE(n,E;p):=\frac{1}{\eta+e^{-\lambda_0(n,E)-\lambda_1(n,E)\epsilon(p)}}
%	\end{equation*}
%	for $(n,E)\in \{\int_{\T^d}(1,\epsilon(p))g(p)dp: g\in L^1(\T^d;(0,\eta^{-1}))\}$ and $p\in \T^d$.
%	
%\end{definition}

%However,  for further derivatives we have to apply Fa\`a di Bruno's formula which is rather long. Since we are only interested in the analytic norms, we can use the following formula for the norm of composition of functions.
Note that this proposition is stronger than its counter part in Proposition \ref{lem.Ff.is.Lipschitz.nochmal}. Therefore, we require a more sophisticated analysis of $\FrhoE$.
\begin{definition}
	Let $a\geq 1$ and $\delta>0$. We define
	%	\begin{align*}
	%	Z_a:=\left\{p\mapsto\frac{1}{\eta+e^{-\lambda_0-\lambda_1\epsilon(p)}}: \lambda_0,\lambda_1\in\R\mbox{ with } \N{\lambda_1}\leq \log a\right\}
	%	\end{align*}
	%	and
	%	\begin{align*}
	%	M_a:=\left\{\int_{\T^d}(1,\epsilon(p))f(p)dp: f\in Z_a \right\}.
	%	\end{align*}
	
	\begin{align*}
	M_a:=\left\{\int_{\T^d}\frac{(1,\epsilon(p))dp}{\eta+e^{-\lambda_0-\lambda_1\epsilon(p)}}: \lambda_0,\lambda_1\in\R\mbox{ with } \N{\lambda_1}\leq \log a\right\}\subset \mathbb R^2
	\end{align*}
	and
	
	\begin{align*}
	\mathcal U_{a,\delta}
	&:=
	%	\left\{y\in \R^2: \exists (m_0,m_1)\in M_a \mbox{ with }|(m_0,m_1)-y|<\delta m_0(1-\eta m_0)\right\}
	%	\\&=
	\bigcup_{(m_0,m_1)\in M_a}B_{\delta m_0(1-\eta m_0)}(m_0,m_1)\supset M_a,
	\end{align*}
	where $B_\theta(y)$ denotes the ball in $\mathbb R^2$ centered at $y$ with radius $\theta$.
\end{definition}
\begin{lemma}\label{Ableitung.FrhoE}
	Let $a\geq1$. There exist constants $A_a,B_a>0$ such that
	
	\begin{align*}
	\N{D_{(n, E)}^iD_p^j\FrhoE(n,E;p)} 
	&\leq
	i!j!A_{a}^{j}\left(\frac{B_{a}}{n(1-\eta n)}\right)^{i}\FrhoE(n,E;p)(1-\eta\FrhoE(n,E;p).
	\end{align*}
	for all $(n,E)\in M_a$, $p\in \T^d$ and $i+j\geq1$. Moreover, if $\eta=0$ these constant may be chosen independently from $a$, i.e., there exist $A,B>0$ such that 
	
	\begin{align*}
	\N{\del_{(n,E)}^i\del_p^j\FrhoE(n,E;p)} 
	&\leq
	i!j!A^{j}\frac{B^i}{n^i}\FrhoE(n,E;p)
	\end{align*}
	for any $i+j\geq1$ and all $(n,E)\in [0,\infty)\times\R$.
\end{lemma}
\begin{proof}
	For a detailed proof see \cite{Bra17} section 5.4.
\end{proof}
%\begin{lemma}
%	Let $a\geq1,\mu>0$ and $(n,E):\R^d\to M_a$ be analytic such that 
%	\[\N{n}_{\dot C^\mu_x}+\N{E}_{\dot C^\mu_x}\leq C n(x)(1-\eta n(x))\]
%	for some $C>0$. Then there exists a $\nu_0\in(0,\mu)$ and a $C_0>0$ such that
%	\begin{equation*}
%	\NN{\FrhoE(n,E)}_{\dot C^\nu_x}\leq \lambda C_0 n(x)(1-\eta n(x))
%	\end{equation*}
%	for all $\nu\in(0,\nu_0)$.
%\end{lemma}	

In the next step, we state the space local version of Lemma \ref{lemma.Verkettung2}, which can be proved exactly like  Lemma \ref{lemma.Verkettung2}.
\begin{lemma}\label{lem.B4}
	Let $x\in V\subset \R^k$ open and let $g:\R^d\times\T^d\to V$, $\phi:V\to \R$ be analytic. Then
	
	\begin{equation*}
	\N{\phi\circ g}_{\dot C^\nu_x}\leq \N{\phi}_{\dot C^\mu_y} \qquad\text{ with }\mu=\N{g}_{\dot C^\nu_x}\text{ and }y=g(x)
	\end{equation*}
	for $\nu>0$.
\end{lemma}
Using this lemma, we can easily find estimates for the derivatives of some functions
\begin{example}\label{exappendix}
	Let 
	
	\[F_{\lambda^0}(x)=\frac{1}{\eta+1+x^2}.\]
	We have
	
	\begin{equation*}
	\N{F_{\lambda^0}}_{\dot C^\nu_x}=\N{\phi\circ (\cdot)^2}_{\dot C^\nu_x}\leq \N{\phi}_{\dot C^\mu_{x^2}} \qquad\text{ with }\mu=\N{(\cdot)^2}_{\dot C^\nu_x}=2\nu|x|+\nu^2
	\end{equation*}
	for $\phi(s)=(\eta+1+s)^{-1}$ according to Lemma \ref{lem.B4}. We have
	$\phi^{(i)}(s)=(-1)^ii!(\eta+1+s)^{-(i+1)}$ which implies that
	
	\begin{align*}
	\N{F_{\lambda^0}}_{\dot C^\nu_x}&=\sum_{i=1}^N\frac{(\nu+x)^{2i}}{i!}|\phi^{(i)}(x^2)|
	\\&=\frac{1}{(\eta+1+x^2)}\sum_{i=1}^N\left(\frac{2\nu|x|+\nu^2}{\eta+1+x^2}\right)^i.
	%=\frac{1}{(\eta+1+x^2)}\frac{2\nu|x|+\nu^2}{(\eta+1+x^2)-2\nu|x|+\nu^2}
	\end{align*}
	Let $\nu:=\frac12\min\{\sqrt\eta,1\}$. Thus,
	
	\begin{align*}
	2\nu|x|-\nu^2\leq |x|-\frac{\eta}2= \frac{1+\eta + x^2}{2}-\frac{(1-|x|)^2}2\leq \frac{1+\eta + x^2}{2}.
	\end{align*}
	This implies that
	
	\begin{align*}
	\N{F_{\lambda^0}^{(a)}(x)}\leq 
	\frac{a!}{\nu^a}\N{F_{\lambda^0}}_{\dot C^\nu_x}
	&\leq \frac{a!}{\nu^a}\frac{1}{\eta+1+x^2}=
	\frac{a!2^a}{\min\{\sqrt\eta,1\}^a}\frac{1}{\eta+1+x^2}
	\end{align*}
	for $\nu=\frac12\min\{\sqrt\eta,1\}$.
\end{example}
\begin{corollary}\label{cor.Verkettung}
	Let $x\in V\subset \R^k$ open and let $m:\R^d\to V$, $\phi:V\times\T^d\to \R$ be analytic. We have
	
	\begin{equation*}
	\NN{\phi\circ m}_{\dot C^\nu_x}\leq \NN{\phi}_{\dot C^\mu_y}\big(1+\NN{m}_{ C^\nu_x}\big) +\mu\NN{\phi}_{ C^0_y}
	\end{equation*}
	with $\mu=\NN{m}_{\dot C^\nu_x}$ and $y=m(x)$ for $\nu>0$. Moreover, assume that $\N{\phi}_{C^{\mu_0}_y}<\infty$  for some $\mu_0>0$. Let $M>0$ and $\bar\mu\in(0,\mu_0)$. Then there exists a constant $C>0$ such that
	
	\begin{equation*}
	\NN{\phi\circ m}_{\dot C^\nu_x}\leq C \NN{m}_{\dot C^\nu_x}.
	\end{equation*}
	for all $\nu>0$ and all $m:\R^d\to V$ being analytic such that $\NN{m}_{ C^\nu_x}\leq M$ and $\NN{m}_{ C^\nu_x}\leq \bar\mu$.
\end{corollary}
\begin{proof}
	Using the chain rule we first compute
	
	\begin{align*}
	\N{\del_x\phi(m,\cdot)}_{C^\nu_x}=\N{\del_1\phi(m,\cdot)\del_x m}_{ C^\nu_x}\leq \N{\del_1\phi(m,\cdot)}_{ C^\nu_x}\N{\del_xm}_{ C^\nu_x}.
	\end{align*}
	Since $\N{f}_{\dot C^\nu_x}=\N{f}_{C^\nu_x}-\NN{f(x,p)}_{L^1_p(\T^d)}$ for $f:\R^d\times\T^d\to\R$ analytic, we have
	\begin{align*}
	\N{\del_x\phi(m,\cdot)}_{\dot C^\nu_x}&\leq \N{\del_1\phi(m,\cdot)}_{ C^\nu_x}\N{\del_xm}_{ C^\nu_x} - \NN{\del_1\phi(m,p)}_{L^1_p(\T^d)}\N{\del_x m(x)}
	% 		\\&=\N{\del_1\phi(m,\cdot)}_{ C^\nu_x}\N{\del_xm}_{\dot C^\nu_x}+\N{\del_1\phi(m,\cdot)}_{\dot C^\nu_x}\N{\del_xm(x)}
	\\&= \N{\del_1\phi(m,\cdot)}_{\dot C^\nu_x}\N{\del_xm}_{ C^\nu_x}+\N{\del_1\phi(m,p)}_{L^1_p(\T^d)}\N{\del_xm(x)}_{\dot C^\nu_x}
	\\&= \N{\del_1\phi(m,\cdot)}_{\dot C^\nu_x}\N{\del_xm}_{ C^\nu_x}+\N{\del_1\phi(m,\cdot)}_{ C^0_x}\N{\del_xm}_{\dot C^\nu_x}
	\end{align*}
	Now we can  conclude the first part of the assertion by Lemma \ref{lem.B4}. With this, the remaining part is a direct consequence of Lemma \ref{lem: nnorm leq C norm 2}.
\end{proof}
\begin{lemma}\label{lem.NN.Frho.nE}
	Let $\eta>0$, $a\geq1$, $C,\nu>0$. There exist  $\delta,\alpha,C_0>0$ and a neighborhood $\mathcal U_a$ of $M_a$ such that for all $x\in\R^d$, $(n,E):\R^d\to \mathcal U_{a,\alpha}$ being analytic in $x$, which satisfy $\NN{(n,E)}_{C^\nu_x}\leq C$ and
	\[\NN{n}_{\dot C^\nu_x}+\NN{E}_{\dot C^\nu_x}\leq \delta n(x)(1-\eta n(x)),\]
	we have  $\NN{\FrhoE(n,E)}_{C^\nu_x}\leq C_0$ and
	
	\begin{equation*}
	\NN{\FrhoE(n,E)}_{\dot C^\nu_x}\leq  C_0 n(x)(1-\eta n(x)).
	\end{equation*}
\end{lemma}
\begin{proof}
	Let $\alpha:=\frac{1}{2B_a}>0$, where $B_a$ is given by Lemma \ref{Ableitung.FrhoE}. For $y\in \mathcal U_{a,\alpha}$ we choose $m=(m_0,m_1)\in M_a$ such that $|(m_0,m_1)-y|<\frac{m_0(1-\eta m_0)}{2B_a}$. Note that $0<m_0<\eta^{-1}$.
	%	\[\mathcal U=\left\{y\in \R^2: (m_0,m_1)\in M_a \mbox{ with }|(m_0,m_1)-y|<\frac{m_0(1-\eta m_0)}{2B_a}\right\},\]
	%	where $B_a$ is given by Lemma \ref{Ableitung.FrhoE}. Assume $(n,E):\R^d\to \mathcal U$ is analytic with $\NN{(n,E)}_{C^\nu_x}\leq C$.  For $y\in \mathcal U$ let $m=(m_0,m_1)\in M_a$ such that $|(m_0,m_1)-y|<\frac{m_0(1-\eta m_0)}{2B_a}$. Then $0<m_0<\eta^{-1}$ and $\alpha:=\frac{m_0(1-\eta m_0)}{2B_a}>0$. 
	Writing $\beta=\frac{m_0(1-\eta m_0)}{2B_a}$, we use Taylor's formula and see that
	
	\begin{align*}
	\N{\FrhoE}_{\dot C^\nu_y} &=\sum_{i+j\geq1}\frac{\nu^{i+j}}{i!j!}\NN{\del^i_y\del^j_p\FrhoE(y,p)}_{L^1_p(\T^d)}
	\\&\leq
	\sum_{i+j\geq1}\frac{(\nu+\beta)^{i}\nu^j}{i!j!}\NN{\del^i_m\del^j_p\FrhoE(m,p)}_{L^1_p(\T^d)}
	\\&\leq
	\sum_{i+j\geq1} (\nu A_a)^j\left(\frac{(\nu+\beta) B_b}{m_0(1-\eta m_0)}\right)^i\NN{\FrhoE(m,p)(1-\eta\FrhoE(m,p))}_{L^1_p(\T^d)}
	\end{align*}
	Now by Jensen's inequality and the Neumann series, we obtain
	
	\begin{align*}
	\N{\FrhoE}_{\dot C^\nu_m} &\leq
	\left(\frac{1}{1-\nu A_a}\frac{m_0(1-\eta m_0)}{m_0(1-\eta m_0)-(\nu+\beta) B_a}-1\right) m_0(1-\eta m_0)
	\end{align*}
	if $\nu A_a<1$ and $(\nu +\beta)B_a< m_0(1-\eta m_0)$. Thus, defining \[\mu_0= \min\{1/(3A_a),m_0(1-\eta m_0)/(3B_a)\},\] we have $\N{\FrhoE}_{\dot C^{\mu_0}_y}\leq 8m_0(1-\eta m_0)$. Therefore, by Lemma \ref{lem: nnorm leq C norm 2}, there exists a constant $C_1>0$ such that
	\begin{align}\label{inproof.NN.FrhoE.Abschaetzung}
	\NN{\FrhoE}_{\dot C^\mu_y}\leq  C_1 \NN{(n,E)}_{\dot C^\nu_x}
	\end{align}
	if $\NN{(n,E)}_{\dot C^\nu_x}\leq \min\{1/(4A_a),m_0(1-\eta m_0)/(4B_a)\}$.
	
	For the next step, we suppose that $(n,E)$ from above fulfills the hypothesis. Then $\NN{(n,E)}_{\dot C^\nu_x}\leq 1/(\eta A_a)n(x)(1-\eta n(x))\leq 1/(4A_a)$ since $n(x)(1-\eta n(x))\leq \eta/4$. Moreover, it holds $\NN{(n,E)}_{\dot C^\nu_x}\leq 1/(4B_a)n(x)(1-\eta n(x))$. 
	%For the next step, let $(n,E):\R^d\to \mathcal U$ satisfy the hypothesis. 
	%	Again by Lemma \ref{lem: nnorm leq C norm 2}, for $\mu_0\in(0,\mu)$ constant $C_2>0$ such that
	%	\begin{align*}
	%		\NN{(n,E)}_{\dot C^\nu_x}\leq \nu \frac{C_2}{C}\left(\N{n}_{\dot C^{\mu}_x}+\N{E}_{\dot C^{\mu}_m}\right)\leq \nu C_2 n(x)(1-\eta n(x)).
	%	\end{align*}
	%	for $\nu\in(0,\mu_0)$. 
	%	In order to combine this with \eqref{inproof.NN.FrhoE.Abschaetzung}, we use Corollary \ref{cor.Verkettung}.
	Thus, it holds \eqref{inproof.NN.FrhoE.Abschaetzung} and in particular
	
	\begin{align*}%\label{inproof.NN.FrhoE.Abschaetzung}
	\NN{\FrhoE}_{\dot C^\mu_y}\leq  C_1 \delta n(x)(1-\eta n(x)).
	\end{align*}
	Finally, we can easily show the remaining estimate by using the inequality
	
	\begin{multline*}
	\NN{\FrhoE(n,E)}_{C^\nu_x}\leq\NN{\FrhoE(n,E)}_{\dot C^\nu_x}+\NN{\FrhoE(n(x),E(x);p)}_{L^1_p(\T^d)} \\+\NN{\del_{(n,E)}\FrhoE(n(x),E(x);p)}_{L^1_p(\T^d)}\NN{\del_x(n(x),E(x))}+\NN{\del_p\FrhoE(n(x),E(x);p)}_{L^1_p(\T^d)}
	\end{multline*}
	and the fact that $\FrhoE$ and $\del_{((n,E),p)}\FrhoE(n,E)$ are bounded (see Lemma \ref{Ableitung.FrhoE}).
	%	Note that there exists a $\nu_0\in(0,\mu)$ such that 
	%	$\nu_0 C_2n(x)(1-\eta n(x))\leq  \min\{1/(4A_a),n(x)(1-\eta n(x))/(4B_a)\}$. 	Hence, we have 
	%		\begin{equation*}
	%		\NN{\FrhoE(n,E)}_{\dot C^\nu_x}\leq \NN{(n,E)}_{\dot C^{\nu}_x}C_1n(x)(1-\eta n(x))\leq\nu C_1C_2 n(x)^2(1-\eta n(x))^2
	%		\end{equation*}
	%		for all $\nu\in(0,\nu_0)$.
\end{proof}

\begin{lemma}\label{lem.FrhoE.Lipschitz}
	Let $\nu>0$ and $a\geq1$. For $C,\nu>0$ let $\alpha,\delta>0$ be as in Lemma \ref{lem.NN.Frho.nE}. Then there exists a constant $C_0>0$ such that for all $x\in\R^d$, $(n_i,E_i):\R^d\to \mathcal U_{a,\alpha}$, $i=0,1$, being analytic in $x$ with $\NN{(n_i,E_i)}_{\dot C^\nu_x}\leq C$ and 
	\[\NN{n_i}_{\dot C^\nu_x}+\NN{E_i}_{\dot C^\nu_x}\leq \frac{\delta}2 n_i(x)(1-\eta n_i(x)) \quad\mbox{for } i=0,1,\]
	we have
	
	\begin{equation*}
	\NN{\FrhoE(n_1,E_1)-\FrhoE(n_0,E_0)}_{ C^\nu_x}\leq C_1 \left(\NN{n_1-n_0}_{ C^\nu_x}+\NN{E_1-E_0}_{ C^\nu_x}\right).
	\end{equation*}
\end{lemma}
\begin{proof}
	Throughout this proof, we fix $x\in \R^d$ and make sure that the constants do not depend explicitly on $x$. To start with, we assume w.l.o.g.~that
	\[n_0(x)(1-\eta n_0(x))\geq n_1(x)(1-\eta n_1(x)).\]
	We define 
	
	\begin{equation*}
	(n_\theta,E_\theta):=\begin{cases}
	(n_0,(1-2\theta)E_0+2\theta E_1)), &\mbox{for }\theta\in\left[0,\frac12\right]\\
	((2-2\theta)n_0+(2\theta-1) n_1),E_1)), &\mbox{for }\theta\in\left[\frac12,1\right].
	\end{cases}
	\end{equation*}
	Hence,  for $\theta\in[0,\frac12]$
	
	\begin{align*}
	\N{n_\theta}_{\dot C^\nu_x}+\N{E_\theta}_{\dot C^\nu_x}
	&\leq  \N{n_0}_{\dot C^\nu_x}+\N{E_0}_{\dot C^\nu_x}+ \N{E_1}_{\dot C^\nu_x}
	\\&\leq \delta n_0(x)(1-\eta n_0(x))=\delta n_\theta(x)(1-\eta n_\theta(x)).
	\end{align*}
	For $\theta\in[\frac12,1]$ it holds
	
	\begin{align*}
	\N{n_\theta}_{\dot C^\nu_x}
	&\leq  (2-2\theta)\N{n_0}_{\dot C^\nu_x}+(2\theta-1)\N{n_1}_{\dot C^\nu_x}
	\\&\leq (2-2\theta) Cn_0(x)(1-\eta n_0(x))+ (2\theta-1)Cn_1(x)(1-\eta n_1(x)).
	\end{align*}
	Since the mapping $t\mapsto t(1-\eta t)$ is concave, we have
	
	\begin{align*}
	\N{n_\theta}_{\dot C^\nu_x}\leq \frac{\delta}{2} n_\theta(x)(1-\eta n_\theta(x)).
	\end{align*}
	Moreover,  we have $ n_1(x)(1-\eta n_1(x))\leq  n_\theta(x)(1-\eta n_\theta(x))$ by construction, which implies
	
	\begin{align*}
	\N{E_\theta}_{\dot C^\nu_x}=\N{E_1}_{\dot C^\nu_x}\leq \frac{\delta}{2}  n_1(x)(1-\eta n_1(x))\leq \frac{\delta}{2}  n_\theta(x)(1-\eta n_\theta(x)).
	\end{align*}
	We summarize that		
	
	\begin{align*}
	\N{n_\theta}_{\dot C^\nu_x}+\N{E_\theta}_{\dot C^\nu_x}
	&\leq\delta n_\theta(x)(1-\eta n_\theta(x))
	\end{align*}
	for all $\theta\in[0,1]$. The formula
	
	\begin{multline*}
	\FrhoE(n_1,E_1)-\FrhoE(n_0,E_0)\\=(E_1-E_0)\int_0^\frac12 \del_{E}\FrhoE(n_\theta,E_\theta)d\theta 
	+(n_1-n_0)\int_\frac12^1 \del_{n}\FrhoE(n_\theta,E_\theta)d\theta 
	\end{multline*}
	and the properties of $(n_\theta,E_\theta)$ we can use the same techniques as in the proof of Lemma \ref{lem.NN.Frho.nE} in order to finish the proof. 
\end{proof}
\begin{proof}[Proof of Proposition \ref{prop.Ff.is.Lipschitz.nochmal}]
	Finally, the proposition is a direct consequence of Lemmas \ref{lem.NN.Frho.nE} and \ref{lem.FrhoE.Lipschitz}, because  we have $\NN{n_f}_{\dot C^\nu_x}+\NN{E_f}_{\dot C^\nu_x}\leq (1+\NN{\epsilon}_{L^1(\T^d)})\NN{f}_{\dot C^\nu_x}$ as well as $\NN{n_f}_{C^\nu_x}+\NN{E_f}_{C^\nu_x}\leq (1+\NN{\epsilon}_{L^1(\T^d)})\NN{f}_{C^\nu_x}$.
\end{proof}

\medskip
% The data information below will be filled by AIMS editorial staff
%Received xxxx 20xx; revised xxxx 20xx.
\medskip


\begin{thebibliography}{99}
	
	\bibitem{AsMe77} [10.1063/1.3037370]
	N.~W. Ashcroft and N.~D. Mermin.
	\newblock {Solid State Physics},
	\newblock \emph{Physics Today}, 30, 1977,
	%	\newblock {doi: 10.1063/1.3037370}.	
	
	\bibitem{ADH15} A. Al-Masoudi, S. D\"orscher, S. H\"afner, U. Sterr, and C. Lisdat.
	Noise and instability of an optical lattice clock. {\em Phys. Rev. A} 92 (2015),
	063814, 7 pages.
	
	\bibitem{BeDe96} N.~B. Abdallah and P. Degond. On a hierarchy of macroscopic 
	models for semiconductors. {\em J. Math. Phys.} 37 (1996), 3308-3333.
	
	\bibitem{Blo05} E.~Bloch. Ultracold quantum gases in optical lattices.
	{\em Nature Physics} 1 (2005), 23-30.
	
	
	\bibitem{Bra17} 
	\newblock M. Braukhoff, 
	\newblock  \emph{ Effective Equations for a Cloud of Ultracold Atoms in an Optical Lattice},
	\newblock  Ph.D thesis, University of Cologne, Germany 2017.
	
	\bibitem{BrJu17} [10.1142/S021820251850015X]
	M. Braukhoff and A. J\"ungel,  Energy-transport systems for optical lattices: derivation, analysis, simulation, \emph{Mathematical Models and Methods in Applied Sciences} (2017), in press. 
	
	\bibitem{BaBe13} [10.3934/krm.2013.6.893]
	C. {Bardos} and N. {Besse}.
	\newblock {The Cauchy problem for the Vlasov-Dirac-Benney equation and related
		issues in fluid mechanics and semi-classical limits.}
	\newblock \emph{{Kinet. Relat. Models}}, 6\penalty0 (4):\penalty0 893--917,
	2013.
	%	\newblock ISSN 1937-5093; 1937-5077/e.
	%	\newblock {doi: 10.3934/krm.2013.6.893}.
	
	\bibitem{BaBe15} [10.1007/978-1-4939-2950-4]
	C. {Bardos} and N. {Besse}.
	\newblock {Hamiltonian structure, fluid representation and stability for the
		Vlasov-Dirac-benney equation.}
	\newblock In \emph{{Hamiltonian partial differential equations and
			applications. Selected papers based on the presentations at the conference on
			Hamiltonian PDEs: analysis, computations and applications, Toronto, Canada,
			January 10--12, 2014}}, pages 1--30. Toronto: The Fields Institute for
	Research in the Mathematical Sciences; New York, NY: Springer, 2015.
	%	\newblock ISBN 978-1-4939-2949-8/hbk; 978-1-4939-2950-4/ebook.
	%	\newblock {doi: 10.1007/978-1-4939-2950-4 1}
	
	\bibitem{BaBe16}
	C. {Bardos} and N. {Besse}.
	\newblock {Semi-classical limit of an infinite dimensional system of nonlinear
		Schr\"odinger equations.}
	\newblock \emph{{Bull. Inst. Math., Acad. Sin. (N.S.)}}, 11\penalty0
	(1):\penalty0 43--61, 2016.
	%	\newblock ISSN 2304-7909; 2304-7895/e.
	
	
	\bibitem{BaNo12} [10.1063/1.4765338]
	C. {Bardos} and A. {Nouri}.
	\newblock {A Vlasov equation with Dirac potential used in fusion plasmas.}
	\newblock \emph{{J. Math. Phys.}}, 53\penalty0 (11):\penalty0 115621, 16, 2012.
	%	\newblock ISSN 0022-2488; 1089-7658/e.
	%	\newblock {doi: 10.1063/1.4765338}.
	
	
	
	\bibitem{DGH15} O.~Dutta, M.~Gajda, P.~Hauke, M.~Lewenstein, D.-S.~L\"uhmann, 
	B.~Malomed, T.~Sowinski, and J.~Zakrzewski.  Non-standard Hubbard models in 
	optical lattices: a review. {\em Rep. Prog. Phys.} 78 (2015), 066001, 47 pages.
	
	%\bibitem{DrJu12} M.~Dreher and A.~J\"ungel. Compact families of piecewise constant 
	%functions in $L^p(0,T;B)$. {\em Nonlin. Anal.} 75 (2012), 3072-3077. 
	
	\bibitem{GNZ09} A.~Griffin, T.~Nikuni, and E.~Zaremba. {\em Bose-Condensed	
		Gases at Finite Temperatures}. Cambridge University Press, Cambridge, 2009.
	
	\bibitem{HaNg15} [10.1007/s00205-016-0985-z]
	D.~{Han-Kwan} and T.~T.~{Nguyen}.
	\newblock Ill-posedness of the hydrostatic {E}uler and singular {V}lasov
	equations.
	\newblock \emph{Arch. Rational Mech. Anal.}, 221(3):1317-1344, 2016.
	
	\bibitem{HaRo16} [10.24033/asens.2313]
	D.~{Han-Kwan} and F. {Rousset}.
	\newblock 
	Quasineutral limit for Vlasov-Poisson with Penrose stable data.
	\newblock \emph{Ann. Sci. École Norm. Sup.}, 49(6):1445-1495, 2016.
	
	\bibitem{JaNo11} [10.1016/j.crma.2011.03.024]
	P.-E. {Jabin} and A.~{Nouri}.
	\newblock {Analytic solutions to a strongly nonlinear Vlasov equation.}
	\newblock \emph{{C. R., Math., Acad. Sci. Paris}}, 349\penalty0
	(9-10):\penalty0 541--546, 2011.
	%	\newblock ISSN 1631-073X.
	%	\newblock {doi: 10.1016/j.crma.2011.03.024}.
	
	\bibitem{Jak04} A.~Jaksch. Optical lattices, ultracold atoms and quantum 
	information processing. {\em Contemp. Phys.} 45 (2004), 367-381.
	
	
	\bibitem{Jue09} [10.1007/978-3-540-89526-8]
	\newblock A.~J\"ungel.  Transport Equations for Semiconductors.
	\emph{Lect. Notes Phys}. 773. Springer, Berlin, 2009.
	
	\bibitem{MoVi11} [10.1007/s11511-011-0068-9]
	C. {Mouhot} and C. {Villani}.
	\newblock {On Landau damping.}
	\newblock \emph{{Acta Math.}}, 207\penalty0 (1):\penalty0 29--201, 2011.
	%	\newblock ISSN 0001-5962; 1871-2509/e.
	%	\newblock {do10.1007/s11511-011-0068-9}.
	
	\bibitem{Ram56} N.~Ramsey. Thermodynamics and statistical mechanics at negative
	absolute temperature. {\em Phys. Rev.} 103 (1956), 20-28.
	
	\bibitem{RMR10} A.~Rapp, S.~Mandt, and A.~Rosch. Equilibration rates and negative
	absolute temperatures for ultracold atoms in optical lattices.
	{\em Phys. Rev. Lett.} 105 (2010), 220405, 4 pages.
	
	\bibitem{SHR12} U.~Schneider, L.~Hackerm\"uller, J.~Ph.~Ronzheimer, S.~Will,
	S.~Braun, T.~Best, I.~Bloch, E.~Demler, S.~Mandt, D.~Rasch, and A.~Rosch.
	Fermionic transport and out-of-equilibrium dynamics in a homogeneous Hubbard
	model with ultracold atoms. {\em Nature Physics} 8 (2012), 213-218.
	
	
\end{thebibliography}
\end{document}